\documentclass[letterpaper, 11pt,  reqno]{amsart}

\usepackage{amsmath,amssymb,amscd,amsthm,amsxtra, esint}
\usepackage{tikz} 
\usepackage{color}
\usepackage{hyperref}

\usepackage{cases}

\usepackage{todonotes}

\usepackage[left=32mm, right=32mm, 
bottom=27mm]{geometry}

\usepackage{tikz} 
\usepackage{relsize}

\usepackage{tikz-cd}

\usepackage{marginnote}
\usepackage{scalerel} 

\usepackage{titlesec}

\titleformat{\section}{\centering\normalfont\bfseries\Large\sffamily}{\thesection.}{1em}{}
\titlespacing*{\section}{1em}{1em}{1em}
\titleformat{\subsection}[runin]{\normalfont\bfseries\large\sffamily}{\thesubsection.}{1em}{}[.]
\titlespacing*{\subsection}{1em}{1em}{1em}

\usetikzlibrary{shapes.misc}
\usetikzlibrary{shapes.symbols}
\usetikzlibrary{decorations}
\usetikzlibrary{decorations.markings}

\newcommand{\pl}{\mathbin{\scaleobj{0.7}{\tikz \draw (0,0) node[shape=circle,draw,inner sep=0pt,minimum size=8.5pt] {\scriptsize $<$};}}}

\newcommand{\ple}{\mathbin{\scaleobj{0.7}{\tikz \draw (0,0) node[shape=circle,draw,inner sep=0pt,minimum size=8.5pt] {\scriptsize $\leqslant$};}}}

\newcommand{\llpvariable}[2]{
    \tikz[baseline=0.01em]{\draw[line width=#1pt] (0, 0) -- (0, 0.6em); \draw[fill=#2, line width=#1pt] (0.1pt, 0.6em) circle (0.3ex);}
}

\newcommand{\llp}{
  \ensuremath{\llpvariable{0.8}{black}}
}

\DeclareRobustCommand{\gobblefour}[4]{}

\allowdisplaybreaks[2]

\usepackage{color}

\definecolor{gr}{rgb}   {0.,   0.69,   0.23 }
\definecolor{bl}{rgb}   {0.,   0.5,   1. }
\definecolor{mg}{rgb}   {0.85,  0.,    0.85}
\definecolor{yl}{rgb}   {0.8,  0.7,   0.}
\definecolor{or}{rgb}  {0.7,0.2,0.2}

\newcommand{\Cc}{\mathcal{C}}

\newcommand{\z}{\mathbf{z}}
\newcommand{\y}{\mathbf{y}}

\newcommand{\CC}{\mathcal{C}}
\newcommand{\CE}{\mathcal{E}}

\usepackage{dsfont}
\newcommand{\IDC}{\mathds{1}}

\newcommand{\bz}{\mathbf{z}}

\newtheoremstyle{mystyle}
{3pt}               
{3pt}               
{\it }                      
{}                      
{\sffamily\bfseries}             
{}                      
{0.5em}                 
{#1 #2{\hspace{0.2cm}--\hspace{-0.2cm}}  }   

\theoremstyle{mystyle}

\newtheorem{thm}{Theorem}
\newtheorem*{thm*}{Theorem}

\newtheorem{cor}[thm]{\hspace{-0.15cm}  {Corollary} }
\newtheorem{lem}[thm]{\hspace{-0.14cm}  {Lemma} }
\newtheorem{prop}[thm]{\hspace{-0.13cm} {Proposition}}
\newtheorem{rem}[thm]{\hspace{-0.15cm} {Remark}}

\DeclareMathOperator*{\Argmin}{Argmin}

\newcommand{\R}{\mathbb{R}}

\newcommand{\T}{\mathbb{T}}

\let\P= \undefined
\newcommand{\P}{\mathbf{P}}

\newcommand{\E}{\mathbb{E}}

\newcommand{\eps}{\varepsilon}

\newcommand{\A}{\mathcal{A}}
\newcommand{\les}{\lesssim}

\newcommand{\Prob}{\mathbb{P}}

\newcommand{\N}{\mathbb{N}}

\renewcommand{\H}{\mathcal{H}}
\newcommand{\D}{\mathcal{D}}

\makeatletter
\def\DeclareSymbol#1#2#3{\expandafter\gdef\csname MH@symb@#1\endcsname{\tikz[baseline=#2, scale=.18]{#3}}}
\def\<#1>{\ensuremath{\mathchoice{\tikzsetnextfilename{macros#1}{\color{black}\csname MH@symb@#1\endcsname}}{\tikzsetnextfilename{macros#1}{\color{black}\csname MH@symb@#1\endcsname}}{\tikzsetnextfilename{macros#1}\scalebox{.7}{\color{black}\csname MH@symb@#1\endcsname}}
{\tikzsetnextfilename{macros#1}\scalebox{.5}{\color{black}\csname MH@symb@#1\endcsname}}}} 
\makeatother 
\DeclareSymbol{1}{0}{\path (-.4,0) -- (.4,0); \draw[line width=0.2mm] (0,0) -- (0,1.2) 
node[circle, fill, draw, solid, inner sep=0pt, minimum size=1.8pt] {};}

\newcommand{\pv}{\hspace{0.05cm};}
\newcommand{\Wick}[1]{\mathbf{:}#1\mathbf{:}}
\newcommand{\de}{\mathrm d}

\usepackage{mathtools}
\numberwithin{equation}{section}
\mathtoolsset{showonlyrefs}

\renewenvironment{proof}{\noindent{\sffamily{\textbf{Proof :}}}}{\begin{flushright}$\square$\end{flushright}}
\newenvironment{proofof}{\noindent{\sffamily{\textbf{Proof }}}}{\begin{flushright}$\square$\end{flushright}}

\usepackage{stmaryrd}

\numberwithin{equation}{section}
\numberwithin{thm}{section}

\begin{document}

\baselineskip = 15pt

\title{Ergodicity of the Anderson $\Phi^4_2$ model}

\author[H.~Eulry and A. Mouzard]
{Hugo Eulry and Antoine Mouzard}

\begin{abstract}
We consider the parabolic stochastic quantization equation associated to the $\Phi^4_2$ model on the torus in a spatial white noise environment. We study the long time behavior of this heat equation with independent multiplicative white noise and additive spacetime white noise, which is a singular SPDE in a singular environement and requires two different renormalization procedures. We prove that the solution is global in time with a strong a priori $L^p$ bound independent of the initial data in $\CC^{-\eps}$ for large $p$. The quenched solution given the environment is shown to be an infinite dimensional Markov process which satisfies the strong Feller property. We prove exponential convergence to a unique invariant measure using a Doeblin criterion for the transition semigroup. 

In particular, our work is a generalization of a previous work by Tsatsoulis and Weber in a case which is not translation invariant hence the method makes no use of the reversibility of the dynamics or the explicit knowledge of the invariant measure and it is therefore in principle applicable to situations where these are not available such as the vector-valued case.
\end{abstract}


\maketitle


\baselineskip = 14pt

\section{Introduction}

Consider the parabolic $\Phi^4_2$ equation
\begin{equation}{\label{AndersonSQE}}
    \partial_t u + \H u + \Wick{u^3} = \sqrt{2}\zeta
\end{equation}
on the torus $\T^2$ driven by the Anderson Hamiltonian 
\begin{equation}
\H=-\Delta+\xi
\end{equation}
where $\xi$ and $\zeta$ are two independent white noises, respectively space and spacetime. This equation falls in the range of singular SPDEs and involves two renormalization procedures associated to each layer of randomness, respectively to define the Anderson Hamiltonian $\H$ and the time evolution of the nonlinearity. In a previous work \cite{EMR24}, we obtained pathwhise local well-posedness of the equation for any initial data $u_0\in\CC^{-\eps}$ with $\eps>0$ small enough as well as probabilistic global well-posedness for initial datum $u_0$ distributed according to the invariant Gibbs measure proportional to
\begin{equation}{\label{AndersonGibbs}}
e^{-\frac{1}{4}\int_{\T^2}\Wick{u^4}}\de\mu^\H(u)
\end{equation}
where $\mu^\H$ is the law of the Anderson Gaussian Free Field associated to $\H$ and $\Wick{u^4}$ the Wick power. In this work, we continue the study of the long time behavior of the dynamic as well as the probabilistic properties of the underlying infinite dimensional stochastic process. We work on the quenched process $u(t)$ in the infinite dimensional space $\CC^{-\eps}$ and prove that there exists a global solution that converges exponentially fast to a unique invariant measure for any initial data $u_0\in\CC^{-\eps}$. We do not use at any point the a priori knowledge of the invariant measure constructed in \cite{EMR24}.

\medskip

The $\Phi_d^4$ model has a long history as a toy model for Quantum Field Theory that we do not intend to explain here, see for example Glimm and Jaffe's book \cite{GJ81} or the recent construction based on SPDEs of $\Phi_3^4$ on the full space by Gubinelli and Hofmanová \cite{GH21} and references therein. The measure was constructed for general polynomial nonlinearity with the so-called $P(\Phi)_2$ model while the construction of the $\Phi_3^4$ measure is already quite difficult, usually in the flat case of $\R^d$ and $\T^d$ which is invariant by translation. The construction of such measures is one of the motivation to the investigation of nonlinear heat equations with additive spacetime white noise following Parisi and Wu \cite{PW81} with stochastic quantization. The idea is to consider a dynamic on an infinite dimensional space and construct the Gibbs measure $\mathcal{Z}^{-1}e^{-\CE(u)}\de u$ as the invariant probability measure of the equation
\begin{equation}
\partial_tu=-\nabla\CE(u)+\sqrt{2}\de W_t
\end{equation}
with $(W_t)_{t\ge0}$ an infinite dimensional Brownian motion. In the case of $\Phi_d^4$ on the torus $\T^d$, the energy is given by
\begin{equation}
\CE(u)=\frac{\lambda}{4}\int_{\T^d}u(x)^4\de x-\frac{m}{2}\int_{\T^d}u(x)^2\de x+\frac{1}{2}\int_{\T^d}|\nabla u(x)|^2\de x
\end{equation}
with coupling constant $\lambda$ and mass $m$ which yields the stochastic quantization equation
\begin{equation}
\partial_tu=\Delta u+mu-\lambda u^3+\sqrt{2}\zeta
\end{equation}
with $\zeta$ a spacetime white noise of parabolic spacetime Hölder regularity $-\frac{2+d}{2}-\kappa$ for any $\kappa>0$, we fix $\lambda=m=1$ in the following discussion. In dimension $d=1$, this equation is well-posed since $u$ is of Hölder regularity $\frac{1}{2}-\kappa$ thanks to Schauder estimates hence the cubic term is well-defined. In dimension $d=2$, the solution $u$ should be a distribution of Hölder regularity $-\kappa<0$ hence the cubic term is ill-defined and requires a renormalization procedure. The first strong solutions to this equation were obtained by Da Prato and Debussche \cite{DD03} with the change of variable
\begin{equation}
u_t=\llp(t)+v_t
\end{equation}
with $\llp(t)=\int_{-\infty}^te^{(t-s)\Delta}\zeta(\de s)$ the stationnary solution to the linear equation. Then $v$ satisfies
\begin{equation}
\partial_tv=\Delta v+v-v^3-3v^2\llp-3v\llp^2-\llp^3
\end{equation}
where the powers $\llp^2$ and $\llp^3$ are still ill-defined. Since $\llp(t)$ is a Gaussian process, its powers can be renormalized as Wick products with $\Wick{\llp^2}$ and $\Wick{\llp^3}$ of Hölder regularity $-\kappa$ for any $\kappa>0$, as irregular as $\llp$ itself. This cancels the roughest term in the equation which makes $v$ is of positive regularity and each term in the equation is well-defined. They prove local well-posedness for the equation in $\CC^{-\eps}$ as well as invariance of the Gibbs measure which can be directly constructed with probabilistic global well-posedness. This is also the content of our previous work \cite{EMR24} for the model considered in this work with a general polynomial nonlinearity. In the more recent work \cite{TW18}, Tsatsoulis and Weber consider the general $\Phi_2^p$ equation and prove global well-posedness in $\CC^{-\eps}$ using a coming down from infinity bound uniform in the initial data. They prove that the solution is a Markov process in $\CC^{-\eps}$ satisfying the strong Feller property. They finally prove a support theorem for the laws of the solutions in the cubic case which implies exponential convergence to a unique invariant measure with Doeblin criterion. In particular, this achieves a construction of the $\Phi_2^4$ measure using only the dynamic. In three dimensions, $\llp$ is of Hölder regularity $-\frac{1}{2}-\kappa$ hence the Wick products $\Wick{\llp^2}$ and $\Wick{\llp^3}$ are of respective Hölder regularity $-1-\kappa$ and $-\frac{3}{2}-\kappa$. This implies that the solution $v$ is a priori of regularity $\frac{1}{2}-\kappa$ making the product $v^2\llp$ and $v\Wick{\llp^2}$ singular. This equation thus falls in the range of singular SPDE. Catellier and Chouk \cite{CC18} proved local well-posedness using paracontrolled calculus following Gubinelli, Imkeller and Perkowski \cite{GIP15}. Mourat and Weber \cite{MW17} proved a similar a priori bound in this setting which implies global well-posedness in $\CC^{-\frac{1}{2}-\eps}$. We finally mention the recent work by Jagannath and Perkowski \cite{JP21} which relies on an exponential transform to prove global well-posedness for the equation on $v$, see also Bailleul, Dang, Ferdinand, and Tô \cite{BDFT23} on compact manifolds in three dimensions, in particular not invariant by translation.

\medskip

In this work, we are interested in the case where the Laplacian is replaced by the Anderson Hamiltonian $\H$ with the random environment $\xi$ a space white noise, which can appear as the scaling limit of Ising model with an i.i.d. environment. The operator itself is a singular stochastic operator that has to be defined with a renormalization procedure as
\begin{equation}
\H:=\lim_{\eps\to0}(-\Delta+\xi_\eps-c_\eps)
\end{equation}
with $\xi_\eps$ a regularization of the noise and $c_\eps$ a logarithmic diverging quantity. The first result of this work is an improvement from probabilistic global well-posedness to deterministic global well-posedness in negative regularity Hölder spaces $\CC^{-\eps}$ for $\eps>0$. This is done using a priori estimates on the $L^p$ norm of the solution for $p$ large enough such that $L^p$ is embbeded in $\CC^{-\eps}$ as in the usual coming down from infinity argument. The difficulty here is similar to $\Phi_3^4$ since the solution $u$ is paracontrolled which makes the a priori estimate hard to get. One first follows Da Prato and Debussche \cite{DD03} and consider
\begin{equation}
u_t=\llp(t)+v_t
\end{equation}
with $\llp(t)=\int_{-\infty}^te^{-(t-s)\H}\zeta(\de s)$ the stochastic convolution associated to the Anderson Hamiltonian. The equation on $v$ is singular and we use the paracontrolled expansion. Following the construction of the Anderson Hamiltonian from \cite{M21}, the solution $v$ is of the form
\begin{equation}
v=v\pl X+w
\end{equation}
with $X=\Delta^{-1}\xi$ and $w$ a smoother remainder. Up to a truncation of the low frequencies of $X$, the map $\Phi(v)=v-v\pl X$ which associates to $v$ its remainder $w$ is invertible which allows to consider the change of variable $v=\Gamma w$ with $\Gamma=\Phi^{-1}$. This map has been an important tool in recent works on evolution PDEs associated to the Anderson Hamiltonian, see \cite{EMR24,GUZ20} and references therein. At the level of the equation $w$, we obtain a priori estimates on its $L^p$ norm which allows to control the initial solution $v=\Gamma w$ and prove global well-posedness, see Corollary \ref{COR:Globalv}. This is close in spirit to the coming down from infinity argument for $\Phi_3^4$ which is described by a singular SPDEs, see Mourrat and Weber \cite{MW17}. This is enough to obtain existence of invariant measures in Proposition \ref{PROP:ExistenceInvariant} with Krylov-Bogoliubov Theorem, this is the approach in \cite{BDFT23} on compact surfaces in three dimensions.

We then study probabilistic properties of the solution as a stochastic process in infinite dimension. The main object of interest is the transition semigroup $(P_t)_{t\ge0}$ defined by
\begin{equation}{\label{MarkovSemigroup}}
P_t\varphi(u_0):=\E \big[\varphi(u(t\pv u_0))\big]
\end{equation}
where $\varphi$ is a bounded Borelian function on $\CC^{-\varepsilon}$ and $u(t\pv u_0)$ denotes the solution of \ref{AndersonSQE} at time $t$ starting from initial data $u_0\in\CC^{-\varepsilon}$. We prove that this is a Markov semigroup satisfying the strong Feller property, that is $P_t\varphi$ is a continuous function of $u_0$ for $t>0$ and $\varphi$ a bounded Borelian function. This is done with a Bismut-Elworthy-Li formula for the semigroup to obtain a uniform bound on $\de(P_t\varphi)$ which implies Lipchitz continuity for $P_t\varphi$. Finally, we prove exponential mixing to a unique invariant measure in Corollary \ref{Cor:ExpMixing} which is the main result of our work. We do not have a strong support theorem for the enhanced noise as obtained in \cite{TW18} due to the roughness of the environment however we prove that it contains zero. This is enough to prove that the solution relaxes to any small ball around zero with positive probability which implies that $0$ is in the support of any invariant measure hence uniqueness with the strong Feller property. In particular, this complete the case of usual $\Phi_2^p$ model for general even integer $p\ge4$ in the work of Tsatsoulis and Weber \cite{TW18} as they only restrict to the cubic model for the support theorem. Our method allows to consider models not invariant by translation such as similar model on any compact surfaces, we restrict to the torus for the sake of simplicity as the invariance is already lost with the Anderson Hamiltonian.

\medskip

In Section \ref{Sec:Renormalization}, we recall the renormalization result from our previous work \cite{EMR24}. We also improve the time regularity of the Wick power and prove that zero is in the support of the enhanced noise. In Section \ref{Sec:globalWP}, we first prove a general local well-posedness result in the same spirit as \cite{EMR24}. We then prove an a priori estimate in $L^p$ for $p$ large enough on the solution independent of $u_0\in\CC^{-\eps}$ which implies global well-posedness. The Markov and strong Feller properties are proved in Section \ref{Sec:strongFeller} via a Bismut-Elworthy-Li formula. In Section \ref{Sec:ergodicity}, we prove that the solution relaxes to $0$ with positive probability which implies exponential convergence to a unique invariant measure using a Doeblin criterion for the transition semigroup.

\medskip

{\bf Acknoledgments :} The authors warmly thank Arnaud Debussche for fruitful and motivating dicussions about the present work. The first author was supported by the ANR projet Smooth “ANR-22-CE40-0017”.

\section{Renormalization and noise regularity}\label{Sec:Renormalization}

In this section, we gather all the needed results on the renormalization of the enhanced noise that will be needed. The stochastic convolution
\begin{equation}
\llp_{\infty,t}=\sqrt{2}\int_{-\infty}^te^{-(t-s)\H}\zeta(\de s)
\end{equation}
with $-\H$ instead of the Laplacian is a stationnary Gaussian process with explicit spacetime covariance given by the Anderson heat kernel. Its law at any fixed $t\in\R$ is given by the Anderson Gaussian free field $\mu^\H$ thus of Hölder spatial regularity $\CC^{-\eps}$ for any $\eps>0$, see Section $5$ in \cite{BDM25}. The renormalization of its Wick power follows from estimates on the functionnal calculus associated to $\H$, we refer to Section $2$ in \cite{EMR24} for details. The Wick power are defined as
$$
\Wick{\llp_{-\infty,t}^n}=\lim_{N\to\infty}\big(H_n(\Pi_N\llp_{-\infty,t},\sigma_N)\big)=\lim_{N\to\infty}\Wick{(\Pi_N\llp_{-\infty,t})^n}
$$
for any $n\ge2$ with $(H_n)_{n\ge0}$ the Hermite polynomials with variance $\sigma_N=\E\big[(\Pi_N\llp_{-\infty,t})^2\big]$ and
$$
\Pi_Nu:=\sum_{k=0}^{N}\left\langle u,\varphi_k^\H\right\rangle\varphi_k^\H
$$ 
the spectral projection on the first $N$ modes of $\H$ for $N\ge0$. The following result corresponds to Proposition $2.4$ in \cite{EMR24}.

\medskip

\begin{prop}\label{Prop:LargeDev}
Let $T>0$, $\varepsilon\in(0,1)$ and $p\geq q\geq1$. For any $n\ge1$, the sequence $(\Wick{(\Pi_N\llp_{-\infty,t})^n})_{N\ge0}$ is a Cauchy sequence in $L^p(\Omega,L^q([0,T],\CC^{-\varepsilon}))$ which converges almost surely to a limit $\Wick{\llp_{\infty,t}^n}$ in $L^q([0,T],\CC^{-\varepsilon})$. Moreover, the map $t\mapsto\llp(t)$ belongs almost surely to $C([0,T],\CC^{-\varepsilon})$ and we have the following tail estimates and deviation bounds:
\begin{align}\label{EstimateTailLoliN}
\mathbb{P}\left(\|\Wick{(\Pi_N\llp_{-\infty,t})^n}\|_{L^q_T\CC_x^{-\varepsilon}}> R\right)\leq C e^{-cR^\frac{2}{n}T^{-\frac{2}{q}}},
\end{align}
and
\begin{align}\label{EstimateDeviationLoliN}
\mathbb{P}\left(\|\Wick{(\Pi_N\llp_{-\infty,t})^n}-\Wick{\llp^n}\|_{L^q_T\CC_x^{-\varepsilon}}> N^{-\kappa}R\right)\leq C e^{-cR^{\frac{2}{n}}T^{-\frac{2}{qn}}},
\end{align}
for some constant $c,C>0$ and exponent $\kappa>0$ that do not depend on $T$, $R$, $p$ and $N$.
\end{prop}

\medskip

To study the probabilistic properties of the solution, we will also need the more general process
$$
\llp_{s,t}=\sqrt{2}\int_{s}^te^{-(t-r)\H}\zeta(\de r)
$$
for any $s\le t$. Since we have
$$
\llp_{s,t}=\llp_{-\infty,t}-e^{-(t-s)\H}\llp_{-\infty,s},
$$
this allows to consider
$$
\Wick{\llp_{s,t}^n}=\sum_{k=0}^n\binom{n}{k}(-1)^k\left(e^{-(t-s)\H}\llp_{-\infty,s}\right)^k\Wick{\llp_{-\infty,t}^{n-k}}
$$
using that $e^{-(t-s)\H}$ is smoothing for $s<t$ for $n\ge2$ as well as algebraic properties of the Hermite polynomials. This only makes sense a priori for the regularized quantities $\Wick{(\Pi_N\llp_{s,t})^n}$ however the formula holds for the limit as proved in the next proposition. For the global well-posedness, we will work with
\begin{equation}
\llp(t)=\sqrt{2}\int_0^te^{-(t-s)\H}\zeta(\de s)
\end{equation}
for $t\ge0$, that is $\llp(t)=\llp_{0,t}$. This gives for $\llp_{s,\cdot}$ similar estimates as for $\llp_{\infty,\cdot}$ and its renormalized powers. Since we are interested in the cubic equation, we consider
\begin{equation}
\z_N:=\Big(\z_N^{(1)},\z_N^{(2)},\z_N^{(3)}\Big)
\end{equation}
the enhanced data with $\z_N^{(1)}=3\Pi_N\llp$ the truncation of $\llp$ with the Wick products
\begin{equation}
\z_N^{(2)}=3\Wick{\big(\Pi_N\llp\big)^2}=3(\Pi_N\llp)^2-3\E\big[\Pi_N\llp^2\big]
\end{equation}
and
\begin{equation}
\z_N^{(3)}=\Wick{\big(\Pi_N\llp\big)^3}=(\Pi_N\llp)^3-3\E\big[\Pi_N\llp^2\big]\Pi_N\llp
\end{equation}
for $N\ge0$ as well as its limit $\z=(\z^{(1)},\z^{(2)},\z^{(3)})$. The following result corresponds to Proposition $2.4$ in \cite{EMR24}, it ensures in particular
\begin{equation}
\z\in C([0,T],\CC^{-\eps})\times L^q([0,T],\CC^{-\eps})\times L^q([0,T],\CC^{-\eps})
\end{equation}
for any $T>0$ with $\eps>0$ and $q\ge1$.

\medskip

\begin{lem}\label{Lem:Wickllp}
Let $\eps>0$ and $q\ge1$ such that $\eps q<1$. For $n\ge1$, we have
\begin{equation}
\lim_{N\to\infty}\E\Big[\|\Wick{\llp_{t,t+\cdot}^n}-\Wick{(\Pi_N\llp_{t,t+\cdot})^n}\|_{L^q([0,T],\CC^{-\eps})}^p\Big]=0
\end{equation}
for any $t\in\R$ and $p\ge1$. In particular, we have
\begin{equation}
\lim_{N\to\infty}\E\Big[\|\z-\z_N\|_{L^q([0,T],\CC^{-\eps})^3}^p\Big]=0
\end{equation}
for any $p\ge1$. Moreover, we have the binomial expansion
\begin{equation}
\Wick{\llp_{s,t}^n}=\sum_{k=0}^n\binom{n}{k}(-1)^k\left(e^{-(t-s)\H}\llp_{-\infty,s}\right)^k\Wick{\llp_{-\infty,t}^{n-k}}
\end{equation}
for any $s<t$.
\end{lem}

\medskip

\begin{proof}
This follows directly from the algebraic properties of the Hermite polynomials as well as the regularization properties of $e^{-t\H}$, see section $2$ in \cite{TW18} for the usual $\Phi_2^n$ model. Since we have
$$
\Wick{(\Pi_N\llp_{s,t})^n}=\sum_{k=0}^n\binom{n}{k}(-1)^k\left(e^{-(t-s)\H}\Pi_N\llp_{-\infty,s}\right)^k\Wick{(\Pi_N\llp_{-\infty,t})^{n-k}}
$$
for any $n\ge1$, the result follows from
\begin{equation}
\Big\|\Wick{(\Pi_N\llp_{s,t})^n}-\Wick{(\Pi_N\llp_{s,t})^n}\Big\|_{\CC^{-\eps}}\lesssim\big(1+\big\|\big(e^{-(t-s)\H}\llp_{-\infty,s}\big)^n\big\|\big)\sum_{k=0}^n\Big\|\Wick{(\Pi_N\llp_{\infty,t})^k}-\Wick{(\Pi_N\llp_{\infty,t})^k}\Big\|_{\CC^{-\eps}}
\end{equation}
using that $\CC^\eps$ is an algebra and the Schauder estimates
\begin{equation}
\big\|e^{-(t-s)\H}\llp_{-\infty,s}\big\|_{\CC^\eps}\lesssim(t-s)^{-\eps}\|\llp_{-\infty,s}\big\|_{\CC^{-\eps}}
\end{equation}
which blows up as $t-s$ goes to zero but is integrable for $\eps q<1$.
\end{proof}

It will be important to consider $\z_{s,t}$ the enhanced noise on the time interval $[s,t]$ for any $s<t$ defined as $\z$ with $\llp_{s,t}$ instead of $\llp$. The following proposition guarantees that the noise is invariant by translation in law. It follows directly from the same property for the regularized data $\z_N$ with the convergence in $L^p(\Omega)$.

\medskip

\begin{lem}\label{Lem:llp}
For any $t\in\R$, the process $\bz_{t,t+\cdot}$ and $\bz$ have the same law.
\end{lem}

\medskip

For the model with the usual Laplacian instead of the Anderson Hamitlonian $\H$, one can prove that the Wick powers converge in $C([0,T],\CC^{-\eps})$ hence continuity in time follows directly. This is not known a priori in our case, we prove the continuity in time directly using Kolmogorov's criterion. The explosion at time zero is due to the initial condition $\llp(0)=0$.

\medskip

\begin{prop}\label{Prop:Regularitynoise}
The enhanced noise $\z=(\z^{(1)},\z^{(2)},\z^{(3)})$ belongs to $C\big((0,T],(\CC^{-\kappa})^3\big)$ for any $T>0$ and $\kappa>0$.
\end{prop}

\medskip

\begin{proof}
We first prove the continuity in time of the stationnary process
\begin{equation}
\llp_{\infty,t}=\sqrt{2}\int_{-\infty}^te^{-(t-s)\H}\zeta(\de s)
\end{equation}
which is a centered Gaussian process with covariance
\begin{align}
\E\big[\llp_{\infty,t}(x)\llp_{\infty,s}(y)\big]&=2\int_{-\infty}^t\int_{-\infty}^s\int_{\T^2\times\T^2}K_{t-t'}(x,x')K_{s-s'}(y,y')\E\big[\zeta(t',x')\zeta(s',y')\big]\de x'\de y'\de t'\de s'\\ 
&=2\int_{-\infty}^{t\wedge s}\int_{\T^2}K_{t-u}(x,z)K_{s-u}(y,z)\de z\de u\\ 
&=2\int_{-\infty}^{t\wedge s}K_{t+s-2u}(x,y)\de u
\end{align}
where $K_t(x,y)$ denotes the kernel of the heat semigroup $e^{-t\H}$. We denote as $\Delta_n$ the Paley-Littlewood projector with associated kernel $\chi_n=2^{2n}\chi(2^n\cdot)$, see Appendix \ref{APP:Besov}. For any random distribution $f$, we have
\begin{align}
\E\big[\|f\|_{\CC^{\alpha}}^p\big]&\lesssim\E\big[\|f\|_{B_{p,p}^{\alpha+\frac{2}{p}}}^p\big]\\ 
&\lesssim\sum_{n\ge0}2^{\alpha np}\E\big[\|\Delta_nf\|_{L^p}^p\big]\\ 
&\lesssim\sum_{n\ge0}2^{\alpha np}\int_{\T^2}\E\big[|\Delta_nf(x)|^p\big]\de x
\end{align}
hence a pointwise bound on $\E\big[|\Delta_nf(x)|^p\big]$ is enough to get Hölder regularity. We have
\begin{align}
\E\big[|\Delta_n\llp_{\infty,t}(x)|^2\big]&=\int_{\T^2\times\T^2}\chi_n(x-y)\chi_n(x-y')\E\big[\llp_{\infty,t}(y)\llp_{\infty,t}(y')\big]\de y\de y'\\
&=2\int_{-\infty}^t\int_{\T^2\times\T^2}\chi_n(x-y)\chi_n(x-y')K_{2t-2u}(y,y')\de y\de y'\de u\\
&=\int_{\T^2\times\T^2}\chi_n(x-y)\chi_n(x-y')G(y,y')\de y\de y'\\
&\lesssim\int_{\T^2\times\T^2}\chi_n(x-y)\chi_n(x-y')\log(|y-y'|)\de y\de y'\\
&\lesssim\log(n)
\end{align}
using the support of $\chi_n$ where $G$ is the Green function of $\H$. For $p\ge2$, Gaussian hypercontractivity yields
\begin{equation}
\E\big[|\Delta_n\llp_{\infty,t}(x)|^p\big]\lesssim\log(n)^{\frac{p}{2}}
\end{equation} 
which gives the Hölder spatial regularity $\CC^{-\kappa}$ for $\llp_{\infty,t}$ with fixed $t\in\R$ and any $\kappa>0$. For the time regularity, we have
\begin{align}
\E\big[&|\Delta_n\llp_{\infty,t}(x)-\Delta_n\llp_{\infty,s}(x)|^2\big]\\ 
&=\int_{\T^2\times\T^2}\chi_n(x-y)\chi_n(x-y')\E\big[(\llp_{\infty,t}(y)-\llp_{\infty,s}(y))(\llp_{\infty,t}(y')-\llp_{\infty,s}(y'))\big]\de y\de y'\\
&=4\int_{\T^2\times\T^2}\chi_n(x-y)\chi_n(x-y')\Big(\int_{-\infty}^{t\wedge s}K_{2(t\wedge s)-2u}(y,y')\de u-\int_{-\infty}^{t\wedge s}K_{t+s-2u}(y,y')\de u\Big)\de y\de y'.
\end{align}
For $t=s+\delta$ with $\delta>0$, we get
\begin{equation}
\E\big[|\Delta_n\llp_{\infty,t}(x)-\Delta_n\llp_{\infty,s}(x)|^2\big]=4\int_{\T^2\times\T^2}\chi_n(x-y)\chi_n(x-y')\Big(\int_0^\delta K_u(y,y')\de u\Big)\de y\de y'
\end{equation}
which we can bound with
\begin{equation}
K_u(y,y')\lesssim\frac{1}{u}e^{-c\frac{|y-y'|^2}{u}}
\end{equation}
with $c>0$ from Proposition \ref{Prop:GaussianBound}. We get
\begin{equation}
\E\big[|\Delta_n\llp_{\infty,t}(x)-\Delta_n\llp_{\infty,s}(x)|^2\big]\lesssim|t-s|^{\alpha_\kappa} 2^{\kappa n}
\end{equation}
for $\alpha_\kappa>0$ which diverges as $\kappa$ goes to $0$. Kolmogorov's criterion following the previous computation with $f=\llp_{\infty,t}-\llp_{\infty,s}$ thus gives that $\llp$ belongs to $\CC([0,T],\CC^{-\kappa})$ with 
\begin{equation}
\E\big[\sup_{t\in[0,T]}\|\llp_{\infty,t}\|_{\CC^{-\kappa}}^p\big]<\infty
\end{equation}
for any $p\ge1$. For the Wick powers, this follows from the same computations using the covariance
\begin{equation}
\E\Big[\Wick{\llp_{\infty,t}(x)^n}\Wick{\llp_{\infty,s}(y)^n}\Big]=n!\E\Big[\llp_{\infty,t}(x)\llp_{\infty,s}(y)\Big]^n
\end{equation}
with Gaussian hypercontractivity in Gaussian multiplicative chaos. In the end, this proves that $\llp$ and its Wick powers belongs to $C(\R,\CC^{-\kappa})$ for any $\kappa>0$. To conclud for the enhanced data $\z$, we use the binomial expression
\begin{equation}
\Wick{\llp(t)^n}=\sum_{k=0}^n\binom{n}{k}(-1)^k\left(e^{-t\H}\llp_{-\infty,0}\right)^k\Wick{\llp_{-\infty,t}^{n-k}}
\end{equation}
with Schauder estimates which gives an explosion as $t$ goes to zero as for the proof of Lemma \ref{Lem:Wickllp}.
\end{proof}

The continuity in time is not important for global well-posedness, it is however needed to prove the next proposition. In \cite{TW18}, Tsatsoulis and Weber proved that the enhanced noise for the $\Phi_2^4$ model on the torus $\T^2$ contains $(0,a,0)$ in its support for any $a\ge0$. Their proof relies strongly on the explicit spectral theory of the Laplacian, it is not clear how to prove such a result in general including for example the case of a general compact surface or more general power. In the last section, we prove that this weak support result is enough to get uniqueness of the invariant measure with exponential ergodicity.

\medskip

\begin{prop}\label{Prop:Supportnoise}
For any $T>t_0>0$, we have
\begin{equation}
\mathbb{P}\Big(\sup_{t\in[t_0,T]}\|\z(t)\|_{(\CC^{-\kappa})^3}\le\eps\Big)>0
\end{equation}
for any $\eps>0$. In particular, $(0,0,0)$ is in the support of $\z$ as a random variable in $C([t_0,T],(\CC^{-\kappa})^3)$.
\end{prop}

\medskip

\begin{proof}
Our strategy of proof is to consider the real random variable
\begin{equation}
Z=\sup_{t\in[t_0,T]}\big(\|\z_t^{(1)}\|_{\CC^{-\kappa}}+\|\z_t^{(2)}\|_{\CC^{-\kappa}}+\|\z_t^{(3)}\|_{\CC^{-\kappa}}\big)
\end{equation}
and prove that
\begin{equation}
\mathbb{P}(Z<\eps)>0
\end{equation}
for any $\eps>0$. This implies that $0$ lies in the support of $Z$ and gives the result. As done by Stroock and Varadhan \cite{SV72}, it is enough to prove that
\begin{equation}
\mathbb{P}(Z<\eps|\A)>0
\end{equation}
for any event $\A$ of positive probability. The conditional Markov inequality
\begin{equation}
\mathbb{P}(Z<\eps|\A)\ge1-\frac{\E[Z|\A]}{\eps}
\end{equation}
gives that the result holds if we construct an event $\A=\A_\eps$ such that $\E[Z|\A]<\eps$. We consider the representation of $\z$ associated to
\begin{equation}
\llp_{\infty,t}=\sum_{n\ge0}X_n(t)e_n
\end{equation}
with $(X_n)_{n\ge0}$ a family of independent stationnary Ornstein-Uhlenbeck process of parameters $(\lambda_n)_{n\ge0}$ the eigenvalues of $\H$. Given any $N\ge0$, the event
\begin{equation}
\A_N:=\Big\{\forall i\in\llbracket0,N\rrbracket,\sup_{t\in[0,T]}|X_i(t)|\le\frac{\eps}{(N+1)\sup_{0\le n\le N}\|e_n\|_{\CC^{\kappa}}}\Big\}
\end{equation}
is indeed of positive probability. Then we consider the decomposition
\begin{equation}
\llp_{\infty,t}=\sum_{n=0}^NX_n(t)e_n+\sum_{n>N}X_n(t)e_n=X_{\le N}(t)+X_{>N}(t)
\end{equation}
which gives for the Wick powers
\begin{equation}
\Wick{\llp_{\infty,t}^2}=X_{\le N}(t)+2X_{\le N}(t)X_{>N}(t)+\Wick{X_{>N}(t)^2}
\end{equation}
and
\begin{equation}
\Wick{\llp_{\infty,t}^3}=X_{\le N}(t)^3+3X_{\le N}(t)^2X_{>N}(t)+3X_{\le N}(t)\Wick{X_{>N}(t)^2}+\Wick{X_{>N}(t)^3}
\end{equation}
for any $t\in[0,T]$. We get
\begin{align}
\|\llp_{\infty,t}\|_{\CC^{-\kappa}}&\le\|X_{\le N}(t)\|_{\CC^{\kappa}}+\|X_{>N}(t)\|_{\CC^{-\kappa}}\\ 
&\le\eps+\|X_{>N}(t)\|_{\CC^{-\kappa}}
\end{align}
as well as
\begin{align} 
\|\Wick{\llp_{\infty,t}^2}\|_{\CC^{-\kappa}}&\le\eps^2+2\eps\|X_{>N}(t)\|_{\CC^{-\kappa}}+\|\Wick{X_{>N}(t)^2}\|_{\CC^{-\kappa}}
\end{align}
and
\begin{equation}
\|\Wick{\llp_{\infty,t}^3}\|_{\CC^{-\kappa}}\le\eps^3+3\eps^2\|X_{>N}(t)\|_{\CC^{-\kappa}}+3\eps\|\Wick{X_{>N}(t)^2}\|_{\CC^{-\kappa}}+\|\Wick{X_{>N}(t)^3}\|_{\CC^{-\kappa}}
\end{equation}
for the Wick powers. Since the random variables $\Wick{X_{>N}(t)^k}$ depends only on $X_i$ for $i>N$, we have
\begin{equation}
\E\Big[\sup_{t\in[0,T]}\|\Wick{X_{>N}(t)^k}\|_{\CC^{-\kappa}}\big|\A_N\Big]=\E\Big[\sup_{t\in[0,T]}\|\Wick{X_{>N}(t)^k}\|_{\CC^{-\kappa}}\Big]
\end{equation}
for any $k\in\{1,2,3\}$. Following the proof of the previous proposition, we have 
\begin{equation}
\lim_{N\to\infty}\E\Big[\sup_{t\in[0,T]}\|\Wick{X_{>N}(t)^k}\|_{\CC^{-\kappa}}\Big]=0
\end{equation}
using the quantitative estimate from Kolmogorov's continuity criterion, see for example Theorem $3.4.16$ from \cite{S24}. Taking $N=N_\eps>0$ large enough depending on $\eps$, we get that
\begin{equation}
\E[Z|\A_N]\le 2\eps+3\eps^2+\eps+4\eps^3+3\eps^2+\eps
\end{equation}
hence 
\begin{equation}
\mathbb{P}\Big(\sup_{t\in[0,T]}\big(\|\llp_{\infty,t}\|_{\CC^{-\kappa}}+\|\Wick{\llp_{\infty,t}^2}\|_{\CC^{-\kappa}}+\|\Wick{\llp_{\infty,t}^3}\|_{\CC^{-\kappa}}\big)\le\eps\Big)>0
\end{equation}
for any $\eps>0$. Again, the result for the enhanced data $\z$ follows from the binomial expression
\begin{equation}
\Wick{\llp(t)^n}=\sum_{k=0}^n\binom{n}{k}(-1)^k\left(e^{-t\H}\llp_{-\infty,0}\right)^k\Wick{\llp_{-\infty,t}^{n-k}}
\end{equation}
with Schauder estimates which gives an explosion as $t$ goes to zero.
\end{proof}

\begin{rem}
In \cite{CF18,TW18}, the authors make strong use of the spectral theory of the Laplacian on the torus to construct a sequence $h_n$ in the Cameron-Martin space that contains nonzero elements. This allows to give a strong support theorem for the enhanced noise thus for solutions of singular SPDEs. While we believe that this approach can be adapted for any potential of positive Hölder regularity, this seems out of reach in the case of the Anderson Hamiltonian or even for the Laplacian of general compact surfaces. However, in our previous work \cite{EMR24}, we constructed the Gibbs measure $\nu^\H$ associated to the Anderson $P(\Phi_2)$ model as an absolutetly continuous measure with respect to the Anderson Gaussian free field $\mu^\H$. Since the Gaussian measure $\mu^\H$ is known to be of full support in $\CC^{-\kappa}$ as its Cameron-Martin space is dense, this is also the case for $\nu^\H$. We do not use this argument in the present work, which does not rely on a direct construction of the measure.
\end{rem}

\section{General well-posedness theory}\label{Sec:globalWP}

In this section, we study the truncated dynamic
\begin{equation}{\label{TruncEquation}}
\partial_t u_N + \H u_N + \Wick{u_N^3}=\sqrt{2}\Pi_N\zeta,\\
\end{equation}
with initial data $u_N(0)=u_0\in\Cc^{-\eps}$ and follow the usual Da Prato-Debussche trick from \cite{DD03} in the context of the usual $\Phi_2^4$ equation. We split the equation into a linear dynamics driven by $\zeta$ and a non-linear dynamics starting from zero initial data. To that end, we write 
\begin{equation}
u_N = \Pi_N\llp + v_N
\end{equation} 
where $\llp$ is the solution to the linear equation
\begin{align*}
    (\partial_t + \H) \llp= \sqrt{2}\zeta
\end{align*}
with $\llp(0)=0$. Let $\z_N=\left(\z_N^{(1)},\z_N^{(2)},\z_N^{(3)}\right)$ be the enhanced data with $\z_N^{(1)}:=3\Pi_N\llp$ the truncation of $\llp$ with the Wick products
\begin{equation}
\z_N^{(2)}:=3\Wick{\big(\Pi_N\llp\big)^2}=3(\Pi_N\llp)^2-3\E\big[\Pi_N\llp^2\big]
\end{equation}
and
\begin{equation}
\z_N^{(3)}:=\Wick{\big(\Pi_N\llp\big)^3}=(\Pi_N\llp)^3-3\Pi_N\llp\E\big[\Pi_N\llp^2\big].
\end{equation}
Then $v_N$ solves the nonlinear equation
\begin{equation}{\label{NonLinearEqv_N}}
\partial_t v_N + \H v_N + v_N^3 + v_N^2\z_N^{(1)}+ v_N\z_N^{(2)} + \z_N^{(3)} =0
\end{equation}
with the same initial data $v_N(0)=u_0\in\CC^{-\eps}$. Provided we can solve \eqref{NonLinearEqv_N} in an appropriate functional space, solutions to \eqref{TruncEquation} will therefore be constructed as $u_N=\Pi_N\llp+v_N$. Given the Anderson Hamiltonian $\H$, the equation is singular due to the irregularity of the spacetime white noise $\zeta$ and the cubic nonlinearity. To construct solutions, we consider here a regularization of the stochastic forcing while a regularization of the nonlinearity was used in our previous work \cite{EMR24}. This choice was motivated there by our goal to prove invariance of the constructed Gibbs measure and obtain probabilistic global well-posedness while it is more convenient to work with a truncated spacetime white noise to obtain the strong Feller property. We start by proving local well-posedness within this framework with the same type or argument as in \cite{EMR24}. We then prove global well-posedness for any initial data in negative Hölder spaces $\CC^{-\eps}$ using a priori estimates.

\subsection{Local theory}{\label{SEC:local}}

We investigate the local well-posedness of \eqref{TruncEquation} for initial data $u_0\in\CC^{-\varepsilon}$ with the change of variable
$$
u_N=\Pi_N\llp+v_N
$$
where $v_N$ satisfies equation \eqref{NonLinearEqv_N}. In view of the Schauder estimate for $\H$ from Proposition \ref{Prop:SchauderH}, the solution $v_N$ is expected to belong to the space
\begin{align}\label{space}
X_T^{-\eps,\sigma}=C([0,T],\CC^{-\eps})\cap C((0,T],\CC^{\sigma})
\end{align}
for $\eps<\sigma<1$ endowed with the norm
$$
\|v\|_{X_T^{-\eps,\sigma}}=\|v\|_{L^\infty_T\CC^{-\varepsilon}}+\sup_{0<t\leq T}t^{\frac{\sigma+\varepsilon}{2}}\|v(t)\|_{\CC^{\sigma}}
$$
that takes account of the blow-up at $t=0$ in $\CC^\sigma$. The unusual restriction $\sigma<1$ instead of $\sigma<2$ in the case of the standard Laplace operator is due to the roughness of the Anderson Hamiltonian. One can prove that the solution belongs to the $L^2$-based Sobolev space $\D^{2-\kappa}$ associated to $\H$ which only embeds in $\CC^{1-\kappa}$ with $\kappa>0$, this structure will be important to obtain a priori estimates. Note that even if \eqref{NonLinearEqv_N} is a truncated equation, by setting
$$
\mathbf{f}(v,\y):=v^3 + v^2\y^{(1)}+ v\y^{(2)} + \y^{(3)},
$$
it falls into the scope of the following more general problem 
\begin{align}\label{model}
    \begin{cases}
    \partial_tv+\H v + \mathbf{f}(v,\y)=0,\\
    v(0)=v_0,
    \end{cases}
\end{align}
for a given vector $\y\in L^q([0;T];\CC^{-\eps})^3$ and initial condition $v_0\in\CC^{-\eps}$. We then have the following local well-posedness result. In particular, the convergence of $v_N$ to a limit solution $v$ follows from the convergence of $\z_N$ to $\z$ since the solution depends continuously on the data.

\medskip

\begin{prop}{\label{Localwpv}}
Let $0<\eps<\sigma< 1$ and $ 1\le q<\infty$ be such that $\frac{\sigma+\eps}{2}\frac{q}{q-1}<\frac1{3}$. Then there exists $\theta,C>0$ such that for any $r,R>0$, there exists $T>0$ of order $\frac{r^\theta}{(1+r+R)^{3\theta}}$ such that for any $v_0\in\CC^{-\eps}$ and $\y\in L^q([0;T],\CC^{-\eps})^3$ satisfying $\|v_0\|_{\Cc^{-\eps}}\le r$ and $\|\y\|_{L^q_{T}\Cc^{-\eps}}\le R$, there exists a unique solution $v$ in $X_T^{-\eps,\sigma}$ to \eqref{model}.
Moreover, we have
\begin{align}\label{localbound}
\|v\|_{X_T^{-\eps,\sigma}}\le 2Cr.
\end{align}
and $v$ depends Lipschitz continuously on both $v_0$ and $\y$.
\end{prop}

\medskip

\begin{proof}
Given the data $(v_0,\y)$, we define the solution mapping $\Phi$ by
$$
\Phi(v)(t):=e^{-t\H}v_0-\int_0^te^{-(t-s)\H}\mathbf{f}(v,\y)(s)\de s
$$
and seek for a solution $v$ as a fixed point of $\Phi$ in the ball of radius $2Cr$ in $X_T^{-\eps,\sigma}$ for some $C>0$ and $T>0$ small enough depending on $r,R$ as in the statement of Proposition~\ref{Localwpv}. Let $v$ be in the ball of radius $2Cr$ in $X_T^{-\eps,\sigma}$. For the norm of negative regularity, we have
\begin{align*}
\big\|\Phi(v)&\big\|_{L^\infty_T\Cc^{-\eps}}\les \|v_0\|_{\Cc^{-\eps}}+\int_0^T\big\|\mathbf{f}(v,\y)(s)\big\|_{\Cc^{-\eps}}\de s\\
&\les \|v_0\|_{\Cc^{-\eps}}+\sum_{j=1}^3\int_0^T\|v(s)\|_{\Cc^{\sigma}}^j\|\y^{(3-j)}(s)\|_{\Cc^{-\eps}}\de s+\int_0^T\|\y^{(3)}(s)\|_{\Cc^{-\eps}}\de s\\
&\les \|v_0\|_{\Cc^{-\eps}}+\sum_{j=1}^3\|v\|_{X_T^{-\eps,\sigma}}^j\|\y^{(3-j)}\|_{L^q_T\Cc^{-\eps}}\Big(\int_0^Ts^{-\big(\frac{\sigma+\eps}2\big)jq'}\de s\Big)^{\frac1{q'}}+T^{\frac1{q'}}\|\y^{(3)}\|_{L^q_T\Cc^{-\eps}}\\
&\le Cr+C'(1+r+R)^{3}T^{\frac1{q'}-3\frac{\sigma+\eps}{2}}
\end{align*}
with $C,C'>0$ independant of $r$, $R$, and $T$ using the expression of $\mathbf{f}$, Schauder estimate for $\H$ and a repeated use of the product estimates in Lemma \ref{LEM:Besov}. Thus, taking $T$ of order $\big(\frac{r}{(1+r+R)^{3}}\big)^{\theta}$ with $\theta=\big(\frac1{q'}-3\frac{\sigma+\eps}{2}\big)^{-1}>0$, we get
\begin{align*}
\big\|\Phi(v)\big\|_{L^\infty_T\Cc^{-\eps}}\le 2Cr.
\end{align*}
Similarly as above, we have
\begin{align*}
\sup_{0<t\le T}t^{\frac{\sigma+\eps}{2}}&\big\|\Phi(v)(t)\big\|_{\Cc^{\sigma}}\les \|v_0\|_{\Cc^{-\eps}}+\sup_{0<t\le T}t^{\frac{\sigma+\eps}{2}}\int_0^t(t-s)^{-\frac{\sigma+\eps}{2}}\|\y^{(3)}(s)\|_{\Cc^{-\eps}}\de s\\
&\quad+\sum_{j=1}^3\sup_{0<t\le T}t^{\frac{\sigma+\eps}{2}}\int_0^t (t-s)^{-\frac{\sigma+\eps}2}\|v^j(s)\y^{(3-j)}(s)\|_{\Cc^{-\eps}}\de s\\
&\les \|v_0\|_{\Cc^{-\eps}}+\Big\{\|\y^{(3)}\|_{L^q_T\Cc^{-\eps}}\sup_{0<t\le T}t^{\frac{\sigma+\eps}{2}}\Big(\int_0^t(t-s)^{-\frac{\sigma+\eps}{2}q'}\de s\Big)^{\frac1{q'}}\\
&\quad+\sum_{j=1}^3\|v\|_{X_T^{-\eps,\sigma}}^j\|\y^{(3-j)}\|_{L^q_T\Cc^{-\eps}}\sup_{0<t\le T}t^{\frac{\sigma+\eps}{2}}\Big(\int_0^t(t-s)^{-\frac{\sigma+\eps}{2}q'}s^{-\big(\frac{\sigma+\eps}2\big)jq'}\de s\Big)^{\frac1{q'}}\\
&\le Cr+C'(1+r+R)^3T^{\frac1{q'}-3\frac{\sigma+\eps}{2}}
\end{align*}
where we used that
\begin{align*}
\int_0^t(t-s)^{-a}s^{-b}\de s\les t^{1-a-b}
\end{align*}   
for any $0<a,b<1$ with $a+b>1$.
Putting the two estimates together with our choice of $T$ yields
\begin{align*}
\big\|\Phi(v)\big\|_{X_T^{-\eps,\sigma}}\le 2Cr
\end{align*}
hence the ball of radius $2Cr$ in $X_T^{-\eps,\sigma}$ is stable under the map $\Phi$ up to the time continuity of $\Phi(v)$. For this point, consider $0<t\leq T$ and $h$ small enough then split the integral as
\begin{align*}
\Phi(v)(t+h)-\Phi(v)(t)&=e^{-t\H}(e^{-h\H}-1)v_0+\int_0^te^{-(t-s)\H}(e^{-h\H}-1)\mathbf{f}(v,\y)(s)\de s\\
&\qquad+\int_t^{t+h}e^{-(t+h-s)\H}\mathbf{f}(v,\y)(s)\de s\\
&=:(1) + (2) + (3).
\end{align*}
We use both the Schauder estimate \ref{Prop:SchauderH} and the continuity of the semigroup given by Proposition \ref{Prop:ContinuityH} to prove the convergence of each term as $h$ goes to $0$. The first term converges to $0$ by continuity of the semigroup in $\CC^\sigma$. For the second term, Schauder estimate for $\H$ gives
\begin{equation}
\left\|(2)\right\|_{\CC^{\sigma}}\lesssim \int_0^t(t-s)^{-\frac{\sigma+\varepsilon}{2}}\|(e^{-h\H}-1)\mathbf{f}(v,\y)(s)\|_{\CC^{-\varepsilon}}\de s
\end{equation}
which again converges to $0$ with the dominated convergence theorem and continuity of the semigroup. The third term is controlled a similar argument, changing variables and proceeding as above we get
\begin{align*}
\left\|(3)\right\|_{\CC^{\sigma}}&\lesssim \int_t^{t+h}(t+h-s)^{-\frac{\sigma+\varepsilon}{2}}\|\mathbf{f}(v,\y)(s)\|_{\CC^{-\varepsilon}}\de s\\
&\lesssim \int_0^{h}\tau^{-\frac{\sigma+\varepsilon}{2}}\|\mathbf{f}(v,\y)(t+h-\tau)\|_{\CC^{-\varepsilon}}\de\tau\\
&\lesssim \left(\int_0^{h}\tau^{-q'\frac{\sigma+\varepsilon}{2}}(t+h-s)^{-q'3\frac{\sigma+\varepsilon}{2}}\de s\right)^{\frac{1}{q'}}\|v\|_{X_T^{-\eps,\sigma}}^3\|\y\|_{L^q_{T_0}\CC^{-\varepsilon}}\\
&\lesssim t^{-3\frac{\sigma+\varepsilon}{2}}h^{\frac{1}{q'}-\frac{\sigma+\varepsilon}{2}}\|v\|_{X_T^{-\eps,\sigma}}^3\|\y\|_{L^q_{T_0}\CC^{-\varepsilon}}
\end{align*}
which converges to $0$ as $h$ goes to $0$. This proves continuity of $\Phi(v)$ as a $\CC^\sigma$-valued map on $(0,T]$. The continuity with values in $\CC^{-\varepsilon}$ at $0$ is obtained via the same arguments using that $u_0\in\CC^{-\eps}$. Altogether, this shows that $\Phi$ maps the ball of radius $2Cr$ of $X_T^{-\eps,\sigma}$ to itself. Similar computations as above yield
\begin{align*}
\big\|\Phi(v_1)-\Phi(v_2)\big\|_{L^\infty_T\Cc^{-\eps}}\leq C'(1+r+R)^2T^{\frac1{q'}-\frac{\sigma+\eps}{2}3}\|v_1-v_2\|_{X_T^{-\eps,\sigma}}
\end{align*}
and
\begin{align*}
\sup_{0<t\le T}t^{\frac{\sigma+\eps}{2}}\big\|\Phi(v_1)(t)-\Phi(v_2)(t)\big\|_{\Cc^{\sigma}}\le C'(1+r+R)^2T^{\frac1{q'}-\frac{\sigma+\eps}{2}3}\|v_1-v_2\|_{X_T^{-\eps,\sigma}}
\end{align*}
for $v_1,v_2$ in the ball of radius $2Cr$ of $X_T^{-\eps,\sigma}$. With our choice of $T$, this shows that $\Phi$ is also a contraction on this ball, thus admits a unique fixed point $v$ which is then a mild solution to \eqref{model} in $X_T^{-\eps,\sigma}$ in the ball of radius $2Cr$. Given different data $(v_0,\y),(\widetilde{v_0},\widetilde{\y})$, similar estimates for solutions $v,\widetilde{v}$ yield
\begin{align*}
\|v-\widetilde{v}\|_{X_T^{-\eps,\sigma}}&\le C\|v_0-\widetilde{v_0}\|_{\Cc^{-\eps}}+C'(1+r+R)^2T^{\frac1{q'}-3\frac{\sigma+\eps}{2}}\|v-\widetilde{v}\|_{X_T^{-\eps,\sigma}}\\
&\qquad+C'(1+r)^2T^{\frac1{q'}-3\frac{\sigma+\eps}{2}}\|\y-\widetilde{\y}\|_{L^q_T\Cc^{-\eps}}
\end{align*}
hence the Lipschitz dependence on $(v_0,\y)$ of the solution.
\end{proof}

The previous proposition holds for very general choice of forcing $\y\in L^q([0;T],\CC^{-\eps})^3$. If we return to the case of interest $\y=\z$, this ensures that equation \eqref{AndersonSQE} is locally well-posed in $\CC^{-\varepsilon}$ with solution in $X_T^{-\eps,\sigma}$. Moreover, the convergenc of the data $\y=\z_N$ to $\z$ in the correct topology ensures the convergence of $v_N$ to $v$ in $X_T^{-\eps,\sigma}$ as $N$ goes to infinity. It yields to following well-posedness result for the initial equation
\begin{equation}
\partial_t+\H u+\Wick{u^3}=\sqrt{2}\zeta
\end{equation}
in the space $\llp+X_T^{-\eps,\sigma}$. In particular, the renormalized product $\Wick{u^3}$ is defined as 
\begin{equation}
\Wick{u^3}=(u-\llp)^3+(u-\llp)^2\z^{(1)}+(u-\llp)\z^{(2)}+\z^{(3)}
\end{equation}
with the enhanced data $\z$ and $u-\llp$ regular enough for the different products to make sense.

\medskip

\begin{thm}\label{THM:LWP}
Let $0<\varepsilon<\frac1{3}$, $\eps<\sigma<1$, and $q>1$ be such that $\frac{\sigma+\eps}{2}\frac{q}{q-1}<\frac1{3}$. Then $\Prob$-almost surely, for any $u_0\in\Cc^{-\eps}$, there exists $T>0$ such that for any $N\in\N^*$,  \eqref{TruncEquation} admits a solution $u_N\in C([0;T];\Cc^{-\eps})$, unique in the affine space $\Pi_N\llp+X_T^{-\eps,\sigma}$. Moreover $u_N$ converges almost surely to a limit $u\in \llp+X_T^{-\eps,\sigma}$ which is the unique solution in the class $\llp+X_T^{-\eps,\sigma}$ to equation \eqref{AndersonSQE}.
\end{thm}

\medskip

\begin{proof}
Let $\eps,\sigma,q$ be as in the statement of Theorem~\ref{THM:LWP}, and recall that $\z$ is the enhanced data and $\z_N$ its approximation. Fix $T>0$ and let $R>0$. For any $N\ge1$ and $1\le k\le 3$, consider the event
\begin{equation}
\Sigma_R^{N,k}:=\Big\{\omega\in\Omega\ ;\ \|\z_N^{(k)}-\z^{(k)}\|_{L^q_{T}\CC^{-\varepsilon}}\leq N^{-\frac{\kappa}{2}}R,\ \|\z_N^{(k)}\|_{L^q_{T}\CC^{-\varepsilon}}\leq R \Big\}
\end{equation}
where $\kappa$ is as in Proposition \ref{Prop:LargeDev}. Then set
$$
\Sigma_R:=\bigcap_{N\in\N^*}\bigcap_{1\leq k\leq 3}\Sigma_R^{N,k}
$$
and
$$
\Sigma:=\liminf_{R\to+\infty}\Sigma_R.
$$
Using the bounds from Proposition \ref{Prop:LargeDev}, we get 
$$
\mathbb{P}\left(\Omega\setminus\Sigma_R\right) \leq C e^{-c\sqrt{R}}
$$
for some constants $c,C>0$ thus Borel-Cantelli Lemma gives that $\Sigma$ is of full probability. Given $R>0$, we get from Proposition~\ref{Localwpv} that there exists a time $T=T(\|u_0\|_{\Cc^{-\eps}},R)>0$ such that for any $\omega\in\Sigma_R$ and $N\geq1$, equation \eqref{model} admits a unique solution $v_N\in X_T^{-\eps,\sigma}$ with data $v_N(0)=u_0$. We then obtain a unique solution $u_N$ to \eqref{TruncEquation} by setting
$$
u_N=\Pi_N\llp + v_N\in \Pi_N\llp + X_T^{-\sigma,\eps}.
$$
Moreover since $\z_N^{(k)}\to \z^{(k)}$ in $L^q([0;T_0];\CC^{-\varepsilon})$, the Lipschitz continuity property of the solution in Proposition~\ref{Localwpv} ensures that $v_N$ converges in $X_T^{-\eps,\sigma}$ to the solution $v$ to \eqref{model} with data $(u_0,\z)$ hence $u_N$ converges to a limit
$$
u=\llp+v\in\llp+X_T^{-\eps,\sigma}
$$
which completes the proof.
\end{proof}

\subsection{A priori estimate}

We now turn to the globalization problem for \eqref{AndersonSQE}, which amounts to prove global existence for the solution $v$ to \eqref{model} for general noise data $\z=(\z^{(1)},\z^{(2)},\z^{(3)})\in C(\R,\CC^{-\varepsilon})\times L^q_{\text{loc}}(\R;\CC^{-\varepsilon})^2$ and initial condition $u_0\in\CC^{-\varepsilon}$ where $\varepsilon$ and $q$ are as in Proposition \ref{Localwpv} to ensure local well-posedness. This in particular includes the case of the fully renormalized equation \eqref{AndersonSQE} but the more general framework allows to play with a range of noise data and initial data, see subsection \ref{SEC:estimates}. In view of the blow-up criterion provided by Proposition \ref{Localwpv}, it is enough to prove that the $\CC^{-\varepsilon}$ norm of $v$ stays bounded uniformly in time. Using the embedding 
$$
L^{3p-2}\hookrightarrow \CC^{-\varepsilon}
$$ 
for $p>\frac{4}{3\varepsilon}$, it comes down to a bound the $L^{3p-2}$ norm of $v$ for $p$ large enough, this is the purpose of the rest of this section. Let $u_0\in\CC^{-\eps}$ and consider $T\in(0,+\infty]$ be the possible infinite blow up time in $\CC^{-\eps}$ of the solution.

We follow the argument from the usual $\Phi_d^4$ equation with a coming down from infinity estimate. Due to the singularity of the environment, we have to work with paracontrolled expansion for the solution hence we are closed to the $\Phi_3^4$ equation with Mourrat and Weber \cite{MW17a} than $\Phi_2^4$ with Tsatsoulis and Weber \cite{TW18}. Indeed, the idea in the usual $\Phi_2^4$ is to test the equation against $v^{3p-2}$ to obtain a priori bounds on the $L^{3p-2}$ norm which is not possible at the level of $u$ due to the irregularity of the linear term $\llp$. When $\H$ is replaced by the Laplacian, one can use the integration by part
\begin{equation}
\langle \Delta v,v^{3p-2}\rangle=-\big\langle|\nabla u|^2,u^{3p-3}\big\rangle
\end{equation}
which is not possible for the Anderson Hamiltonian. To overcome this issue, we use that $\H$ is a perturbation of the Laplacian which requires to work at the level of paracontrolled remainder $w=\Gamma^{-1}v$ implicitly defined by
\begin{equation}
v=v\pl X + w
\end{equation}
for $v$ in the Sobolev spaces associated to $\H$ where $X=\Delta^{-1}\xi\in\CC^{1-\kappa}$. In the following, we test the equation satified by $w$ against $w^{3p-2}$ a obtain a priori estimates. Using that $\Gamma$ is continuous from $L^p$ to itself, this gives the needed bound on $v$ and will give global well-posedness. In fact, it will be important to introduce a truncation parameter $n>0$ and rather consider the ansatz
$$
v=v\pl X_{>n} + w^{(n)}
$$
where $X_{>n}=\Delta_{>n}X$. Such truncation is already important in recent construction of the Anderson Hamiltonian, see \cite{BDM25,GUZ20,M21} and Appendix 
\ref{App:the_anderson_hamiltonian} for the full story. With this truncation, we are equiped with the one-to-one mapping $\Gamma_n:w^{(n)}\mapsto v$ defined as the inverse mapping of 
\begin{equation}
\Phi_n(v)=v-v\pl X_{>n}
\end{equation}
which is an invertible perturbation of the identity for $n=n(X)>0$ large enough. An important property is that $w^{(n)}$ inherits the same regularity properties as $w$, uniformly in $n$. The construction of $\H$ also relies heavily on the fact that one can tune the paraproduct $\pl$ so that $v$ is controlled by $w$ in appropriate norms while the converse being always true, see Lemma \ref{Lem:truncation}. We claim that it is still the case for $w^{(n)}$ and that the bound holds uniformly in $n$. In particular, we will drop the parameter $n$ in the notation $w^{(n)}$ in the following to lighten the computations, it will be fix large enough depending only on $X$ and universal constants.

\medskip

\begin{lem}{\label{Controlvwunif}}
For any $n\geq1$ and $p>1$, there is a constant $C>0$ that does not depend on $n$ such that
$$
C^{-1}\|v\|_{L^p}\leq \|w^{(n)}\|_{L^p} \leq C\|v\|_{L^p}
$$
with $v=\Gamma_nw$.
\end{lem}

\medskip

\begin{proof}
This follows from the continuity of $\Gamma_n$ and $\Phi_n$. The constant is uniform in $n$ since both applications converge to the identity as $n$ goes to infinity. Thus this is a direct consequence of Lemma \ref{Lem:truncation} and the fact that $\|X_{>n}\|_{B^{1-\kappa}_{q,\infty}}$ is bounded by $\|X\|_{B^{1-\kappa}_{q,\infty}}$ uniformly in $n$, for any $\kappa>0$ and $q\geq1$.
\end{proof}

This truncation will be important to make cubic terms involving $X_{>n}$ small for $n$ large enough. The construction of $\H$ with this paracontrolled expression gives a remainder operator such that
$$
\Gamma_n^{-1}\H\Gamma_nw= -\Delta w + \xi\ple w + R_n(w)
$$
and $R_n:B^{\kappa}_{p,\infty}\to B^{-\kappa'}_{p,\infty}$ is continuous. While the norm of $R_n$ as an operator diverges when $n$ is large, it will be harmless as $R_n$ is a linear term that will be made small using Young's estimate, no matter the truncation. Then equation \eqref{model} in the new variable $v=\Gamma_nw$ becomes 
\begin{align}{\label{Eqw}}
    \partial_t w -\Delta w + \xi\ple w + w^3 + Q(v,w,\bz) + R_n(w) = 0
\end{align}
where $Q(v,\bz) = Q_1(v,w)+Q_2(v,\bz)+Q_3(v,\bz)$ with
\begin{align*}
Q_1(v,w)&=-v^3\pl X_{>n} +(v\pl X_{>n})^3+ 3(v\pl X_{>n})w^2+3(v\pl X_{>n})^2w,\\
Q_2(v,\bz)&=v^2\bz^{(1)} + v\bz^{(2)} + \bz^{(3)},\\
Q_3(v,\bz)&=-(v^2\bz^{(1)} + v\bz^{(2)} + \bz^{(3)})\pl X_{>n}.
\end{align*}
We will keep $v$ in our notations to make things easier to read, but one should really keep in mind that $\Gamma_n$ is a one-to-one mapping between $w$ and $v$ so that we could drop the explicit dependency on $v$ and replace it by $\Gamma_n w$. It is crucial noting that the leading power $w^3$ comes with a positive sign. Note that while $Q_2$ and $Q_3$ are at most quadratic in $w$, $Q_1$ gathers all other cubic term in $w$, it will be controled using the parameter $n>0$ in the following. Assume first for the sake of the exposition that $w_0\in L^{3p-2}$. Testing equation \eqref{Eqw} against $w^{3p-3}$, we get
\begin{align*}
&\frac{1}{3p-2}\big(\|w_t\|_{L^{3p-2}}^{3p-2}-\|w_0\|_{L^{3p-2}}^{3p-2}\big) + (3p-3)\int_0^t\||\nabla w_s|^2w_s^{3p-4}\|_{L^1}\de s + \int_0^t\|w_s\|_{L^{3p}}^{3p}\de s\\
&\hspace{1cm} = -\int_0^t\langle \xi\ple w_s,w_s^{3p-3}\rangle\de s -\int_0^t\langle Q(v_s,w_s,\bz_s),w_s^{3p-3}\rangle\de s-\int_0^t\langle R_n(w_s),w_s^{3p-3}\rangle\de s
\end{align*}
for any positive time $t<T$ before the possible blow up. Since $w$ is the remainder in the paracontrolled expansion of $v$, it has enough regularity for $\nabla w$ to make sense hence $\||\nabla w_s|^2w_s^{3p-4}\|_{L^1}$ is a \textit{good} term since $p$ is an even integer, that is non-negative. Note that $Q$ gathers nonlinear contributions of $v$ where $Q_1$ has all the cubic dependencies on $v$ recalling that $w$ is linear in $v$, but only through the regularizing operator $\cdot\pl X_{>n}$. As for $Q_2$ and $Q_3$, they only involve lower order powers of $v$ while $R$ only depends linearly on $v$. The key idea is that the higher regularity norms of $w$ only appear integrated in time, we can then leverage the issue of estimating these terms by using the mild formulation of the equation 
\begin{align*}
w_t=e^{t\Delta}w_0 - \int_0^t e^{(t-s)\Delta}\big(\xi\ple w_s + w_s^3 + Q(v_s,w_s,\bz_s) + R_n(w_s)\big)\de s
\end{align*}
and the fact that $e^{(t-s)\Delta}$ is a regularizing operator, at the price of a diverging power of $t-s$. From there we will be able to close the estimate and obtain a \textit{coming down from infinity} estimate. We now prove estimates on all the needed term to get a first a priori estimate with Proposition \ref{LpBound}. Remember that $p$ is an even integer supposed to be large enough and $T\in(0,\infty]$ is the maximal time of existence provided by the local existence Theorem \ref{Localwpv}.

\medskip

\begin{lem}{\label{LEM:EstimateNoiseParaW}}
For any $\kappa>0$ and $0<s<t<T$, we have
$$
\left|\int_s^t\langle \xi\ple w_r,w_r^{3p-3}\rangle\de r\right|\lesssim \left(\int_s^t\|w_r\|_{L^{3p}}^{3p}\de r\right)^{\frac{p-1}{p}}\left(\int_s^t\|w_r\|^p_{B^{1+2\kappa}_{p,\infty}}\de r\right)^{\frac{1}{p}}
$$
where the implicit constant depends only on $p,\kappa$ and $\xi$.
\end{lem}

\medskip

\begin{proof}
Since the space white noise $\xi$ belongs to $\CC^{-1-\kappa}$ for any $\kappa>0$, we have
\begin{align*}
\left|\int_0^t\langle \xi\ple w_s,w_s^{3p-3}\rangle\de s\right|&\lesssim\int_0^t\|w_s^{3p-3}\|_{L^{\frac{3p}{3p-3}}}\|\xi\ple w_s\|_{L^p}\de s\\
&\lesssim\int_0^t\|w_s\|_{L^{3p}}^{3p-3}\|w_s\|_{B^{1+2\kappa}_{p,\infty}}\de s\\
&\lesssim \left(\int_0^t\|w_s\|_{L^{3p}}^{3p}\de s\right)^{\frac{p-1}{p}}\left(\int_0^t\|w_s\|^p_{B^{1+2\kappa}_{p,\infty}}\de s\right)^{\frac{1}{p}}
\end{align*}
using paraproduct and resonant term estimates as well as Hölder inequality and the relation
$$
1=\frac{p-1}{p}+\frac{1}{p}=\frac{3p-3}{3p}+\frac{1}{p}
$$
for $p\ge1$ which will be used multiple times in the following proofs.
\end{proof}

\begin{lem}{\label{LEM:EstimateR}}
For any $0<s<t<T$, we have
$$
\left|\int_s^t\langle R_n(w_r),w_r^{3p-3}\rangle\de r\right|\lesssim_n \left(\int_s^t\|w_r\|_{L^{3p}}^{3p}\de r\right)^{\frac{p-1}{p}}\left(\int_s^t\|w_r\|^p_{B^{\varepsilon}_{p,\infty}}\de r\right)^{\frac{1}{p}}.
$$
\end{lem}

\medskip

\begin{proof}
Lemma \ref{Lem:RemainderH} gives that the remainder term is a bounded operator
$$
R_n:B_{q,\infty}^{\kappa}\to B_{q,\infty}^{-\kappa'}
$$
for arbitrary $\kappa,\kappa'>0$ and $q\in[1,\infty]$. This yields
\begin{align*}
\big|\langle R_n(w),w^{3p-3}\rangle\big|&\lesssim \|R_n(w)\|_{B^{-\kappa}_{\frac{3p}{2},1}}\|w^{3p-3}\|_{B^{\kappa}_{\frac{3p}{3p-2},\infty}}\\
&\lesssim \|R_n(w)\|_{B^{-\kappa}_{\frac{3p}{2},\infty}}\|w^{3p-3}\|_{B^{\kappa}_{\frac{3p}{3p-2},\infty}}\\
&\lesssim_n \|w\|_{B^{\kappa}_{\frac{3p}{2},\infty}}\|w^{3p-3}\|_{B^{\kappa}_{\frac{3p}{3p-2},\infty}}
\end{align*}
and we have
\begin{align*}
\|w^{3p-3}\|_{B^{\kappa}_{\frac{3p}{3p-2},\infty}}&\lesssim \|w^{3p-4}\|_{L^{\frac{3p}{3p-4}}}\|w\|_{B^\kappa_{\frac{3p}{2},\infty}}\\
&\lesssim \|w\|_{L^{3p}}^{3p-4}\|w\|_{B^\kappa_{\frac{3p}{2},\infty}}
\end{align*}
using the nonlinear estimate from Lemma \ref{LEM:Besov}.
We interpolate $B^{\kappa}_{\frac{3p}{2},\infty}$ between $L^{\frac{3p}{2}}$ and $B^{2\kappa}_{\frac{3p}{2},\infty}$ using Lemma \ref{LEM:InterpolationBesov}.
\begin{align*}
\|w\|_{B^{\kappa}_{\frac{3p}{2},\infty}}&\lesssim \|w\|_{B^{2\kappa}_{\frac{3p}{2},\infty}}^{\frac{1}{2}}\|w\|_{L^{\frac{3p}{2}}}^{\frac{1}{2}}.
\end{align*} 
The estimate follows by Hölder inequality with
\begin{align*}
\big|\int_s^t\langle R_n(w_r),w_r^{3p-3}\rangle\de r\big|&\lesssim_n \int_s^t \|w_r\|_{L^{3p}}^{3p-3}\|w_r\|_{B^{1+\varepsilon}_{p,\infty}}\de r\\
&\lesssim_n \left(\int_s^t\|w_r\|^{3p}_{L^{3p}}\de r\right)^{\frac{p-1}{p}}\left(\int_s^t\|w_r\|^p_{B^{2\kappa}_{p,\infty}}\de r\right)^{\frac{1}{p}}
\end{align*}
and setting $2\kappa=\varepsilon$.
\end{proof}

\begin{lem}{\label{LEM:EstimateQ1}}
For any $\delta>0$, there exists $n_0=n_0(X,\delta,p)\ge0$ such that for any $n\geq n_0$ and $0<s<t<T$, we have
$$
\left|\int_s^t\langle Q_1(v_r,w_r),w_r^{3p-3}\rangle\de r\right|\leq \delta \int_s^t\|w_r\|_{L^{3p}}^{3p}\de r.
$$
\end{lem}

\medskip

\begin{proof}
Fix $\delta>0$ and recall that
\begin{equation}
Q_1(v,w)=-v^3\pl X_{>n} +(v\pl X_{>n})^3+ 3(v\pl X_{>n})w^2+3(v\pl X_{>n})^2w.
\end{equation}
We have
\begin{align*}
\left|\int_s^t\langle (v_r\pl X_{>n})w_r^2,w_r^{3p-3}\rangle\de r\right|&\lesssim\int_s^t\|v_r\pl X_{>n}\|_{L^{3p}}\|w_r^{3p-1}\|_{L^{\frac{3p}{3p-1}}}\de r\\ 
&\lesssim\|\cdot\pl X_{>n}\|_{L^{3p}\to L^{3p}}\int_s^t\|v_r\|_{L^{3p}}\|w_r\|_{L^{3p}}^{3p-1}\de r\\
&\lesssim\|X_{>n}\|_{\CC^\kappa}\int_s^t\|v_r\|_{L^{3p}}\|w_r\|_{L^{3p}}^{3p-1}\de r
\end{align*}
for any $\kappa\in(0,1)$ using Hölder inequality. Similar computations give
\begin{align*}
\left|\int_s^t\langle Q_1(v_r,w_r),w_r^{3p-3}\rangle\de r\right|&\lesssim \|X_{>n}\|_{\CC^\kappa}\int_s^t\|v_r\|_{L^{3p}}^3\|w_r\|_{L^{3p}}^{3p-3}\de r\\
&\qquad +\sum_{k=1}^3\|X_{>n}\|_{\CC^\kappa}^k\int_s^t\|v_r\|_{L^{3p}}^k\|w_r\|_{L^{3p}}^{3p-k}\de r
\end{align*}
where the implicit constant does not depend on $n$. Recall from Lemma \ref{Lem:truncation} that $v$ is controlled by $w$ uniformly in $n$. Thus, choosing $n$ large enough depending on $\delta,p$ and $X$, we get
$$
\left|\int_s^t\langle Q_1(v_r,w_r),w_r^{3p-3}\rangle\de r\right|\leq \delta\int_s^t\|w_r\|_{L^{3p}}^{3p}\de r
$$
which completes the proof.
\end{proof}

For the other nonlinear terms, note that only the regularizing properties of $\cdot\pl X_{>n}$ are used. As such, the estimates all end up being uniform in the truncation parameter $n$.

\medskip

\begin{lem}{\label{LEM:EstimateQ2}}
For any $\delta>0$, there exists a constant $c>0$ that depends only on $\delta,p,\varepsilon$ such that for any $0<s<t<T$, we have
$$
\left|\int_s^t\langle Q_2(v_r,\bz_r),w_r^{3p-3}\rangle\de r\right|\leq \delta\int_s^t\|w_r\|_{L^{3p}}^{3p}\de r + cK_t\big(1+\int_s^t\|w_r\|_{B_{p,\infty}^{2\varepsilon}}^p\de r\big)
$$
where the constant $K_t$ is defined by
\begin{align}{\label{DEF:K}}
K_t = \sup_{0\leq r\leq t}\|\bz_r^{(1)}\|_{\CC^{-\varepsilon}}^{2p} + \int_0^t\|\bz^{(2)}_r\|_{\CC^{-\varepsilon}}^{3p}\de r + \int_0^t\|\bz^{(3)}_r\|_{\CC^{-\varepsilon}}^{3p}\de r.
\end{align}
\end{lem}

\medskip

\begin{proof}
The major issue with
\begin{equation}
Q_2(v,\bz)=v^2\bz^{(1)} + v\bz^{(2)} + \bz^{(3)}
\end{equation}
is handling the quadratic part $v^2\bz^{(1)}$, to that end we further expand $Q_2$ as 
\begin{align*}
Q_2&=(v\pl X_{>n})^2\bz^{(1)} + w^2\bz^{(1)} + 2w(v\pl X_{>n})\bz^{(1)} + (v\pl X_{>n})\bz^{(2)} + w\bz^{(2)} + \bz^{(3)}
\end{align*}
since $v=\Gamma_nw$ is paracontrolled to take advantage of the regularizing effect of $\cdot\pl X_{>n}$. Starting with the $w^2$ term, we get
$$
|\langle w^2\bz^{(1)},w^{3p-3}\rangle|\lesssim \|\bz^{(1)}\|_{\CC^{-\varepsilon}}\|w^{3p-1}\|_{B^\varepsilon_{1,\infty}}
$$
since $\bz^{(1)}$ has spatial Hölder regularity $-\varepsilon$. Using the power product rule Lemma \ref{LEM:Besov} and interpolation Lemma \ref{LEM:InterpolationBesov} in Besov spaces, we bound $w^{3p-1}$ as
\begin{align*}
\|w^{3p-1}\|_{B^\varepsilon_{1,\infty}}&\lesssim\|w^{3p-2}\|_{L^{\frac{3p}{3p-2}}}\|w\|_{B^\varepsilon_{\frac{3p}{2},\infty}}\\
&\lesssim\|w\|_{L^{3p}}^{3p-\frac{3}{2}}\|w\|_{B^{2\varepsilon}_{p,\infty}}^{\frac{1}{2}}.
\end{align*} 
Using Hölder's estimate, we get
\begin{align*}
\left|\int_s^t\langle w_r^2\bz_r^{(1)},w_r^{3p-3}\rangle\de r\right|&\lesssim \int_s^t\|\bz^{(1)}\|_{\CC^{-\varepsilon}}\|w_r\|_{L^{3p}}^{3p-\frac{3}{2}}\|w_r\|_{B^{2\varepsilon}_{p,\infty}}^{\frac{1}{2}}\de r\\
&\lesssim\sup_{0\leq r\leq t}\|\bz_r^{(1)}\|_{\CC^{-\varepsilon}}\Big(\int_s^t\|w_r\|_{L^{3p}}^{3p}\de r\Big)^{\frac{3p-\frac{3}{2}}{3p}}\Big(\int_s^t\|w_r\|_{B^{2\varepsilon}_{p,\infty}}^{p}\de r\Big)^{\frac{1}{2p}}
\end{align*}
and thanks to the embedding $B^{2\varepsilon}_{p,\infty}\hookrightarrow L^{3p}$ together with Young's estimate, we have
\begin{align*}
\left|\int_s^t\langle w_r^2\bz_r^{(1)},w_r^{3p-3}\rangle\de r\right|&\leq \int_s^t\|\bz^{(1)}\|_{\CC^{-\varepsilon}}\|w_r\|_{L^{3p}}^{3p-\frac{3}{2}}\|w_r\|_{B^{2\varepsilon}_{p,\infty}}^{\frac{1}{2}}\de r\\
&\leq\delta \int_s^t\|w_r\|_{L^{3p}}^{3p}\de r + \big(\sup_{0\leq r\leq t}\|\bz_r^{(1)}\|_{\CC^{-\varepsilon}}\big)^{p}c(\delta)\int_s^t\|w_r\|^p_{B^{2\varepsilon}_{p,\infty}}\de r
\end{align*}
since $\bz^{(1)}$ is in $C_t\CC^{-\varepsilon}$. Turning to the $(v\pl X_{>n})^2$ term, we make use of the regularizing effects of the operator $\cdot\pl X_{>n}$
\begin{align*}
\left|\langle(v\pl X_{>n})^2\bz^{(1)},w^{3p-3}\rangle\right|&\lesssim \|\bz^{(1)}\|_{\CC^{-\varepsilon}}\|(v\pl X_{>n})^2\|_{\CC^{\varepsilon}}\|w^{3p-3}\|_{B^\varepsilon_{1,\infty}}\\
&\lesssim \|\bz^{(1)}\|_{\CC^{-\varepsilon}}\|v\|_{L^{3p}}^2\|w^{3p-3}\|_{B^\varepsilon_{1,\infty}}\\
&\lesssim \|\bz^{(1)}\|_{\CC^{-\varepsilon}}\|w\|_{L^{3p}}^2\|w^{3p-3}\|_{B^\varepsilon_{1,\infty}}.
\end{align*}
As before, the power product and interpolation rules yield 
\begin{align*}
\|w^{3p-3}\|_{B^\varepsilon_{1,\infty}}&\lesssim\|w^{3p-4}\|_{L^{\frac{3p}{3p-4}}}\|w\|_{B^\varepsilon_{\frac{3p}{4},\infty}}\\
&\lesssim \|w\|_{L^{3p}}^{3p-4+\frac{1}{2}}\|w\|_{B^{2\varepsilon}_{\frac{3p}{7},\infty}}^{\frac{1}{2}}
\end{align*}
and we conclude using once again Young estimate. Same goes for the mixed term $w(v\pl X_{>n})$
\begin{align*}
\left|\langle w(v\pl X_{>n})\bz^{(1)},w^{3p-3}\rangle\right|&\lesssim \|\bz^{(1)}\|_{\CC^{-\varepsilon}}\|v\pl X_{>n}\|_{\CC^{\varepsilon}}\|w^{3p-2}\|_{B^\varepsilon_{1,\infty}}\\
&\lesssim \|\bz^{(1)}\|_{\CC^{-\varepsilon}}\|v\|_{L^{3p}}\|w^{3p-2}\|_{B^\varepsilon_{1,\infty}}\\
&\lesssim \|\bz^{(1)}\|_{\CC^{-\varepsilon}}\|w\|_{L^{3p}}\|w^{3p-2}\|_{B^\varepsilon_{1,\infty}}\\
&\lesssim \|\bz^{(1)}\|_{\CC^{-\varepsilon}}\|w\|_{L^{3p}}^{3p-2}\|w\|_{B^\varepsilon_{p,\infty}}\\
&\lesssim \|\bz^{(1)}\|_{\CC^{-\varepsilon}}\|w\|_{L^{3p}}^{3p-\frac{3}{2}}\|w\|_{B^{2\varepsilon}_{\frac{3p}{5},\infty}}^{\frac{1}{2}}
\end{align*}
and we finally obtain that
\begin{align*}
\left|\int_s^t\left\langle v_r^2\bz^{(1)},w_r^{3p-3}\right\rangle\de r\right|&\leq \delta \int_s^t\|w_r\|^{3p}_{L^{3p}}\de r+\big(\sup_{0\leq r\leq t}\|\bz_r^{(1)}\|_{\CC^{-\varepsilon}}\big)^{2p}c(\delta)\int_s^t\|w_r\|^p_{B^{2\varepsilon}_{p,\infty}}\de r.
\end{align*}
The bound with the linear part or the constant part all follow from similar estimates, the only difference is that $\bz^{(2)}$ and $\bz^{(3)}$ are not continuous in time but rather $L^q$ for any $q<+\infty$. This is however not an issue as the homogeneity in $w$ in these terms is at most 1. For the linear terms, we get
\begin{align*}
\left|\langle (v\pl X_{>n})\bz^{(2)},w^{3p-3}\rangle\right|&\lesssim \|\bz^{(2)}\|_{\CC^{-\varepsilon}}\|v\pl X_{>n}\|_{\CC^{\varepsilon}}\|w^{3p-3}\|_{B^\varepsilon_{1,\infty}}\\
&\lesssim \|\bz^{(2)}\|_{\CC^{-\varepsilon}}\|v\|_{L^{3p}}\|w^{3p-3}\|_{B^\varepsilon_{1,\infty}}\\
&\lesssim \|\bz^{(2)}\|_{\CC^{-\varepsilon}}\|w\|_{L^{3p}}\|w^{3p-3}\|_{B^\varepsilon_{1,\infty}}\\
&\lesssim \|\bz^{(2)}\|_{\CC^{-\varepsilon}}\|w\|_{L^{3p}}^{3p-3}\|w\|_{B^\varepsilon_{\frac{3p}{4},\infty}}\\
&\lesssim \|\bz^{(2)}\|_{\CC^{-\varepsilon}}\|w\|_{L^{3p}}^{3p-\frac{5}{2}}\|w\|_{B^{2\varepsilon}_{\frac{3p}{5},\infty}}^{\frac{1}{2}}\\
&\leq \delta\|w\|_{L^{3p}}^{3p}+c(\delta)\|\bz^{(2)}\|_{\CC^{-\varepsilon}}^{\frac{6p}{5}}\|w\|_{B^{2\varepsilon}_{p,\infty}}^{\frac{3p}{5}}
\end{align*}
and 
\begin{align*}
\left|\langle w\bz^{(2)},w^{3p-3}\rangle\right|&\lesssim \|\bz^{(2)}\|_{\CC^{-\varepsilon}}\|w^{3p-2}\|_{B^\varepsilon_{1,\infty}}\\
&\lesssim \|\bz^{(2)}\|_{\CC^{-\varepsilon}}\|w\|_{L^{3p}}^{3p-3}\|w\|_{B^\varepsilon_{p,\infty}}\\
&\lesssim \|\bz^{(2)}\|_{\CC^{-\varepsilon}}\|w\|_{L^{3p}}^{3p-\frac{5}{2}}\|w\|_{B^{2\varepsilon}_{p,\infty}}^{\frac{1}{2}}\\
&\leq \delta\|w\|_{L^{3p}}^{3p}+c(\delta)\|\bz^{(2)}\|_{\CC^{-\varepsilon}}^{\frac{6p}{5}}\|w\|_{B^{2\varepsilon}_{p,\infty}}^{\frac{3p}{5}}.
\end{align*}
Integrating over time and applying Young's inequality, we obtain
\begin{align*}
\left|\int_s^t\left\langle v_r\bz_r^{(2)},w_r^{3p-3}\right\rangle\de r\right|\leq \delta\int_s^t\|w_r\|_{3p}^{3p}\de r + c(\delta)\int_0^t\|\bz^{(2)}_r\|_{\CC^{-\varepsilon}}^{3p}\de r \big(1+\int_s^t\|w_r\|_{B^{2\varepsilon}_{p,\infty}}\de r\big).
\end{align*}
As for the linear term, this is basically the same idea, while being slightly easier as we do not need to use an interpolation inequality
\begin{align*}
\left|\langle \bz^{(3)},w^{3p-3}\rangle\right|&\lesssim \|\bz^{(3)}\|_{\CC^{-\varepsilon}}\|w^{3p-3}\|_{B^\varepsilon_{1,\infty}}\\
&\lesssim \|\bz^{(3)}\|_{\CC^{-\varepsilon}}\|w\|_{L^{3p}}^{3p-4}\|w\|_{B^\varepsilon_{p,\infty}}\\
&\lesssim \delta\|w\|_{L^{3p}}^{3p}+c(\delta)\|\bz^{(3)}\|_{\CC^{-\varepsilon}}^{\frac{3p}{4}}\|w\|_{B^{\varepsilon}_{p,\infty}}^{\frac{3p}{4}}.
\end{align*}
Altogether, we obtain 
$$
\left|\int_s^t\left\langle Q_2(v_r,\bz_r),w_r^{3p-3}\right\rangle\de r\right|\leq \delta\int_s^t\|w_r\|_{L^{3p}}^{3p}\de r + c(\delta)K_t\big(1+\int_s^t\|w_r\|_{B^{2\varepsilon}_{p,\infty}}^p\de r).
$$
\end{proof}

\begin{lem}{\label{LEM:EstimateQ3}}
For any $\delta>0$, there exists a constant $c>0$ that depends only on $\delta,p,\varepsilon$ such that for any $0<s<t<T$, we have
$$
\left|\int_s^t\langle Q_3(v_r,\bz_r),w_r^{3p-3}\rangle\de r\right|\leq \delta\int_s^t\|w_r\|_{L^{3p}}^{3p}\de r + cK_t\big(1+\int_s^t\|w_r\|_{B_{p,\infty}^{2\varepsilon}}^p\de r\big)
$$
where $K_t$ is defined in \eqref{DEF:K}.
\end{lem}

\medskip

\begin{proof}
This is an easier counterpart of the previous proof as the regularizing operator $w\mapsto w\pl X_{>n}$ allows to estimate $w^{3p-3}$ in $L^{\frac{3p}{3p-3}}$ directly
$$
\langle Q_3(v,\bz),w^{3p-3}\rangle\lesssim \|v^2\bz^{(1)} + v\bz^{(2)} + \bz^{(3)}\|_{\CC^{-\varepsilon}}\|w\|_{L^{3p}}^{3p-3}
$$
and the first term on the right-hand side is dealt with as before.
\end{proof}

All in all, we obtain the following first a priori bound for $w$ taking $\delta<1$. In particular, it requires a regularity exponent strictly greater than one due to the roughness of the environment, the product with the space white noise $\xi\in\CC^{-1-\kappa}$ has to be well-defined.

\medskip

\begin{prop}{\label{LpBound}}
Let $p>\frac{4}{3\varepsilon}$ be an even integer and $n$ large enough depending on $p$. For $0<s<t<T$, we have
\begin{align*}
\|w_t\|_{L^{3p-2}}^{3p-2} + \int_s^t\|w_r\|_{L^{3p}}^{3p}\de r \lesssim \|w_s\|_{L^{3p-2}}^{3p-2} + K_t\big(1+ \int_s^t\|w_r\|_{B^{1+\varepsilon}_{p,\infty}}^p\de r\big)
\end{align*}
where the implicit constant depends only on $p$ and $\varepsilon$ and $K_t>0$ is given by \eqref{DEF:K}.
\end{prop}

\medskip

\begin{proof}
This is only a matter on putting together the estimates from Lemma \ref{LEM:EstimateNoiseParaW} to \ref{LEM:EstimateQ3}. Recall that testing the equation for $w$ \eqref{Eqw} against $w^{3p-3}$ and integrating in time yields
\begin{align*}
&\frac{1}{3p-2}\big(\|w_t\|_{L^{3p-2}}^{3p-2}-\|w_s\|_{L^{3p-2}}^{3p-2}\big) + (3p-3)\int_s^t\||\nabla w_r|^2w_r^{3p-4}\|_{L^1}\de r + \int_s^t\|w_r\|_{L^{3p}}^{3p}\de r\\
&\hspace{1cm} = -\int_s^t\langle \xi\ple w_r,w_r^{3p-3}\rangle\de r -\int_s^t\langle Q(v_r,w_r,\bz_r),w_r^{3p-3}\rangle\de r-\int_s^t\langle R_n(w_r),w_r^{3p-3}\rangle\de r
\end{align*}
for any $0<s<t<T$. Note that the term involving $\nabla w$ is nonnegative so that, using Hölder estimate on Lemma \ref{LEM:EstimateNoiseParaW} and \ref{LEM:EstimateR} as well as Lemma \ref{LEM:EstimateQ1}, \ref{LEM:EstimateQ2}, \ref{LEM:EstimateQ3} we end up with
\begin{align*}
&\frac{1}{3p-2}\|w_t\|_{L^{3p-2}}^{3p-2} + \int_s^t\|w_r\|_{L^{3p}}^{3p}\de r\leq \frac{1}{3p-2}\|w_s\|_{L^{3p-2}}^{3p-2}+\delta \int_s^t\|w_r\|_{L^{3p}}^{3p}\de r + c_\delta K_t\big(1+\int_s^t\|w_r\|_{B^{1+\varepsilon}_{p,\infty}}^p\de r\big)
\end{align*}
for any $\delta>0$. Choosing $\delta$ small enough, we can absorb the integrated $L^{3p}$-norm on the right hand side, hence the the estimate.
\end{proof}

\subsection{Higher regularity estimates on \texorpdfstring{$w$}{w}}

As we intend to use a \textit{coming down from infinity} argument, it remains to bound the rightmost term in order to close the estimate. This follows from the mild formulation of the equation.

\medskip

\begin{lem}{\label{EstimateHigherReg}}
For any $0<s<t<T$, we have
\begin{align*}
\int_s^t\|w_r\|^p_{B^{1+4\varepsilon}_{p,\infty}}\de r\lesssim 1 + \int_s^t\|e^{(r-s)\Delta}w_s\|^p_{B^{1+4\varepsilon}_{p,\infty}}\de r+ \tilde{K}_t\int_s^t\|w_r\|_{L^{3p}}^{3p}\de r
\end{align*}
with $\tilde{K_t}>0$ a constant of the form
\begin{align}{\label{DEF:tildeK}}
\tilde{K}_t=\big(\sup_{0\leq r\leq t}\|\bz^{(1)}\|_{\CC^{-\varepsilon}} + t^{a}\|\bz^{(2)}\|_{L_t^{q}\CC^{-\varepsilon}} + t^{a}\|\bz^{(3)}\|_{L_t^{q}\CC^{-\varepsilon}}\big)^b
\end{align}
for exponents $a,b>0$ and $q>1$ all depending only on $p$ and $\varepsilon$.
\end{lem}

\medskip

\begin{rem}
In the following, the exact exponents in $\tilde{K}$ may differ from one line to the other, but we still write $\tilde{K}$ for a constant having the same form as \eqref{DEF:tildeK}. Note that $\tilde{K}_t$ depends on the randomness of $\bz$ only on the interval $[0,t]$, this will come in handy later on.
\end{rem}

\medskip

\begin{proof}
Set $\gamma=1+4\varepsilon$ and consider the mild formulation
$$
w_t=e^{(t-s)\Delta}w_s - \int_s^te^{(t-r)\Delta}\left(\xi\ple w_r + R_n(w_r) + w_r^3 + Q(v_r,w_r,\bz_r)\right)\de r
$$
to take advantage of the regularizing effects of the heat operator $e^{(t-r)\Delta}$. Provided $p$ is large enough, and for any $\kappa<\frac{\gamma}{2}$, we have
\begin{align*}
\int_s^t\|e^{(t-r)\Delta}\big(\xi\ple w_r\big)\|_{B^\gamma_{p,\infty}}\de r&\lesssim\int_s^t(t-r)^{-\frac{\gamma+\kappa}{2}}\|\xi\ple w_r\|_{B_{p,\infty}^{-\kappa}}\de r\\
&\lesssim\int_s^t(t-r)^{-\frac{\gamma+\kappa}{2}}\|w_r\|_{B_{p,\infty}^{1+2\kappa}}\de r\\
&\lesssim\Big(\int_s^t(t-r)^{-\frac{(\gamma+\kappa)p}{2(p-1)}}\de r\Big)^{\frac{p-1}{p}}\Big(\int_s^t\|w_r\|_{B_{p,\infty}^{1+2\kappa}}^p\de r\Big)^{\frac{1}{p}}
\end{align*}
using Schauder estimates and Hölder inequality. We also get
\begin{align*}
\int_s^t\|e^{(t-r)\Delta}R_n(w_r)\|_{B^\gamma_{p,\infty}}\de r&\lesssim \left(\int_s^t(t-r)^{-\frac{(\gamma +\kappa)p}{2(p-1)}}\de r\right)^{\frac{p-1}{p}}\left(\int_s^t\|w_r\|_{B^{1+\kappa}_{p,\infty}}^p\de r\right)^{\frac{1}{p}}
\end{align*} 
and for the cubic part
\begin{align*}
\int_s^t\|e^{(t-r)\Delta}w_r^3\|_{B^\gamma_{p,\infty}}\de r \lesssim \left(\int_s^t(t-r)^{-\frac{\gamma p}{2(p-1)}}\de r\right)^{\frac{p-1}{p}}\left(\int_s^t\|w_r\|_{L^{3p}}^{3p}\de r\right)^{\frac{1}{p}}
\end{align*}
with similar arguments. Now for the collection of terms appearing in $Q$, recall that 
\begin{equation}
Q_1(v,w)=-v^3\pl X_{>n} +(v\pl X_{>n})^3+ 3(v\pl X_{>n})w^2+3(v\pl X_{>n})^2w
\end{equation} 
gathers cubic homogeneity terms which can be bounded using Hölder inequality and Lemma \ref{Lem:truncation} as
\begin{align*}
\|Q_1(w)\|_{L^p}&\lesssim\|v^3\|_{L^p} + \|v\|_{L^{3p}}^3 + \|v\|_{L^{3p}}\|w^2\|_{L^{\frac{3p}{2}}}+\|v\|_{L^{3p}}^2\|w\|_{L^{3p}}\\
&\lesssim\|w\|_{L^{3p}}^3
\end{align*}
using the continuity of $\Gamma_n$ hence
\begin{align*}
\int_s^t\|e^{(t-r)\Delta}Q_1(v_r,w_r)\|_{B^\gamma_{p,\infty}}\de r\lesssim \left(\int_s^t(t-r)^{-\frac{\gamma p}{2(p-1)}}\de r\right)^{\frac{p-1}{p}}\left(\int_s^t\|w_r\|_{L^{3p}}^{3p}\de r\right)^{\frac{1}{p}}
\end{align*}
integrating against $e^{(t-r)\Delta}$. Making use of the product rule Lemma \ref{LEM:Besov} to handle the quadratic terms, we get
\begin{align*}
\|Q_2(v,\bz)\|_{B^{-\varepsilon}_{p,\infty}}&\lesssim \|\bz^{(1)}\|_{\CC^{-\varepsilon}}\|v\|_{B^{2\varepsilon}_{\frac{3p}{2},\infty}}\|v\|_{L^{3p}} + \|\bz^{(2)}\|_{\CC^{-\varepsilon}}\|v\|_{B^{2\varepsilon}_{p,\infty}} + \|\bz^{(3)}\|_{\CC^{-\varepsilon}}\\
&\lesssim \|\bz^{(1)}\|_{\CC^{-\varepsilon}}\|w\|_{B^{2\varepsilon}_{\frac{3p}{2},\infty}}\|w\|_{L^{3p}} + \|\bz^{(2)}\|_{\CC^{-\varepsilon}}\|w\|_{B^{2\varepsilon}_{p,\infty}} + \|\bz^{(3)}\|_{\CC^{-\varepsilon}}.
\end{align*}
The analysis of the two last terms on the right-hand side is straightforward, we investigate the quadratic part. Integrating against ${e^{(t-r)\Delta}}$, we have
\begin{align*}
&\int_s^t\|e^{(t-r)\Delta}\bz^{(1)}_rv_r^2\|_{B^\gamma_{p,\infty}}\de r\lesssim\int_s^t(t-r)^{-\frac{\gamma+\varepsilon}{2}}\|\bz^{(1)}_r\|_{\CC^{-\varepsilon}}\|w_r\|_{B^{2\varepsilon}_{\frac{3p}{2},\infty}}\|w_r\|_{L^{3p}}\de r\\
&\hspace*{2cm}\lesssim \left(\int_s^t(t-r)^{-\frac{3(\gamma+\varepsilon) p}{2(3p-4)}}\|\bz^{(1)}_r\|_{\CC^{-\varepsilon}}^{\frac{3p}{3p-4}}\de r\right)^{\frac{3p-4}{3p}}\left(\int_s^t\|w_r\|_{B^{2\varepsilon}_{\frac{3p}{2},\infty}}^p\de r\right)^{\frac{1}{p}}\left(\int_s^t\|w_r\|_{L^{3p}}^{3p}\de r\right)^{\frac{1}{3p}}
\end{align*}
using Hölder estimate. For $\eps$ small enough and $p$ large enough, we have $\alpha_1:=\frac{3(1+5\varepsilon) p}{2(3p-4)}<1$ and the $\bz^{(1)}$ integral can be treated as
\begin{align*}
\int_s^t(t-r)^{-\frac{3(\gamma+\varepsilon) p}{2(3p-4)}}\|\bz^{(1)}_r\|_{\CC^{-\varepsilon}}^{\frac{3p}{3p-4}}\de r&\leq \frac{1}{1-\alpha_1}(t-s)^{1-\alpha_1}\sup_{0\leq r\leq t}\|\bz_r^{(1)}\|_{\CC^{\varepsilon}}^{\frac{3p}{3p-4}}.
\end{align*}
Again, interpolation Lemma \ref{LEM:InterpolationBesov} between $L^{3p}$ and $B^{4\varepsilon}_{p,\infty}$ yields
$$
\|w\|_{B^{2\varepsilon}_{\frac{3p}{2}},\infty}\lesssim\|w\|_{L^{3p}}^{\frac{1}{2}}\|w\|^{\frac{1}{2}}_{B^{4\varepsilon}_{p,\infty}}
$$
thus we get
\begin{align*}
\left(\int_s^t\|w_r\|_{B^{2\varepsilon}_{\frac{3p}{2},\infty}}^p\de r\right)^{\frac{1}{p}}\left(\int_s^t\|w_r\|_{L^{3p}}^{3p}\de r\right)^{\frac{1}{3p}}&\lesssim \left(\int_s^t\|w_r\|_{B^{4\varepsilon}_{p,\infty}}^{\frac{p}{2}}\|w_r\|_{L^{3p}}^{\frac{p}{2}}\de r\right)^{\frac{1}{p}}\left(\int_s^t\|w_r\|_{L^{3p}}^{3p}\de r\right)^{\frac{1}{3p}}\\
&\lesssim \left(\int_s^t\|w_r\|_{B^{4\varepsilon}_{p,\infty}}^{\frac{3p}{5}}\de r\right)^{\frac{5}{6p}}\left(\int_s^t\|w_r\|_{L^{3p}}^{3p}\de r\right)^{\frac{1}{2p}}\\
&\lesssim \left(\int_s^t\|w_r\|_{B^{4\varepsilon}_{p,\infty}}^p\de r\right)^\frac{1}{p}+\left(\int_s^t\|w_r\|_{L^{3p}}^{3p}\de r\right)^{\frac{1}{p}}
\end{align*}
using Young's estimate. The very same argument can be used to deal with the lower homogeneity terms of $Q_2$ and $Q_3$, although we make use of $\bz^{(2)}$ and $\bz^{(3)}$ being $L^q$ in time for any $q<+\infty$ to ensure that integrals of the type 
$$
\int_s^t(t-r)^{-\sigma}\|\bz^{(i)}_r\|_{\CC^{-\varepsilon}}\de r
$$
for $\sigma\in(0,1)$ and $i\in\{1,2\}$ are finite and uniformly bounded in $s$ by some time integral of $\|\bz^{(2)}\|_{\CC^{-\varepsilon}}$ and $\|\bz^{(3)}\|_{\CC^{-\varepsilon}}$. This proves that 
$$
\|w_t\|_{B^\gamma_{p,\infty}}\lesssim 1 + \|e^{(t-s)\Delta}w_s\|_{B^\gamma_{p,\infty}} + \tilde{K}_t\left(\int_s^t\|w_r\|_{B^{1+\kappa}_{p,\infty}}^p\de r\right)^\frac{1}{p}+\tilde{K}_t\left(\int_s^t\|w_r\|_{L^{3p}}^{3p}\de r\right)^{\frac{1}{p}}
$$
for any $0<\gamma<2$, $\kappa>0$ $p>1$ large enough, where the implicit constant depends only on $\varepsilon,\kappa,p$. Since $\gamma=1+4\varepsilon$, it follows from interpolation and Young inequalities that
\begin{align*}
\|w\|_{B^{1+\varepsilon}_{p,\infty}}^p&\lesssim \|w\|_{B^{1+4\varepsilon}_{p,\infty}}^{\frac{1+\varepsilon}{1+4\varepsilon}p}\|w\|_{L^{p}}^{\frac{3\varepsilon}{1+4\varepsilon}p}\\
&\lesssim \|w\|_{B^{1+4\varepsilon}_{p,\infty}}^{\frac{1+\varepsilon}{1+3\varepsilon}p}+\|w\|_{L^{3p}}^{3p}.
\end{align*}
Setting $\kappa=\varepsilon$ and integrating over time finally yield
\begin{align*}
\int_s^t\|w_r\|^p_{B^{1+4\varepsilon}_{p,\infty}}\de r\lesssim 1 + \int_s^t\|e^{(r-s)\Delta}w_s\|^p_{B^{1+4\varepsilon}_{p,\infty}}\de r+ \tilde{K}_t\int_s^t\|w_r\|_{L^{3p}}^{3p}\de r + \int_s^t\|w_r\|^{\frac{1+\varepsilon}{1+3\varepsilon}p}_{B^{1+4\varepsilon}_{p,\infty}}\de r.
\end{align*}
Since $\frac{1+\varepsilon}{1+4\varepsilon}<1$, the last term on the right-hand side can be absorbed by the left-hand side using Young's estimate, which concludes the proof.
\end{proof}

It only remains to turn the previous results into a proper integral inequality, this is the purpose of the following proposition.

\medskip

\begin{prop}\label{Prop:aprioriestimate}
Let $p>\frac{4}{3\varepsilon}$ be an even integer and $n$ large enough depending on $p$. For $0<s<t<T$, we have
\begin{align}{\label{APrioriEstimatew}}
\|w_t\|_{L^{3p-2}}^{3p-2}\le (c+\tilde{K}_t)\Big(1+\|w_s\|_{L^{3p-2}}^{3p-2}+\|w_s\|_{B^{1+4\varepsilon}_{p,\infty}}^{\frac{3p-2}{3}}\Big)
\end{align}
and its integrated counterpart
\begin{align*}
\int_s^t\left(\|w_r\|_{L^{3p-2}}^{3p-2}+\|w_r\|_{B^{1+4\varepsilon}_{p,\infty}}^{\frac{3p-2}{3}}\right)^\frac{3p}{3p-2}\de r\leq (c+\tilde{K}_t)\Big(1 +\|w_s\|_{L^{3p-2}}^{3p-2}+\|w_s\|_{B^{1+4\varepsilon}_{p,\infty}}^{\frac{3p-2}{3}}\Big)
\end{align*}
for a positive constant $c=c(p,\varepsilon,X)>0$ independent of time.
\end{prop}

\medskip

\begin{proof}
We collect the two main ingredients provided by the previous Lemmas \ref{LpBound} and \ref{EstimateHigherReg}
$$
\|w_t\|_{L^{3p-2}}^{3p-2} + \int_s^t\|w_r\|_{L^{3p}}^{3p}\de r \lesssim \|w_s\|_{L^{3p-2}}^{3p-2} + \int_s^t\|w_r\|_{B^{1+\varepsilon}_{p,\infty}}^p\de r
$$
and
$$
\int_s^t\|w_r\|^p_{B^{1+4\varepsilon}_{p,\infty}}\de r\lesssim 1 + \int_s^t\|e^{(r-s)\Delta}w_s\|^p_{B^{1+4\varepsilon}_{p,\infty}}\de r+ \int_s^t\|w_r\|_{L^{3p}}^{3p}\de r.
$$
Now set $\gamma:=(1+4\varepsilon)\frac{p-1}{p}$ so that 
\begin{align*}
\|w\|_{B^\gamma_{p,\infty}}^p&\lesssim\|w\|_{B^{1+4\varepsilon}_{p,\infty}}^{p-1}\|w\|_{L^p}\\
&\lesssim \|w\|_{B^{1+4\varepsilon}_{p,\infty}}^{\frac{3p-2}{3}}+\|w\|_{L^p}^{3p-2}.
\end{align*}
From the estimate in Lemma \ref{EstimateHigherReg}, we then get
\begin{align*}
\int_s^t\|w_r\|^p_{B^{1+\varepsilon}_{p,\infty}}\de r&\leq c + \delta\int_s^t\|w_r\|_{L^{3p}}^{3p}\de r + \delta \int_s^t\|e^{(r-s)\Delta}w_s\|_{B^{1+4\varepsilon}_{p,\infty}}^p\de r\\
&\leq c + \delta\int_s^t\|w_r\|_{L^{3p}}^{3p}\de r + \delta \|w_s\|_{B^{\gamma}_{p,\infty}}^p\\
&\leq c + \delta\int_s^t\|w_r\|_{L^{3p}}^{3p}\de r + c(\delta)\big(\|w_s\|_{B^{1+4\varepsilon}_{p,\infty}}^{\frac{3p-2}{3}}+\|w_s\|_{L^p}^{3p-2}\big).
\end{align*}
Pluging this in the right-hand side of Lemma \ref{LpBound}, we get
\begin{align*}
\|w_t\|_{L^{3p-2}}^{3p-2}+\int_s^t\|w_r\|_{L^{3p}}^{3p}\de r+\int_s^t\|w_r\|_{B^{1+4\varepsilon}_{p,\infty}}^p\de r\lesssim 1+\|w_s\|_{L^{3p-2}}^{3p-2}+\|w_s\|_{B^{1+4\varepsilon}_{p,\infty}}^{\frac{3p-2}{3}}
\end{align*}
where the implicit constant depends only on $\varepsilon, p$ and $\bz$ which complete the proof. As each term on the left hand side is non negative, this yields 
\begin{align}
\|w_t\|_{L^{3p-2}}^{3p-2}\lesssim 1+\|w_s\|_{L^{3p-2}}^{3p-2}+\|w_s\|_{B^{1+4\varepsilon}_{p,\infty}}^{\frac{3p-2}{3}}
\end{align}
and
\begin{align*}
\int_s^t\left(\|w_r\|_{L^{3p-2}}^{3p-2}+\|w_r\|_{B^{1+4\varepsilon}_{p,\infty}}^{\frac{3p-2}{3}}\right)^\frac{3p}{3p-2}\de r\lesssim 1 +\|w_s\|_{L^{3p-2}}^{3p-2}+\|w_s\|_{B^{1+4\varepsilon}_{p,\infty}}^{\frac{3p-2}{3}},
\end{align*}
hence the estimate. Note that the implicit constants in these estimates only involve the noise terms $\bz$ through the use of all the previous lemmata, hence the particular form $c+\tilde{K}_t$.
\end{proof}

\subsection{Global well-posedness and estimates on the solution}{\label{SEC:estimates}}

Recall the following comparison result, Lemma $7.3$ from Mourrat and Weber \cite{MW17a}. This lemma ensures that the behavior of $F$ away from 0 is controlled uniformly in the initial condition, this is the so-called \textit{coming down from infinity}.

\medskip

\begin{lem}{\label{ComingDownFromInfinity}}
Let $\tau>s\ge0$ and $F:[s,\tau)\to\R_+$ be continuous such that for any $s\leq s_1<s_2<\tau$, we have
$$
\int_{s_1}^{s_2}F^\lambda(r)\de r\le cF(s_1)
$$
for fixed $\lambda>1$ and $c>0$. Then there exists $N\in\N^*$ and $s=t_0<t_1<\cdots<t_N=\tau$ such that 
$$
F(t_n)\lesssim c^{\frac{1}{\lambda-1}}(t_{n+1}-s)^{-\frac{1}{\lambda-1}}
$$
for $0\le n\le N-1$ where the implicit constant depends only on $\lambda$.
\end{lem}

\medskip

\begin{proof}
We construct the sequence $t_0,\cdots, t_N$ by induction. Define $t_0=s$ and $$t_1^*=t_0+c2^\lambda F(t_0)^{1-\lambda},$$ if $t_1^*\geq \tau$ then set $N=1$ and $t_1=\tau$ so that 
$$
F(t_0)=F(s)=(t_1^*-s)^{-\frac{1}{\lambda-1}}c^{\frac{1}{\lambda-1}}2^{\frac{\lambda}{\lambda-1}}\lesssim (t_1-s)^{-\frac{1}{\lambda-1}}c^{\frac{1}{\lambda-1}}.
$$
Otherwise, taking $s_1=s$ and $s_2=t_1^*$, we have
$$
(t_1^*-s) \min_{s\leq r\leq t_1^*} F^\lambda(r)\leq \int_s^{t_1^*} F^\lambda(r)\de r\leq c F(s)
$$
thus
$$
\min_{s\leq r\leq t_1^*} F^\lambda(r) \leq 2^{-\lambda}F^\lambda(s)
$$
and choose $t_1$ as 
$$
t_1=\inf\Argmin_{[t_0,t_1^*]}F
$$
so that in the end
$$
F(t_1)\leq \frac{1}{2}F(t_0)\quad\text{and}\quad t_1-t_0\leq c2^{\lambda}F^{1-\lambda}(t_0). 
$$
Assume that $t_0<t_1<\cdots<t_n<\tau$ have been constructed as above and set 
$$
t_{n+1}^*=t_n+c2^{\lambda}F^{1-\lambda}(t_n).
$$
Either $t_{n+1}^*\geq\tau$, in which case we can set $N=n+1$ and $t_{n+1}=\tau$ and conclude, or $t_{n+1}^*<\tau$ and we set 
$$
t_{n+1}=\inf\Argmin_{[t_n,t_{n+1}^*]}F
$$
so that, similarly to the initial step, 
\begin{equation}{\label{EQ:ProofComingDown}}
F(t_{n+1})\leq \frac{1}{2}F(t_0)\quad\text{and}\quad t_{n+1}-t_n\leq c2^{\lambda}F^{1-\lambda}(t_n). 
\end{equation}
Note that the construction procedure terminates after a finite number of steps, indeed, by construction
$$
t_{n+1}^*-t_n=c2^{\lambda}F(t_n)^{1-\lambda}\geq 2^{\lambda-1}(t_n^*-t_{n-1})
$$
thanks to the first part of \eqref{EQ:ProofComingDown}, thus the time step is increased by a factor $2^{\lambda-1}>1$ at each step of the procedure, thus $t_{n+1}$ necessarily exceeds $\tau$ after a finite number of steps.\\ 
As for the estimate, it suffices to write the increments as a telespoci sum
$$
t_{n+1}-s=\sum_{k=0}^{n}(t_{k+1}-t_k)\leq c2^\lambda\sum_{k=0}^{n}F^{1-\lambda}(t_k)
$$
where the last estimate follows by construction. Moreover, iterating the first part of \eqref{EQ:ProofComingDown} and recallig that $1-\lambda<0$, we get
$$
F^{1-\lambda}(t_k)\leq 2^{(1-\lambda)}F^{1-\lambda}(t_{k+1})\leq 2^{(1-\lambda)(n-k)}F^{1-\lambda}(t_n)
$$
for any $k\leq n$. Thus
\begin{align*}
t_{n+1}-s&=c2^\lambda F^{1-\lambda}(t_n)\sum_{k=0}^{n}2^{(1-\lambda)(n-k)}\\
&\leq c2^\lambda F^{1-\lambda}(t_n)\sum_{k\geq 0}2^{k(1-\lambda)}
\end{align*}
where the last series is convergent. Thus 
$$
F(t_n)\lesssim c^{\frac{1}{\lambda-1}}(t_{n+1}-s)^{-\frac{1}{\lambda-1}}
$$
where the implicit constant depends only on $\lambda$ and the proof is complete.
\end{proof}

\medskip

\begin{prop}{\label{ComingDown}}
Let $p>\frac{4}{3\varepsilon}$ be an even integer and $n$ large enough. We have
$$
(1\wedge\sqrt{t})\|w_t\|_{L^{3p-2}}\leq C\wedge \tilde{K}_t\leq C\wedge \tilde{K}_T
$$
for all $0<t<T$. Moreover, the right-hand side does not depend on $u_0$ and only depends on $\bz$ on the time interval $[0,T]$. 
\end{prop}

\medskip

\begin{proof}
Let $0<s<t<T$ and consider 
$$
F:r\in[s,T)\mapsto \|w_r\|_{L^{3p-2}}^{3p-2}+\|w_r\|_{B^{1+4\varepsilon}_{p,\infty}}^{\frac{3p-2}{3}}
$$
which satisfies the assumptions of Lemma \ref{ComingDownFromInfinity}, this is the a priori estimates from Proposition \ref{Prop:aprioriestimate}. Define also
$$
\tau:=\inf\left\{r\geq s, F(r)\leq 1\right\}\wedge\, T
$$
so that $F(r)\geq1$ for any $r\in[s,\tau]$. If $\tau\geq t$, then applying Lemma \ref{ComingDownFromInfinity}, we obtain $s=t_0<\cdots<t_{N+1}=\tau$ such that
$$
F(t_n)\lesssim (t_{n+1}-s)^{-\frac{3p-2}{2}}
$$
since $\lambda=\frac{3p}{3p-2}>1$. Let $n$ be such that $t_n\leq t<t_{n+1}$, then using \eqref{APrioriEstimatew}, we get
\begin{align}
\|w_t\|_{L^{3p-2}}^{3p-2}&\lesssim 1+F(t_n)\\ 
&\lesssim F(t_n)\\ 
&\lesssim (t_{n+1}-s)^{-\frac{3p-2}{2}}\\
&\lesssim (t-s)^{-\frac{3p-2}{2}}
\end{align}
since $t_n\le t\le\tau$. When $\tau< t$, note that the previous case ensure that 
$$
\|w_{\tau}\|_{L^{3p-2}}^{3p-2}\lesssim (\tau-s)^{-\frac{3p-2}{2}},
$$
so that applying \eqref{APrioriEstimatew} again, we get
$$
\|w_t\|_{L^{3p-2}}^{3p-2}\lesssim 1+F(\tau)\lesssim (\tau-s)^{-\frac{3p-2}{2}}.
$$
Taking $s=\frac{t}{2}$, we finally obtain
$$
(1\wedge\sqrt{t})\|w_t\|_{L^{3p-2}}\leq C\vee \tilde{K}_t\leq C\vee \tilde{K}_T
$$
for $0<t<T$ where the dependency of the upper-bound comes from the fact that the $\bz$ terms in $\tilde{K}_t$ only take into account what happens between the initial time and the current time, being 0 and $t$ respectively.
\end{proof}

With the embedding $L^{3p-2}\hookrightarrow \CC^{-\varepsilon}$ for $p>\frac{4}{3\varepsilon}$ and the equivalence of the norms $\|v_t\|_{L^{3p-2}}$ and $\|w_t\|_{L^{3p-2}}$, this is enough to prove global well-posedness. Recall that the constant $\widetilde K_T>0$ is defined by
\begin{equation}
\widetilde K_T=\big(\sup_{0\leq t\leq T}\|\bz^{(1)}\|_{\CC^{-\varepsilon}} + t^{a}\|\bz^{(2)}\|_{L_t^{q}\CC^{-\varepsilon}} + t^{a}\|\bz^{(3)}\|_{L_t^{q}\CC^{-\varepsilon}}\big)^b
\end{equation}
and almost surely finite for any $T>0$.

\medskip

\begin{cor}{\label{COR:Globalv}}
Let $0<\varepsilon<\frac{1}{3}$ and $u_0\in\CC^{-\varepsilon}$. Consider $v$ the solution to
\begin{align*}
\begin{cases}
\partial_tv+\H v + v^3 + 3v^2\z^{(1)} + 3v\z^{(2)}+\z^{(3)}=0\\
v_{|t=0}=u_0
\end{cases}    
\end{align*}
obtained in Proposition \ref{Localwpv} and $T>0$ the maximal time of existence. We have
\begin{equation}
(1\wedge\sqrt{t})\|v(t)\|_{\CC^{-\varepsilon}}\leq 1+\widetilde K_T
\end{equation}
for any $0\le t<T$. In particular, this implies that $T=+\infty$ thus the solution is global.
\end{cor}

\medskip

\begin{proof}
The estimate follows directly from Proposition \ref{ComingDown} applied with the enhanced data $\z$, together with Lemma \ref{Controlvwunif} and the embedding $L^{3p-2}\hookrightarrow \Cc^{-\varepsilon}$ for $p$ large enough. The fact that the solution is global is then a consequence of the blow-up criterion in the local well-posedness theory \ref{Localwpv}.
\end{proof}

To conclude this section, we give a number of probabilistic bounds on the solutions that follow from the previous estimates. This will be important in the following section to study the solution as an infinite dimensional stochastic process.

\medskip

\begin{cor}{\label{COR:GlobalEstimatev}}
For any $0<\varepsilon<\frac{1}{3}$, $u_0\in\CC^{-\varepsilon}$ and $p\geq1$, there exists positive constants $C,\gamma$ such that
$$
\E\left[\sup_{0\leq t\leq 1}\|v(t)\|_{\CC^{-\varepsilon}}^p\right]\leq C(1+\|u_0\|_{\CC^{-\varepsilon}})^\gamma.
$$
\end{cor}

\medskip

\begin{proof}
Let $T$ be the time obtained by the fixed point argument in \ref{Localwpv}
$$
T=c\left(\frac{\|u_0\|_{\CC^{-\varepsilon}}}{(1+\|u_0\|_{\CC^{-\varepsilon}}+\|\z\|_{L^q\big((0,1);\CC^{-\varepsilon}\big)^3})^3}\right)^{-\frac{1}{\frac{1}{q'}-3\frac{\sigma+\varepsilon}{2}}}
$$
for some constant $c>0$ and 
$$
t^*=1\wedge T\wedge c(1+\|u_0\|_{\CC^{-\varepsilon}})^{-\frac{2}{\frac{1}{q'}-3\frac{\sigma+\varepsilon}{2}}}.
$$
Then, for $0<t<t^*$, the fixed point bound \eqref{localbound} yields
$$
\|v(t)\|_{\CC^{-\varepsilon}}\leq C(1+\|u_0\|_{\CC^{-\varepsilon}}).
$$
For $t^*\leq t\leq1$, we use Proposition \ref{ComingDown} to get
\begin{align*}
\|v(t)\|_{\CC^{-\varepsilon}}&\leq \big(1\vee {t^*}^{-\frac{1}{2}}\big)\tilde{K}_1\\
&\leq c^{-\frac{1}{2}}(1+\|u_0\|_{\CC^{-\varepsilon}})^{\frac{1}{\frac{1}{q'}-3\frac{\sigma+\varepsilon}{2}}}\tilde{K}_1
\end{align*}
and we conclude by gluing the estimates on $[0,t^*]$ and $[t^*,1]$ together, taking the $p$-th power and the expectation, using that $\tilde{K}_1$ has finite moments.
\end{proof}

Note that another choice of ansatz for the solution $u$ to \eqref{AndersonSQE} can be made
$$
u=\tilde{v}+\llp+e^{-t\H}u_0
$$
so that $\tilde{v}$ satisfies equation \eqref{model} with null initial condition and shifted noise data $\tilde{\z}=(\tilde{\z}^{(1)},\tilde{\z}^{(2)},\tilde{\z}^{(3)})$ with
\begin{align*}
\tilde{\z}^{(1)}&:=3\z^{(1)}+3e^{-t\H}u_0,\\
\tilde{\z}^{(2)}&:=3\z^{(2)}+6e^{-t\H}u_0\z^{(1)} + 3\big(e^{-t\H}u_0\big)^2,\\
\tilde{\z}^{(3)}&:=\z^{(3)}+3(e^{-t\H}u_0)^2\z^{(1)} +3e^{-t\H}u_0\z^{(2)} + \big(e^{-t\H}u_0\big)^3
\end{align*}
and falls under the scope of Proposition \ref{Localwpv} and \ref{COR:Globalv}. The main interest of $\tilde{v}$ is that it has better regularity since $\tilde{v}(0)=0$ is of better regularity than $u_0\in\CC^{-\eps}$. Of course, one has to be careful here since the data $\tilde{\z}$ depends on the initial data, this will be useful to prove Corollary \ref{COR:SquareSobolevuN}.

\medskip

\begin{lem}{\label{LEM:EstimateL2tildev}}
Let $0<\varepsilon<\frac{1}{3}$, $u_0\in\CC^{-\varepsilon}$ and $\tilde{v}$ be the solution to
\begin{align}{\label{EQ:tildev}}
\begin{cases}
\partial_t\tilde{v}+\H \tilde{v} + \mathbf{f}(\tilde{v},\tilde{\z}) =0\\
\tilde{v}_{|t=0}=0.             
\end{cases}        
\end{align}
For any $p\geq1$, we have
$$
\E\left[\sup_{0\leq t\leq 1} \left(\|\tilde{v}(t)\|_{L^2}^2+\int_0^t\|\tilde{v}(s)\|_{\D^1}^2\de s\right)^p\right]\leq C(1+\|u_0\|_{\CC^{-\varepsilon}})^\gamma
$$
for some constants $C,\gamma>0$ depending only on $\varepsilon$ and $p$.
\end{lem}

\medskip

\begin{proof}
For $t>0$, define the energy
$$
\Psi(t):=\frac{1}{2}\|\tilde{v}_t\|_{L^2}^2+\frac{1}{2}\int_0^t\|\tilde{v}_s\|_{\D^1}^2\de s.
$$
Note that running the same fixed-point argument as in Proposition \ref{Localwpv} in $L^2$ based Sobolev spaces ensures that $\Psi$ is differentiable in time, with derivative given by
\begin{align*}
\Psi'(t)&= -\frac{1}{2}\|\tilde{v}_t\|_{\D^1}^2-\|\tilde{v}_t\|_{L^4}^4-\int_{\T^2}\tilde{v}_t^3\tilde{\z}^{(1)}_t\de x-\int_{\T^2}\tilde{v}_t^2\tilde{\z}^{(2)}_t\de x-\int_{\T^2}\tilde{v}_t\tilde{\z}^{(3)}_t\de x.
\end{align*}
Note that this is essentially the same as testing the equation \ref{EQ:tildev} against $\tilde{v}$. It remains to bound the cross product terms to control them with the $\D^1$ and $L^4$ norms of $\tilde{v}$. Note that each integral can be bounded as a $\CC^{-\varepsilon},W^{\varepsilon,1}$ duality pairing with
$$
\left|\int_{\T^2}\tilde{v}_t^{j+1}\tilde{\z}^{(3-j)}_t\de x\right|\lesssim \|\tilde{v}_t^{j+1}\|_{W^{\varepsilon,1}}\|\tilde{\z}^{(3-j)}_t\|_{\CC^{-\varepsilon}}
$$
where $j\in\{0,1,2\}$. For $\eta\in(0,1)$ small enough, interpolation and Young inequalities yield
\begin{align*}
\|\tilde{v}_t^{j+1}\|_{W^{\varepsilon,1}}&\lesssim \|\tilde{v}^{j+1}_t\|_{W^{1-\eta,1}}^{\frac{\varepsilon}{1-\eta}}\|\tilde{v}^{j+1}_t\|_{L^1}^{\frac{1-\eta-\varepsilon}{1-\eta}}\\
&\lesssim \|\tilde{v}^{j+1}_t\|_{W^{1-\eta,1}}^{\frac{\varepsilon}{1-\eta}r}+\|\tilde{v}^{j+1}_t\|_{L^1}^{\frac{1-\eta-\varepsilon}{1-\eta}r'}
\end{align*}
where $r\geq1$ will be chosen later and $\frac{1}{r}+\frac{1}{r'}=1$. Using the embedding $H^{1-\eta}\hookrightarrow W^{1-\eta,1}$ and $L^4\hookrightarrow L^{j}$ for $j\in\{0,1,2\}$, the power product can be bounded using Lemma \ref{LEM:prodHugo} to get
\begin{align*}
\|\tilde{v}_t^{j+1}\|_{W^{\varepsilon,1}}&\lesssim \|\tilde{v}_t\|_{H^{1-\frac{\eta}{2^j}}}^{(j+1)\frac{\varepsilon}{1-\eta}r}+\|\tilde{v}_t\|_{L^{j+1}}^{(j+1)\frac{1-\eta-\varepsilon}{1-\eta}r'}\\
&\lesssim \|\tilde{v}_t\|_{\D^1}^{(j+1)\frac{\varepsilon}{1-\eta}r}+\|\tilde{v}_t\|_{L^4}^{(j+1)\frac{1-\eta-\varepsilon}{1-\eta}r'}.
\end{align*}
We need to tune the exponents so that the $\D^1$ norm (resp. $L^4$ norm) is at most to the power $2$ (resp. $4$), hence the following conditions
$$
(j+1)\frac{\varepsilon}{1-\eta}r <2\qquad\text{and}\qquad (j+1)\frac{1-\eta-\varepsilon}{1-\eta}r'<4
$$
for $j\in\{0,1,2\}$. The second one happens as soon as $r>3$, i.e. $r'<\frac{4}{3}$, and the first one follows by taking $\varepsilon$ small enough. Hence, using again Young's inequality, we have
$$
\left|\int_{\T^2}\tilde{v}_t^{j+1}\tilde{\z}^{(3-j)}_t\de x\right|\leq C_1(\delta)\|\tilde{\z}^{(3-j)}_t\|_{\CC^{-\varepsilon}}^{q_1}+\delta^{q_1'}\|\tilde{v}_t\|_{\D^1}^2+C_2(\delta)\|\tilde{\z}^{(3-j)}_t\|_{\CC^{-\varepsilon}}^{q_2}+\delta^{q_2'}\|\tilde{v}_t\|_{L^4}^4
$$
for any $\delta>0$ and 
$$
\frac{1}{q_1'}=\frac{(j+1)r}{2}\frac{\varepsilon}{1-\eta} \qquad\text{and}\qquad \frac{1}{q_2'}=\frac{(j+1)r'}{4}\frac{1-\eta-\varepsilon}{1-\eta}.
$$
Taking $\delta$ small enough ensures that the cross product terms can be absorbed so that
$$
\Psi'(t)\lesssim \|\tilde{\z}^{(1)}_t\|_{\CC^{-\varepsilon}}^{\tilde{q}_1}+\|\tilde{\z}^{(2)}_t\|_{\CC^{-\varepsilon}}^{\tilde{q}_2}+\|\tilde{\z}^{(3)}_t\|_{\CC^{-\varepsilon}}^{\tilde{q}_3}
$$
for some exponents $\tilde{q}_1,\tilde{q}_2,\tilde{q}_3\geq1$. Since $\Psi(0)=0$, integration in time yields 
$$
\sup_{0\leq t\leq1}\Psi(t)\lesssim \int_0^1\left(\|\tilde{\z}^{(1)}_s\|_{\CC^{-\varepsilon}}^{\tilde{q}_1}+\|\tilde{\z}^{(2)}_s\|_{\CC^{-\varepsilon}}^{\tilde{q}_2}+\|\tilde{\z}^{(3)}_s\|_{\CC^{-\varepsilon}}^{\tilde{q}_3}\right)\de s
$$
and taking the $p$-th moment, the estimate follows from the fact that $\tilde{z}\in L^p\Big(\Omega;L^q\big((0,1);\CC^{-\varepsilon}\big)\Big)$ for any $p,q\geq1$.
\end{proof}

This in particular implies that $\tilde{v}$ has positive spatial Sobolev regularity, uniformly in time.

\medskip

\begin{prop}{\label{PROP:BoundTildev}}
Let $0<\varepsilon<\sigma$ be small enough, $u_0\in\CC^{-\varepsilon}$ and $p\geq1$. We have
$$
\E\left[\sup_{0\leq t\leq 1} \|\tilde{v}(t)\|_{H^\sigma}^p\right]\leq C(1+\|u_0\|_{\CC^{-\varepsilon}}^\gamma)
$$
for some $C,\gamma>0$ depending only on $\varepsilon, \sigma$ and $p$.
\end{prop}

\medskip

\begin{proof}
Write the mild formulation for $\tilde{v}$
$$
\tilde{v}_t=-\int_0^te^{-(t-s)\H}\left(\tilde{v}_s^3 + \tilde{v}_s^2\tilde{\z}^{(1)}_s+\tilde{v}_s\tilde{\z}_s^{(2)}+\tilde{\z}_s^{(3)} \right)\de s.
$$
Taking the $H^\sigma$ norm inside and using Schauder estimates for $\H$, this reads
$$
\|\tilde{v}_t\|_{H^\sigma}\lesssim\int_0^t(t-s)^{-\frac{\sigma+\varepsilon}{2}}\left(\|\tilde{v}_s^3\|_{H^{-\varepsilon}} + \|\tilde{v}_s^2\tilde{\z}^{(1)}_s\|_{H^{-\varepsilon}}+\|\tilde{v}_s\tilde{\z}_s^{(2)}\|_{H^{-\varepsilon}}+\|\tilde{\z}_s^{(3)}\|_{H^{-\varepsilon}}\right)\de s.
$$
We now handle each term seperately using the embedding $H^{\sigma+\frac{1}{2}}\hookrightarrow B^\sigma_{4,4}$ and interpolation in Besov spaces. We have
\begin{align*}
\|\tilde{v}_t^3\|_{H^{-\varepsilon}}&\lesssim \|\tilde{v}_t^2\|_{H^{\sigma}}\|\tilde{v}_t\|_{C^{-\varepsilon}}\\
&\lesssim \|\tilde{v}_t\|_{H^{\sigma+\frac{1}{2}}}^2\|\tilde{v}_t\|_{C^{-\varepsilon}}\\
&\lesssim \|\tilde{v}_t\|_{H^{1-\delta}}^{\frac{2\sigma+1}{1-\delta}}\|\tilde{v}_t\|_{L^2}^{\frac{2-2\delta-(2\sigma+1)}{1-\delta}}\|\tilde{v}_t\|_{C^{-\varepsilon}}\\
&\lesssim \|\tilde{v}_t\|_{\D^1}^{\frac{2\sigma+1}{1-\delta}}\|\tilde{v}_t\|_{L^2}^{\frac{2-2\delta-(2\sigma+1)}{1-\delta}}\|\tilde{v}_t\|_{C^{-\varepsilon}}
\end{align*}
for any $0<\delta<\frac{1}{2}-\sigma$ using the embedding $\D^1\hookrightarrow H^{1-\delta}$. Set $\delta=\frac{1-4\sigma}{3}$ so that $\frac{2\sigma+1}{1-\delta}=\frac{3}{2}$. Young's inequality then yields
$$
\|\tilde{v}_t^3\|_{H^{-\varepsilon}}\lesssim \|\tilde{v}_t\|_{\D^1}^{\frac{7}{4}} + \|\tilde{v}_t\|_{L^2}^{\alpha} + \|\tilde{v}_t\|_{C^{-\varepsilon}}^\beta
$$
for some exponents $\alpha,\beta>0$. Note that the exact form of $\alpha,\beta$ do not matter as the $L^2$ and $\CC^{-\varepsilon}$ norms of $\tilde{v}$ are handled uniformly in time by Corollary \ref{COR:GlobalEstimatev} and Lemma \ref{LEM:EstimateL2tildev}, our concern is on the exponent on the $\D^1$ norm which has to be less than $2$ as we can only bound its $L^2$ norm in time. Integrating in time and using Hölder inequality, the first term yields
\begin{align*}
\int_0^t(t-s)^{-\frac{\sigma+\varepsilon}{2}}\|\tilde{v}_t^3\|_{H^{-\varepsilon}}\de s&\lesssim \int_0^t(t-s)^{-\frac{\sigma+\varepsilon}{2}}\left(\|\tilde{v}_s\|_{\D^1}^{\frac{7}{4}} + \|\tilde{v}_s\|_{L^2}^{\alpha} + \|\tilde{v}_s\|_{C^{-\varepsilon}}^\beta\right)\de s\\
&\lesssim \left(\int_0^t\|\tilde{v}_s\|_{\D^1}^{2}\de s\right)^{\frac{7}{8}} + \sup_{0\leq t\leq 1} \big(\|\tilde{v}_s\|_{L^2}^\alpha+\|\tilde{v}_s\|_{C^{-\varepsilon}}^\beta\big)
\end{align*}
provided $\sigma+\varepsilon<\frac{1}{4}$ so that every time integral is bounded uniformly in $t\in[0,1]$.
A similar treatment allows to deal with the square term with
\begin{align*}
\|\tilde{v}_t^2\tilde{\z}^{(1)}_t\|_{H^{-\varepsilon}}&\lesssim \|\tilde{v}_t\|_{H^{\sigma+\frac{1}{2}}}^2\|\tilde{\z}^{(1)}_t\|_{C^{-\varepsilon}}\\
&\lesssim \|\tilde{v}_t\|_{H^{1-\delta}}^{\frac{2\sigma+1}{1-\delta}}\|\tilde{v}_t\|_{L^2}^{\frac{2-2\delta-(2\sigma+1)}{1-\delta}}\|\tilde{\z}_t\|_{C^{-\varepsilon}}\\
&\lesssim \|\tilde{v}_t\|_{\D^1}^{\frac{2\sigma+1}{1-\delta}}\|\tilde{v}_t\|_{L^2}^{\frac{2-2\delta-(2\sigma+1)}{1-\delta}}\|\tilde{\z}^{(1)}_t\|_{C^{-\varepsilon}}\\
&\lesssim\|\tilde{v}_t\|_{\D^1}^{\frac{7}{4}} + \|\tilde{v}_t\|_{L^2}^{\alpha} + \|\tilde{\z}^{(1)}_t\|_{C^{-\varepsilon}}^\beta\\
&\lesssim \|\tilde{v}_t\|_{\D^1}^{\frac{7}{4}} + \|\tilde{v}_t\|_{L^2}^{\alpha} + \|\z^{(1)}_t\|_{C^{-\varepsilon}}^\beta+ \|u_0\|_{\CC^{-\varepsilon}}^\beta
\end{align*}
and the linear one is more straightforward
\begin{align*}
\|\tilde{v}_t\tilde{\z}^{(2)}_t\|_{H^{-\varepsilon}}&\lesssim \|\tilde{v}_t\|_{H^{\sigma}}\|\tilde{\z}^{(2)}_t\|_{C^{-\varepsilon}}\\
&\lesssim \|\tilde{v}_t\|_{\D^1}\|\tilde{\z}^{(2)}_t\|_{C^{-\varepsilon}}\\
&\lesssim \|\tilde{v}_t\|_{\D^1}^2+\|\tilde{\z}^{(2)}_t\|_{C^{-\varepsilon}}^2\\
&\lesssim \|\tilde{v}_t\|_{\D^1}^2+\|\z^{(2)}_t\|_{C^{-\varepsilon}}^2+t^{-(\sigma+\varepsilon)}\|u_0\|_{\CC^{-\varepsilon}}^2\|\z_t^{(1)}\|_{\CC^{-\varepsilon}}^2+t^{-2(\sigma+\varepsilon)}\|u_0\|_{\CC^{-\varepsilon}}^4\\
&\lesssim \|\tilde{v}_t\|_{\D^1}^2+\|\z^{(2)}_t\|_{C^{-\varepsilon}}^2+\|\z_t^{(1)}\|_{\CC^{-\varepsilon}}^4+t^{-2(\sigma+\varepsilon)}\|u_0\|_{\CC^{-\varepsilon}}^4.
\end{align*}
Integrating in time and using Hölder inequality as before, this yields
$$
\|\tilde{v}_t\|_{H^\sigma}\lesssim \left(\int_0^1\|\tilde{v}_s\|_{\D^1}^2\de s\right)^\alpha + \sup_{0\leq t\leq1} \|\tilde{v}_t\|_{L^2}^\beta+\sup_{0\leq t\leq1}\|\tilde{v}_t\|_{C^{-\varepsilon}}^\gamma + \|\z\|_{L^q((0,1);\CC^{-\varepsilon})^3} + \|u_0\|_{\CC^{-\varepsilon}}^\eta
$$
for some exponents $\alpha, \beta, \gamma, q$ large enough. The conclusion then follows using Corollary \ref{COR:GlobalEstimatev} and Lemma \ref{LEM:EstimateL2tildev} for the first three terms, and finiteness of the moments of $\|\z\|_{L^q((0,1);\CC^{-\varepsilon})^3}$.
\end{proof}

\begin{rem}
Note that \ref{COR:Globalv}, \ref{COR:GlobalEstimatev}, \ref{PROP:BoundTildev} and  are stated for the solution to the limiting equation, but they all hold true for the renormalized solution $v_N$ and $\tilde{v}_N$, with the right-hand side being uniform in $N$.
\end{rem}

\section{Strong Feller property for Anderson \texorpdfstring{$\Phi^4_2$}{PHI42} model}\label{Sec:strongFeller}

We put the previous estimates to use to investigate the properties of the transition semigroup of equation \eqref{AndersonSQE}
\begin{align*}
P_t\varphi(u_0):=\E\Big[\varphi\big(u(t\pv u_0)\big)\Big]
\end{align*} 
where $u(t\pv u_0)$ denote the solution at time $t$ starting from $u_0\in\CC^{-\eps}$. We proved in the previous section that solutions are global, we work here for any time $t\ge0$. In this section, it is convenient to consider
\begin{equation}
\llp_{s,t}=\sqrt{2}\int_s^te^{-(t-r)\H}\zeta(\de r)
\end{equation}
for $s,t\ge0$ to use that the independence of time increments of $\zeta$ as well as Markov property, we use the same notation for $\z=(\z^{(1)},\z^{(2)},\z^{(3)})$. We begin this section with preleminary results that will be important in the following.

\medskip

\begin{lem}{\label{LEM:Shiftedv}}
For any $u_0\in\CC^{-\varepsilon}$, $t\geq0$ and $h>0$, we have
$$
u(t+h\pv u_0)=\llp_{t,t+h}+v_{t,t+h}
$$
where $v_{t,t+\cdot}$ solves a shifted version of equation \eqref{model} with $\y=\z_{t,t+\cdot}$ and $v_{t,t}=u(t\pv u_0)$, that is the mild formulation
$$
v_{t,t+h}=e^{-h\H}u(t\pv u_0)-\int_0^he^{-(h-r)\H}\sum_{j=0}^3v_{t,t+r}^{j}\bz_{t,t+r}^{(3-j)}\de r.
$$
\end{lem}

\medskip

\begin{proof}
The mild formulation of equation \eqref{model} with $\y=\z$ and initial data $u_0$ gives
\begin{align*}
v(t+h\pv u_0)&=e^{-(t+h)\H}u_0-\int_0^{t+h}e^{-(t+h-r)\H}\sum_{j=0}^3v_r^{j}\bz_r^{(3-j)}\de r\\
&=e^{-h\H}\Big(e^{-t\H}u_0-\int_0^{t}e^{-(t-r)\H}\sum_{j=0}^3v_r^{j}\bz_r^{(3-j)}\de r\Big)-\int_t^{t+h}e^{-(t+h-r)\H}\sum_{j=0}^3v_r^{j}\bz_r^{(3-j)}\de r\\
&=e^{-h\H}\big(u(t\pv u_0)-\llp(t)\big)-\int_0^he^{-(h-r)\H}\sum_{j=0}^3v_{t+r}^{j}\bz_{t+r}^{(3-j)}\de r
\end{align*}
hence
\begin{align}
u(t+h\pv u_0)=\llp(t+h)-e^{-h\H}\llp(t)+e^{-h\H}u(t\pv u_0)-\int_0^he^{-(h-r)\H}\sum_{j=0}^3v_{t+r}^{j}\bz_{t+r}^{(3-j)}\de r
\end{align}
recalling that $u=\llp+v$. We have
\begin{align}
\llp(t+h)-e^{-h\H}\llp(t)&=\sqrt{2}\int_0^{t+h}e^{-(t+h-r)\H}\zeta(\de r)-\sqrt{2}\int_0^te^{-(t+h-r)\H}\zeta(\de r)\\ 
&=\sqrt{2}\int_t^{t+h}e^{-(t+h-r)\H}\zeta(\de r)\\ 
&=\llp_{t,t+h}
\end{align}
hence the expression
\begin{equation}
v_{t,t+h}=e^{-h\H}u(t\pv u_0)-\int_0^he^{-(h-r)\H}\sum_{j=0}^3v_{t+r}^{j}\bz_{t+r}^{(3-j)}\de r.
\end{equation}
From the binomial expansion Lemma \ref{Lem:Wickllp} and the relation between $v_t$ and $v_{t,t+h}$, we get
$$
\sum_{j=0}^3v_{t+r}^{j}\bz_{t+r}^{(3-j)}=\sum_{j=0}^3v_{t,t+r}^{j}\bz_{t,t+r}^{(3-j)}
$$
and the expected expression of $v_{t,t+h}$ follows.
\end{proof}

\begin{cor}{\label{BoundMoments}}
Let $u(\cdot\pv u_0)=\llp + v(\cdot\pv u_0)$ be the solution to equation \eqref{AndersonSQE}. Then for any $\varepsilon>0$ small enough and $p$ large enough, we have
$$
\sup_{u_0\in\CC^{-\varepsilon}}\sup_{t\geq0}\big(1\wedge t^{\frac{3p-2}{2}}\big)\E\left[\|u(t\pv u_0)\|_{\CC^{-\varepsilon}}^{3p-2}\right]<+\infty.
$$
\end{cor}

\medskip

\begin{proof}
First note that
$$
\sup_{u_0\in\CC^{-\varepsilon}}\sup_{0\leq t\leq 1}\big(1\wedge t^{\frac{3p-2}{2}}\big)\E\left[\|u(t\pv u_0)\|_{\CC^{-\varepsilon}}^{3p-2}\right]<+\infty
$$
is a direct consequence of Proposition \ref{ComingDown}. The bound when the supremum ranges over all positive times does not follow directly however as we cannot bound the expectation of $\sup_{t\geq0}\tilde{K}_t$. We are able to circumvent this issue using that the noise data $\bz$ has the same law on intervals of the same size and the bound in Proposition \ref{ComingDown} does not depend on the initial condition. Let $t\ge1$ and observe that for $u_0\in\CC^{-\varepsilon}(\T^2)$, we have
\begin{equation}
u(t\pv u_0)=\llp_{t-1,t}+v_{t-1,t}
\end{equation} 
where $v_{t-1,t}$ is given by Lemma \ref{LEM:Shiftedv}. In particular, $v_{t-1,t-1+\cdot}$ solves the equation \eqref{model} with noise $\y=\bz_{t-1,t-1+\cdot}$ and initial data $u(t-1\pv u_0)$. Proposition \ref{ComingDown} with the continuity of $\Gamma^{-1}$ gives
\begin{equation}
\|v_{t-1,t}\|_{L^{3p-2}}\lesssim 1\wedge c\tilde{K}_{t-1,t}
\end{equation} 
with
$$
\tilde{K}_{t-1,t}=\big(\sup_{0\leq r\leq 1}\|\bz_{t-1,t-1+r}^{(1)}\|_{\CC^{-\varepsilon}} + \|\bz_{t-1,t-1+r}^{(2)}\|_{L_r^{{q}}\big([0,1],\CC^{-\varepsilon}\big)} + \|\bz_{t-1,t-1+r}^{(3)}\|_{L_r^{{q}}\big([0,1],\CC^{-\varepsilon}\big)}\big)^{b}.
$$
From Lemma \ref{Lem:llp}, the law of $\bz_{t-1,t-1+\cdot}$ does not depend on $t$, one can take the expectation and get
\begin{align*}
\E\left[\|u(t\pv u_0)\|_{\CC^{-\varepsilon}}^{3p-2}\right]&\lesssim \E\left[\|\llp_{t-1,t}\|_{\CC^{-\varepsilon}}^{3p-2}\right]+\E\left[\|\tilde{v}_{t-1,t}\|_{\CC^{-\varepsilon}}^{3p-2}\right]\\
&\lesssim \E\left[\|\llp_{t-1,t}\|_{\CC^{-\varepsilon}}^{3p-2}\right]+\E\left[\|\tilde{v}_{t-1,t}\|_{L^{3p-2}}^{3p-2}\right]
\end{align*}
where the right-hand side is bounded uniformly both in $t$ and $u_0$. This gives
$$
\sup_{u_0\in\CC^{\varepsilon}}\sup_{t\leq r\leq t+1}\big(1\wedge r^{\frac{3p-2}{2}}\big)\E\left[\|u(r\pv u_0)\|_{\CC^{-\varepsilon}}^{3p-2}\right]<+\infty
$$
which completes the proof since $t\ge1$ is arbitrary.
\end{proof}

\subsection{The transition semi-group and Feller property}

As we intend to prove Markov property for $t\mapsto u(t\pv \cdot)$, we give a more precise description of $\tilde{v}$ defined in Lemma \ref{LEM:Shiftedv}. For $h\geq0$, we write $v(h\pv v_0\pv \mathbf{y})$ the solution to \eqref{model} at time $h$, starting from $v_0\in\CC^{-\eps}$ and driven by noise data $\mathbf{y}$. This notation will simplify the proof of the following proposition. It proves that $t\mapsto P_t$ is indeed a semigroup which gives the Markov property as direct consequence.

\medskip

\begin{prop}{\label{Prop:Markov}}
Let $u_0\in\CC^{-\varepsilon}$ and $t\ge0$. Then for any $\varphi\in C_b\big(\CC^{-\varepsilon}\big)$, we have
$$
\E\Big[\varphi\big(u(t+h\pv u_0)\big)\big|\mathcal{F}_t\Big]=P_h\varphi\big(u(t\pv u_0)\big)
$$
for $h\geq0$.
\end{prop}

\medskip

\begin{proof}
In view of Lemma \ref{LEM:Shiftedv}, we have
$$
u(t+h\pv u_0)=\llp_{t,t+h}+v\big(h;u(t\pv u_0)\pv\bz_{t,t+\cdot}\big)
$$
where $u(t\pv u_0)$ is $\mathcal{F}_t$-measurable while $\bz_{t,t+\cdot}$ is $\mathcal{F}_t$-independent, see Lemma \ref{Lem:llp}. Thus, according to Proposition 1.12 in \cite{DZ14}, we have
\begin{align*}
\E\left[\varphi\big(u(t+h\pv u_0)\big)\middle|\mathcal{F}_t\right]&=\E\left[\varphi\big(\llp_{t,t+h}+v\big(h\pv u(t\pv u_0)\pv \bz_{t,t+\cdot}\big)\big)\middle|\mathcal{F}_t\right]\\
&=\E\left[\varphi\big(\llp_{t,t+h}+v\big(h\pv u(t\pv u_0)\pv \bz_{t,t+\cdot}\big)\big)\right]
\end{align*}
where the expectation is with respect to the data $\bz_{t,t+\cdot}$. Moreover, since $\bz_{t,t+\cdot}$ and $\bz$ have the same law with again Lemma \ref{Lem:llp}, we get
\begin{align*}
\E\left[\varphi\big(\llp_{t,t+h}+v\big(h\pv v_0\pv \bz_{t,t+\cdot}\big)\big)\right]&= \E\left[\varphi\big(\llp_{0,h}+v\big(h\pv v_0\pv \bz\big)\big)\right]\\
&=\E\left[\varphi\big(u(h\pv v_0)\big)\right]\\
&=P_h\varphi\big(v_0\big)
\end{align*}
for any $v_0\in\CC^{-\varepsilon}$. Applying this in $v_0=u(t\pv u_0)$ completes the proof.
\end{proof}

The Feller property for the semigroup $(P_t)_{t\geq0}$ follows as a direct consequence of Proposition \ref{Localwpv}. Indeed, not only do we obtain local well-posedness, but also that the solution $v$ to \eqref{model} depends continuously on the initial data in $\CC^{-\varepsilon}$. In the strong Feller property, one only assumes the functionnal to be a Borelian bounded function on $\CC^{-\eps}$, this will be proved in the rest of the section.

\medskip

\begin{cor}{\label{COR:Feller}}
The semigroup $(P_t)_{t\geq0}$ has the Feller property, that is
\begin{equation}
\forall\varphi\in C_b\big(\CC^{-\varepsilon}\big),P_t\varphi\in C_b\big(\CC^{-\varepsilon}\big)
\end{equation}
for any $t\ge0$.
\end{cor}

\medskip

We conclude this subsection with the existence of invariant measures for the semigroup $(P_t)_{t\geq0}$. We prove an ergodic theorem for the dual semigroup $(P^*_t)_{t\geq0}$ acting on the set of probability measures on $\CC^{-\varepsilon}$. The proof relies on a compactness argument hence uniqueness of an invariant measure is not ensured at this stage.

\medskip

\begin{prop}{\label{PROP:ExistenceInvariant}}
Let $u_0\in\CC^{-\varepsilon}(\T^2)$. There is a increasing sequence of diverging times $(t_k)_{k\ge0}$ and a probability measure $\mu_{u_0}$ on $\CC^{-\varepsilon}(\T^2)$ such that
$$
\lim_{k\to\infty}\frac{1}{t_k}\int_0^{t_k}P_r^*\delta_{u_0}\de r=\mu_{u_0}
$$
where the convergence holds weakly. As such, $\mu_{u_0}$ is an invariant measure for $(P_t)_{t\geq0}$ acting on $C_b\big(\CC^{-\varepsilon}\big)$.
\end{prop}

\medskip

\begin{proof}
We intend to use Krylov-Bogoliubov theorem as given in theorem 3.1.1 of \cite{DZ96}. Set
$$
R_t:=\frac{1}{t}\int_0^{t}P_r^*\delta_{u_0}\de r
$$ 
for $t>0$, it is enough to prove that the sequence $(R_{t_k})_{k\in\N}$ is uniformly tight as a sequence of probability measures on $\CC^{-\varepsilon}$ for an increasing sequence of diverging times $(t_k)_{k\ge0}$. In particular, we are not bothered by the behavior of $R_t$ close to $0$ since only when large times matter. Let $p>1$ be as in Proposition \ref{BoundMoments}, $t\geq0$ and $R>0$. We begin with the following observation that arises from the definition of $P^*$. Using Markov and Jensen's inequality, we have
\begin{align*}
R_t\left(\{u_0\in\CC^{-\varepsilon_0}\,,\ \|u_0\|_{\CC^{-\varepsilon_0}}>R\}\right)&=\frac{1}{t}\int_0^t\mathbb{P}\left(\|u(r\pv u_0)\|_{\CC^{-\varepsilon_0}}>R\right)\de r\\
&\lesssim \frac{1}{Rt}\int_0^t\E\left[\|u(r\pv u_0)\|_{\CC^{-\varepsilon_0}}^{3p-2}\right]^{\frac{1}{3p-2}}\de r.
\end{align*}
Then Proposition \ref{BoundMoments} yields
\begin{align*}
R_t\left(\{u_0\in\CC^{-\varepsilon_0}\,,\ \|u_0\|_{\CC^{-\varepsilon_0}}>R\}\right)\lesssim \frac{1+t}{Rt}.
\end{align*}
Given any sequence of increasing diverging times $(t_k)_{k\ge0}$, we have
\begin{equation}
\sup_{k\ge0}\frac{1+t_k}{Rt_k}<\infty
\end{equation}
while the set $\{u_0\in\CC^{-\varepsilon_0}\,,\ \|u_0\|_{\CC^{-\varepsilon_0}}\leq R\}$ is a compact subset of $\CC^{-\varepsilon}$ for any $\varepsilon>\varepsilon_0$ hence $(R_{t_k})_{k\in\N}$ is tight. As such, there exists a limiting measure $\mu_{u_0}$ on $\CC^{-\varepsilon}$ such that $R_{t_k}$ converges weakly to $\mu_{u_0}$ up to extraction. Theorem 3.1.1 in \cite{DZ96} then ensures that $\mu_{u_0}$ is an invariant measure to the semigroup $(P_t)_{t\geq0}$ acting on $C_b\big(\CC^{-\varepsilon}\big)$.
\end{proof}

\subsection{Investigating the strong Feller property}

As far as the strong Feller property is concerned, the idea is to obtain a Lipschitz bound on the semigroup $P_t$ which will ensure that continuity is recovered even for merely borelian bounded test functions $\varphi$. This leads to investigate the behavior of $u$ as a function of its initial condition $u_0$. For the computations to make sense, it is necessary to work at the level of the truncated equation so that the Wick power is an actual polynomial rather than a renormalized object. Thus, let $P_t^N$ be the transition semigroup associated to equation \eqref{TruncEquation}
\begin{equation}
\partial_t u_N + \H u_N + \Wick{u_N^3} = \sqrt{2}\Pi_N\zeta
\end{equation}
defined as
\begin{equation}
P_t^N(\varphi)(u_0)=\E\big[\varphi(u_N(t\pv u_0)\big]
\end{equation} 
for $u_0\in\CC^{-\eps}$. The main result of this section is the following local Lipschitz bound on the differential of the truncated semigroup which gives the strong Feller property with Corollary \ref{COR:StrongFeller}.

\medskip

\begin{prop}{\label{PROP:DiffPt}}
Let $\varepsilon,\sigma,\delta>0$ and $p\geq 1$ such that \ref{Gronwall_eta}, \ref{COR:GlobalEstimatev} and \ref{COR:SquareSobolevuN} hold true. Assume also that $\delta+p\sigma+(p-1)\varepsilon + (\sigma-\varepsilon) <1$. Then there exists some constant $\gamma>0$ such that for any $N\geq1$, $\varphi\in C_b\left(\CC^{-\varepsilon}\right)$, $u_0\in \CC^{-\varepsilon}$ and $h\in L^2\cap \CC^{-\delta}$, we have
$$
\left|\de(P_t^N\varphi)(u_0)\cdot h\right|\lesssim \|\varphi\|_{L^\infty}\|h\|_{\CC^{-\delta}}\big(1+t^{-\frac{1+(p+1)\sigma+p\varepsilon+\delta}{2}}\big)(1+\|u_0\|_{\CC^{-\varepsilon}})^{\gamma}
$$
for any $t>0$.
\end{prop}

\medskip

We first investigate the behavior of $u$ under variation of the initial condition. Recall equation \ref{TruncEquation}, if we set 
\begin{align}{\label{DEF_eta}}
\eta^h_N =\de(u_N)(u_0)\cdot h
\end{align}
for $u_0,h\in\CC^{-\varepsilon}$, then $\eta^h_N$ solves the linear equation
\begin{align}{\label{EQ_eta}}
\begin{cases}
\partial_t \eta^h_N + \H \eta^h_N + 3\Wick{u_N^2}\,\eta^h_N =0\\
{\eta^h_N}_{|t=0}=h.
\end{cases}   
\end{align}
The following lemma is what calls for the use of an alternative semi group in order to kill the blow-up obtained using Grönwall's estimate following \cite{DD20}.

\medskip

\begin{lem}{\label{Gronwall_eta}}
Let $\sigma>\varepsilon$, $\delta>0$ small enough and $p\geq1$ large enough satisfying $2\sigma+\varepsilon+\delta+1<\frac{2(p-1)}{p}$. Then there exists a constant $c>0$, uniform in $t\in(0,1]$, $N\in\N$, $u_0\in\CC^{-\varepsilon}$ and $h\in\CC^{-\delta}$ such that
$$
\|\eta^h_N(t)\|_{\CC^{\sigma}}\lesssim t^{-(\sigma+\delta)/2}\|h\|_{\CC^{-\delta}}e^{c\int_0^t\|\Wick{u_N^2(s)}\|_{H^{-\varepsilon}}^p\de s}.
$$
\end{lem}

\medskip

Note that $\eta^h_N$ has positive spatial regularity as soon as $t>0$, this is due to the regularizing effects of $e^{-t\H}$ with Schauder estimates, hence the diverging factor as $t$ goes to $0$. Also, a more direct bound would yield a norm $\CC^{-\varepsilon}$ in the exponential as well as a lighter condition on the exponents however, as we want to use a Bismuth-Elworthy-Li formula in the following, it is more convenient to work with a Hilbert space norm of $u_N$ for differentiation matters.

\medskip

\begin{proof}
Write the mild formulation of equation \eqref{EQ_eta} and take the $\CC^{\sigma}$ norm, this yields
\begin{align*}
\|\eta^h_N(t)\|_{\CC^{\sigma}}&\lesssim \|e^{-t\H}h\|_{\CC^{\sigma}}+\int_0^t\left\|e^{-(t-s)\H}\big(\Wick{u_N^2(s)}\eta^h_N(s)\big)\right\|_{\CC^{\sigma}}\de s\\
&\lesssim t^{-(\sigma+\delta)/2}\|h\|_{\CC^{-\delta}}+\int_0^t(t-s)^{-(\sigma+\varepsilon)/2}\left\|e^{-\frac{t-s}{2}\H}\big(\Wick{u_N^2(s)}\eta^h_N(s)\big)\right\|_{\CC^{-\varepsilon}}\de s\\
&\lesssim t^{-(\sigma+\delta)/2}\|h\|_{\CC^{-\delta}}+\int_0^t(t-s)^{-(\sigma+\varepsilon)/2}\left\|e^{-\frac{t-s}{2}\H}\big(\Wick{u_N^2(s)}\eta^h_N(s)\big)\right\|_{H^{1-\varepsilon}}\de s\\
&\lesssim t^{-(\sigma+\delta)/2}\|h\|_{\CC^{-\delta}}+\int_0^t(t-s)^{-(\sigma+\varepsilon+1)/2}\left\|\Wick{u_N^2(s)}\right\|_{H^{-\varepsilon}}\left\|\eta^h_N(s)\right\|_{\CC^{\sigma}}\de s
\end{align*}
using Schauder estimates and Besov embeddings. Setting $\lambda(t):=t^{(\sigma+\delta)/2}\|\eta^h_N(t)\|_{\CC^{\sigma}}$, Hölder inequality yields
\begin{align*}
\lambda(t)&\lesssim\|h\|_{\CC^{-\delta}}+\int_0^t(t-s)^{-(\sigma+\varepsilon+1)/2}s^{-(\sigma+\delta)/2}\left\|\Wick{u_N^2(s)}\right\|_{H^{-\varepsilon}}\lambda(s)\de s\\
&\lesssim \|h\|_{\CC^{-\delta}}+\left(\int_0^t(t-s)^{-\frac{p(\sigma+\varepsilon+1)}{2(p-1)}}s^{-\frac{p(\sigma+\delta)}{2(p-1)}}\de s\right)^\frac{p-1}{p}\left(\int_0^t\left\|\Wick{u_N^2(s)}\right\|_{H^{-\varepsilon}}^p\lambda(s)^p\de s\right)^\frac{1}{p}
\end{align*}
for any $p\geq1$. Choosing $p$ such that $2\sigma+\varepsilon+\delta+1<\frac{2(p-1)}{p}$, the first integral is bounded uniformly in $t\in(0,1]$ so that Grönwall inequality on $\lambda(t)^p$ yields the bound
$$
\|\eta^h_N(t)\|_{\CC^{\sigma}}\lesssim t^{-(\sigma+\delta)/2}\|h\|_{\CC^{-\delta}}e^{c\int_0^t\|\Wick{u_N^2(s)}\|_{H^{-\varepsilon}}^p\de s}.
$$
\end{proof}

\subsection{Estimates on the auxiliary semigroup}

Because of the exponential on the right-hand side of Lemma \ref{Gronwall_eta}, we will rather consider the Feynman-Kac semigroup $S^N_t$ defined by
\begin{align}{\label{DEF:S_N}}
S_t^N\varphi(u_0):=\E\Big[e^{-\int_0^tV\left(u_N(s\pv u_0)\right)\de s}\varphi\big(u_N(s\pv u_0)\big)\Big]    
\end{align}
for a well-chosen non-negative functional $V$ of the solution $u_N$. In particular, the Kolmogorov equation for $S_t^N$ is the same as for $P_t^N$ with an additional shift by the potential. At the level of the mild formulation, this reads
\begin{align}{\label{relationPtSt}}
P_t^N\varphi = S_t^N\varphi + \int_0^t S_{t-s}^N\big(VP_s^N\varphi\big)\de s.
\end{align}
Thus, while we hope to obtain better estimats on $S_t^N$ thanks to the exponential factor, this will be sufficient for establishing similar bounds on $P_t^N$. The method is taken from the recent work \cite{DD20} by Da Prato and Debussche on the usual $\Phi_2^4$ equation on the torus. We refer to Da Prato and Zabczyk's book \cite{DZ97} and references therein for more details on the Feynman Kac semigroups. In view of Lemma \ref{Gronwall_eta}, set 
$$
V(u_N):=\tilde{c}\|\Wick{u_N^2}\|_{H^{-\varepsilon}}^p
$$
for a positive constant $\tilde{c}$ to be chosen later, thus following moments bound on $V$ holds. This is where it better to work in an Hilbert space since we will differentiate $V$.

\medskip

\begin{cor}{\label{COR:SquareSobolevuN}}
Let $N\geq1$ and $\sigma>\varepsilon>0$ be small enough so that \ref{COR:Globalv}, \ref{COR:GlobalEstimatev}, \ref{PROP:BoundTildev}. For $u_0\in\CC^{-\varepsilon}$ and $p\geq 1$, there exist a constant $\gamma>0$ depending on $\varepsilon$ and $p$ such that
$$
\E\Big[\|\Wick{u_N(t)^2}\|_{H^{-\varepsilon}}^p\Big]\lesssim t^{-\frac{\sigma+\varepsilon}{2}p}(1+\|u_0\|_{\CC^{-\varepsilon}})^\gamma
$$
uniformly in $t\in(0,1)$ and $N\geq1$.
\end{cor}

\medskip

\begin{proof}
Recall that $u_N=\Pi_N\llp+v_N$ as well as 
$$
\tilde{v}_N=v_N - e^{-t\H}u_0
$$
which has better regularity estimate given by Lemma \ref{LEM:EstimateL2tildev}. Using the product rule on Besov space, we get
\begin{align*}
\|\Wick{u_N(t)^2}\|_{H^{-\varepsilon}}&\lesssim \|v_N(t)^2\|_{H^{-\varepsilon}} + \|v_N(t)\Pi_N\llp(t)\|_{H^{-\varepsilon}} + \|\Wick{(\Pi_N\llp(t))^2}\|_{H^{-\varepsilon}}\\
&\lesssim \|v_N(t)\|_{H^\sigma}\|v_N(t)\|_{\CC^{-\varepsilon}} + \|v_N(t)\|_{H^\sigma}\|\Pi_N\llp(t)\|_{\CC^{-\varepsilon}} + \|\Wick{(\Pi_N\llp(t))^2}\|_{\CC^{-\varepsilon}}\\
&\lesssim \big(t^{-\frac{\sigma+\varepsilon}{2}}\|u_0\|_{\CC^{-\varepsilon}}+\|\tilde{v}_N(t)\|_{H^\sigma}\big)\big(\|v_N(t)\|_{\CC^{-\varepsilon}}+\|\Pi_N\llp(t)\|_{\CC^{-\varepsilon}}\big) + \|\Wick{(\Pi_N\llp(t))^2}\|_{\CC^{-\varepsilon}}\\
&\lesssim \|\Pi_N\llp(t)\|_{\CC^{-\varepsilon}}^2 + \|\Wick{(\Pi_N\llp(t))^2}\|_{\CC^{-\varepsilon}} + \|\tilde{v}_N(t)\|_{H^\sigma}^2 + \|v_N(t)\|_{\CC^{-\varepsilon}}^2\\
&\qquad + t^{-\frac{\sigma+\varepsilon}{2}}\|u_0\|_{\CC^{-\varepsilon}}\big(\|v_N(t)\|_{\CC^{-\varepsilon}}+\|\Pi_N\llp(t)\|_{\CC^{-\varepsilon}}\big).
\end{align*}
Taking the $p$-th moment, we can conclude using Corollary \ref{COR:GlobalEstimatev} for the Hölder norm and Proposition \ref{PROP:BoundTildev} for the Sobolev norm.
\end{proof}

Let us recall the following Bismuth-Elworthy-Li formula from \cite{DZ97}.

\medskip

\begin{prop}{\label{PROP:BEL}}
For any $\varphi\in C_b\left(\CC^{-\varepsilon}\right)$, the semigroup $S_T^N\varphi$ is differentiable in any direction $h\in L^2$ and we have
\begin{align*}
\de(S_t^N\varphi)(u_0)\cdot h=\E\Big[&e^{-\int_0^tV(u_N(s\pv u_0))\de s}\varphi\big(u_N(t\pv u_0)\big)\times\\ 
&\Big(\frac{1}{t}\int_0^t\eta_N^h(s)\zeta(\de s)+\int_0^t(1-\frac{s}{t})\de V(u_N(s\pv u_0))\cdot\eta^h_N(s)\de s\Big)\Big].
\end{align*}  
\end{prop}

\medskip

This is precisely where working at the level of the approximated solution $u_N$ rather than the limiting one $u$ is needed. To the best of our knowledge, for the identity above to hold it is required that the Feynman Kac potential is a Lipschitz function of the solution. This is however not an issue as all the estimates are uniform in $N\ge1$. Abusing notations, we denote by $\|\varphi\|_{L^\infty}$ the infimum of $\{M\geq 0\,,\ \forall u_0\in\CC^{-\varepsilon}\,,\ |\varphi(u_0)|\leq M\}$ for $\varphi\in\mathcal{B}_b$ the space of Borelian bounded functions on $\CC^{-\eps}$.

\medskip

\begin{prop}{\label{PROP:Estimate_S^Nvarphi}}
Let $\varepsilon,\sigma,\delta>0$ and $p\geq 1$ such that \ref{Gronwall_eta}, \ref{COR:GlobalEstimatev} and \ref{COR:SquareSobolevuN} hold true. Assume also that $\delta+p\sigma+(p-1)\varepsilon<2$, then there exists a constant $\gamma>0$ such that for any $N\geq1$, $\varphi\in C_b\left(\CC^{-\varepsilon}\right)$, $u_0\in \CC^{-\varepsilon}$ and $h\in L^2\cap \CC^{-\delta}$, we have
$$
\left|\de(S_t^N\varphi)(u_0)\cdot h\right|\leq C\|\varphi\|_{L^\infty}\|h\|_{\CC^{-\delta}}\left(t^{-\frac{1+\sigma+\delta}{2}}+(1+\|u_0\|_{\CC^{-\varepsilon}})^\gamma\right).
$$
\end{prop}

\medskip

\begin{proof}
Let $u_0\in\CC^{-\varepsilon}$ and $h\in L^2\cap\CC^{-\delta}$ as well as $\sigma, \delta, p$ be as in Lemma \ref{Gronwall_eta}. Using Proposition \ref{PROP:BEL}, we  bound 
$$
\mathbf{A}:=\E\Big[e^{-\int_0^tV(u_N(s\pv u_0))\de s}\varphi\big(u_N(t\pv u_0)\big)\frac{1}{t}\int_0^t\eta_N^h(s)\zeta(\de s)\Big]
$$
and
$$
\mathbf{B}:=\E\Big[e^{-\int_0^tV(u_N(s\pv u_0))\de s}\varphi\big(u_N(t;u_0)\big)\int_0^t\big(1-\frac{s}{t}\big)\de V\big(u_N(s\pv u_0)\big)\cdot\eta_N^h(s)\de s\Big]
$$
separately. For the first term, we have
\begin{align*}
\mathbf{A}&\leq \frac{1}{t}\|\varphi\|_{L^\infty}\E\Big[\Big|e^{-\int_0^tV(u_N(s\pv u_0))\de s}\int_0^t\eta_N^h(s)\zeta(\de s)\Big|^2\Big]^\frac{1}{2}\\
&\leq \frac{1}{t}\|\varphi\|_{L^\infty}\E\Big[\int_0^te^{-2\int_0^sV(u_N(r\pv u_0))\de r}\big|\eta_N^h(s)\big|^2\de s\Big]^\frac{1}{2}
\end{align*}
Cauchy-Schwarz inequality and Itô formula. Choosing the constant $\tilde{c}$ in the definition of $V$ large enough, Lemma \ref{Gronwall_eta} yields
\begin{align*}
\E\left[\int_0^te^{-2\int_0^sV(u_N(r\pv u_0))\de r}\left|\eta_N^h(s)\right|^2\de s\right]^\frac{1}{2}&\lesssim \|h\|_{\CC^{-\delta}}\left(\int_0^ts^{\sigma+\delta}\de s\right)^{\frac{1}{2}}\\
&\lesssim t^{\frac{1-\sigma-\delta}{2}}\|h\|_{\CC^{-\delta}}   
\end{align*}
so that
$$
\mathbf{A}\lesssim t^{-\frac{1+\sigma+\delta}{2}}\|\varphi\|_{L^\infty}\|h\|_{\CC^{-\delta}}.
$$
For the second term, we have
$$
\mathbf{B}\leq \|\varphi\|_{L^\infty}\E\left[e^{-2\int_0^tV(u_N(s\pv u_0))\de s}\left|\int_0^t\big(1-\frac{s}{t}\big)\de V\big(u_N(s\pv u_0)\big)\eta_N^h(s)\de s\right|^2\right]^{\frac{1}{2}}
$$
using Cauchy-Schwarz inequality once again. Recall that 
$$
V(u_N(s\pv u_0))=\tilde{c}\|\Wick{u_N(s;u_0)^2}\|_{H^{-\varepsilon}}^p
$$ 
so that we can easily differentiate $V$, that is
\begin{align*}
\de V\big(u_N(s\pv u_0)\big)\cdot\eta_N^h(s)=2p\|\Wick{u_N(s\pv u_0)^2}\|_{H^{-\varepsilon}}^{p-2}\left\langle\Wick{u_N(s\pv u_0)^2},u_N(s\pv u_0)\eta_N^h(s)\right\rangle_{H^{-\varepsilon}}
\end{align*}
and therefore
$$
\mathbf{B}\leq \|\varphi\|_{L^\infty}\E\left[e^{-2\int_0^tV(u_N(s\pv u_0))\de s}\left(\int_0^t\|\Wick{u_N(s\pv u_0)^2}\|_{H^{-\varepsilon}}^{p-1}\|u_N(s\pv u_0)\|_{H^{-\varepsilon}}\|\eta_N^h(s)\|_{\CC^{}}\de s\right)^2\right]^{\frac{1}{2}}.
$$
Using once again Lemma \ref{Gronwall_eta} to handle the $\eta_N^h$ term as well as Minkowski and Cauchy-Schwarz inequalities we get
\begin{align*}
\mathbf{B}&\lesssim \|\varphi\|_{L^\infty}\|h\|_{\CC^{-\delta}}\E\Big[\Big(\int_0^t\|\Wick{u_N(s\pv u_0)^2}\|_{H^{-\varepsilon}}^{p-1}\|u_N(s\pv u_0)\|_{H^{-\varepsilon}}s^{-\frac{\sigma+\delta}{2}}\de s\Big)^2\Big]^{\frac{1}{2}}\\
&\lesssim \|\varphi\|_{L^\infty}\|h\|_{\CC^{-\delta}}\int_0^t\E\Big[\|\Wick{u_N(s\pv u_0)^2}\|_{H^{-\varepsilon}}^{2p-2}\|u_N(s\pv u_0)\|^2_{H^{-\varepsilon}}\Big]^{\frac{1}{2}}s^{-\frac{\sigma+\delta}{2}}\de s\\
&\lesssim \|\varphi\|_{L^\infty}\|h\|_{\CC^{-\delta}}\int_0^t\E\Big[\|\Wick{u_N(s\pv u_0)^2}\|_{H^{-\varepsilon}}^{4p-4}\Big]^{\frac{1}{4}}\E\Big[\|u_N(s\pv u_0)\|^4_{H^{-\varepsilon}}\Big]^{\frac{1}{4}}s^{-\frac{\sigma+\delta}{2}}\de s.
\end{align*}
Using the moments bounds \ref{COR:GlobalEstimatev} and \ref{COR:SquareSobolevuN}, this yields
$$
\mathbf{B}\lesssim \|\varphi\|_{L^\infty}\|h\|_{\CC^{-\delta}}(1+\|u_0\|_{\CC^{-\varepsilon}})^\gamma\int_0^ts^{-\frac{\sigma+\delta}{2}-(p-1)\frac{\sigma+\varepsilon}{2}}\de s
$$
and provided $\delta+p\sigma+(p-1)\varepsilon<2$, the time integral is bounded uniformly in $[0,1]$ and the estimate follows.
\end{proof}

Note that in view of the definition of $S_t^N$ \eqref{DEF:S_N} and the bound \ref{COR:SquareSobolevuN}, $S_t^N$ can be extended from bounded borelian functions $\varphi$ on $\CC^{-\varepsilon}$ to functions under the form $V\varphi$, hence the following counterpart to the previous estimate.

\medskip

\begin{prop}{\label{PROP:Estimate_S^NVvarphi}}
Let $\varepsilon,\sigma,\delta>0$ and $p\geq 1$ such that \ref{Gronwall_eta}, \ref{COR:GlobalEstimatev} and \ref{COR:SquareSobolevuN} hold true. Assume also that $\delta+p\sigma+(p-1)\varepsilon<2$, then there exists some a $\gamma>0$ such that for any $N\geq1$, $\varphi\in C_b\left(\CC^{-\varepsilon}\right)$, $u_0\in \CC^{-\varepsilon}$ and $h\in L^2\cap \CC^{-\delta}$, we have
$$
\left|\de(S_t^N(V\varphi))(u_0)\cdot h\right|\leq C\|\varphi\|_{L^\infty}\|h\|_{\CC^{-\delta}}t^{-\frac{1+(p+1)\sigma+p\varepsilon+\delta}{2}}(1+\|u_0\|_{\CC^{-\varepsilon}})^{\gamma}.
$$
\end{prop}

\medskip

\begin{proof}
The proof follows the exact same lines as the one above, up to the multiplicative term that appears when using Cauchy-Schwarz inequality $\E\left[V\big(u_N(t\pv u_0)\big)^2\right]^{\frac{1}{2}}$ when estimating $\mathbf{A}$ and $\mathbf{B}$. This is handled by Corollary \ref{COR:SquareSobolevuN} and the estimate follows from 
\begin{align*}
\mathbf{A}+\mathbf{B}&\lesssim \|\varphi\|_{L^\infty}\|h\|_{\CC^{-\delta}}\left(t^{-\frac{1+\sigma+\delta}{2}}+(1+\|u_0\|_{\CC^{-\varepsilon}})^\gamma\right)\E\left[\|\Wick{u_N(t\pv u_0)^2}\|_{H^{-\varepsilon}}^2\right]^{\frac{1}{2}}\\
&\lesssim \|\varphi\|_{L^\infty}\|h\|_{\CC^{-\delta}}t^{-p\frac{\sigma+\varepsilon}{2}}\left(t^{-\frac{1+\sigma+\delta}{2}}+(1+\|u_0\|_{\CC^{-\varepsilon}})^\gamma\right)(1+\|u_0\|_{\CC^{-\varepsilon}})^{\gamma'}\\
&\lesssim \|\varphi\|_{L^\infty}\|h\|_{\CC^{-\delta}}t^{-\frac{1+(p+1)\sigma+p\varepsilon+\delta}{2}}(1+\|u_0\|_{\CC^{-\varepsilon}})^{\gamma''}
\end{align*}
for yet another exponent $\gamma''>0$.
\end{proof}

We can now exploit relation \eqref{relationPtSt} in order to obtain similar estimates on the transition semigroup $P_t^N$ and prove our main result. 

\medskip

\begin{proofof}\textbf{\sffamily of Proposition \ref{PROP:DiffPt} :}
Note that the assumption on the different parameters indeed allows to use the previous proposition. Thus, differentiating relation \eqref{relationPtSt} in the direction $h$ and injecting the estimates of Propositions \ref{PROP:Estimate_S^Nvarphi} and \ref{PROP:Estimate_S^NVvarphi}, we get
\begin{align*}
\left|\de(P_t^N\varphi)(u_0)\cdot h\right|&\leq \left|\de(S_t^N\varphi)(u_0)\cdot h\right| + \int_0^t\left|\de\big(S_{t-s}^N(VP_s^N\varphi)\big)(u_0)\cdot h\right|\de s\\
&\lesssim \|\varphi\|_{L^\infty}\|h\|_{\CC^{-\delta}}\left(t^{-\frac{1+\sigma+\delta}{2}}+(1+\|u_0\|_{\CC^{-\varepsilon}})^\gamma\right)\\
&\qquad +\int_0^t\|P_s\varphi\|_{L^\infty}(t-s)^{-\frac{1+(p+1)\sigma+p\varepsilon+\delta}{2}}\de s\|h\|_{\CC^{-\delta}}(1+\|u_0\|_{\CC^{-\varepsilon}})^{\gamma}\\
&\lesssim \|\varphi\|_{L^\infty}\|h\|_{\CC^{-\delta}}\left(t^{-\frac{1+\sigma+\delta}{2}}+(1+\|u_0\|_{\CC^{-\varepsilon}})^\gamma\right)\\
&\qquad +\int_0^t(t-s)^{-\frac{1+(p+1)\sigma+p\varepsilon+\delta}{2}}\de s\|\varphi\|_{L^\infty}\|h\|_{\CC^{-\delta}}(1+\|u_0\|_{\CC^{-\varepsilon}})^{\gamma}
\end{align*}
where in the last inequality we used
$$
\|P_t\varphi\|_{L^\infty}\leq \|\varphi\|_{L^\infty}
$$
which follows directly from the definition of the semigroup. The assumption on the indices ensures the time integral converges and is bounded uniformly in $t\in[0,1]$, hence the claim.
\end{proofof}

\begin{cor}{\label{COR:StrongFeller}}
Let $\varepsilon,\sigma,\delta>0$ and $p\geq 1$ such that \ref{PROP:DiffPt} holds true. Then there exists a constant $\gamma>0$ such that for any bounded borelian function $\varphi$ on $\CC^{-\varepsilon}$, we have
$$
\left|(P_t\varphi)(u_0)-(P_t\varphi)(\tilde{u}_0)\right|\lesssim \|\varphi\|_{L^\infty}\|u_0-\tilde{u}_0\|_{\CC^{-\delta}}\big(1+t^{-\frac{1+(p+1)\sigma+p\varepsilon+\delta}{2}}\big)(1+\|u_0\|_{\CC^{-\varepsilon}}+\|\tilde{u}_0\|_{\CC^{-\varepsilon}})^{\gamma}
$$
for any $t>0$, initial data $u_0, \tilde{u}_0\in \CC^{-\varepsilon}$ and direction $h\in\CC^{-\delta}$. In particular, the semigroup $(P_t)_{t\ge0}$ has the strong Feller property, that is
\begin{equation}
\forall\varphi\in \mathcal{B}_b\big(\CC^{-\varepsilon}\big),P_t\varphi\in C_b\big(\CC^{-\varepsilon}\big)
\end{equation}
for any $t>0$.
\end{cor}

\medskip

\begin{proof}
First, since $L^2\cap\CC^{-\delta}$ is dense in $\CC^{-\delta}$, Proposition \ref{PROP:DiffPt} ensures that $h\mapsto\de(P_t^N)(u_0)\cdot h$ extends uniquely as a bounded linear form on $\CC^{-\varepsilon}$ that satisfies the same bound. Moreover, using the mean value inequality, the estimate holds true for $P_t^N\varphi$ with the right-hand side being uniform in $N\geq1$. Letting $N\to\infty$, the convergence part of Proposition \ref{THM:LWP} ensures 
$$
P_t^N(\varphi)\to P_t\varphi
$$ 
pointwisely on $\CC^{-\varepsilon}$ hence the bound and the strong Feller property.
\end{proof}

\section{Ergodicity for Anderson \texorpdfstring{$\Phi^4_2$}{PHI42} model}\label{Sec:ergodicity}

In order to obtain exponential mixing from the strong Feller property, Tsatsoulis and Weber \cite{TW18} proved a support theorem for the enhanced noise which implies a support result for the solution, following arguments by Chouk and Friz \cite{CF18}. This argument relies strongly on a construction on the torus with the explicit eigenfunctions of the Laplacian. While they prove that $(0,a,0)$ is in the support of their enhanced noise for any $a\ge0$, we only proved that $(0,0,0)$ is in the support of $\z$ here, this is the content of Proposition \ref{Prop:Supportnoise}. We prove in the following that this is still enough to get exponential mixing and uniqueness of the invariant measure. Since the equation with $\z=0$ relaxes to $0$ exponentially fast, the following proposition states that this still holds if the enhanced noise is small in $L^\infty([0,T],\CC^{-\kappa})^3$ which happends with positive probability. In particular, the speed of relaxation is independent of the initial data as the coming down from infinity a priori bound is uniform in $u_0\in\CC^{-\kappa}$. See also Chapter $5$ in Da Prato and Zabczyk \cite{DZ96} for more details on this.

\medskip

\begin{prop}\label{Prop:ApproximateControlability}
For any $\eps>0$, there exists $T=T(\eps)>0$ of order $|\log(\eps)|$ such that
\begin{equation}
\|v(T\pv u_0)\|_{\CC^{1-\kappa}}\le\eps
\end{equation}
for any $\kappa>0$ and $u_0\in\CC^{-\kappa}$ with positive probability.
\end{prop}

\medskip

\begin{proof}
Let $\eps>0$ and $T>0$. We work conditionnaly to 
\begin{equation}
\|\z\|_{L^q([0,1],(\CC^{-\kappa})^3)}+\sup_{t\in[1,T]}\|\z_t\|_{(\CC^{-\kappa})^3}+\|\z\|_{L^q([0,T],\CC^{-\kappa})^3}\le\eps
\end{equation}
which is an event of positive probability with Proposition \ref{Prop:Supportnoise}. Recall that $v$ satisfies the equation
\begin{equation}
\partial_tv+\H v+v^3+v^2\z^{(1)}+v\z^{(2)}+\z^{(3)}=0
\end{equation}
with $v(0)=u_0\in\CC^{-\kappa}$. First, the coming down from infinity with Corollary \ref{COR:Globalv} gives the bound
\begin{equation}
\|v(1)\|_{\CC^{-\kappa}}\lesssim\|v(1)\|_{L^p}\lesssim1
\end{equation}
independant of $u_0$ with $p\ge1$ large enough. Since the noise is bounded, the well-posedness result from Proposition \ref{Localwpv} with $\y=\z$ gives
\begin{equation}
\sup_{t\in[2,T]}\|v(t)\|_{\CC^{\sigma}}\lesssim c
\end{equation}
for any $\sigma<1$ with a constant $c>0$. Testing the equation for $v$ against $v$ yields
\begin{align*}
\frac{1}{2}\partial_t\|v\|_{L^2}^2&=-\langle\H v,v\rangle-\|v\|_{L^4}^4-\langle v^3,\z^{(1)}\rangle-\langle v^2,\z^{(2)}\rangle-\langle v,\z^{(3)}\rangle\\ 
&\le-\lambda_0\|v\|_{L^2}^2+\|v\|_{\CC^\sigma}^3\|\z^{(1)}\|_{\CC^{-\kappa}}+\|v\|_{\CC^\sigma}^2\|\z^{(2)}\|_{\CC^{-\kappa}}+\|v\|_{\CC^\sigma}\|\z^{(3)}\|_{\CC^{-\kappa}}\\
&\le-\lambda_0\|v\|_{L^2}^2+c\eps
\end{align*}
which gives
\begin{equation}
\|v(t)\|_{L^2}^2\le \|v(2)\|_{L^2}^2e^{-2\lambda_0(t-2)}+\frac{c\eps}{\lambda_0}\big(1-e^{-2\lambda_0(t-2)}\big)
\end{equation}
for $t\ge2$ with a Gröwall type argument. Taking $T=T(\eps)>0$ large enough, we get
\begin{equation}
\|v(T)\|_{L^2}^2\le \eps\big(1+\frac{c}{\lambda_0}\big)
\end{equation}
which is arbitrary small as $\eps>0$ goes to zero. This completes the proof since
\begin{align}
\|v(T)\|_{C^{1-\kappa+\kappa\theta-2\theta}}&\lesssim\|v(T)\|_{B_{\frac{2}{\theta},\frac{2}{\theta}}^{(1-\kappa)(1-\theta)}}\\
&\lesssim\|v(T)\|_{L^2}^\theta\|v(T)\|_{\CC^{1-\kappa}}^{1-\theta}\\
&\lesssim \eps^{\frac{\theta}{2}}
\end{align}
for any $\theta\in(0,1)$ using interpolation and embbeding of Besov spaces.
\end{proof}

\begin{rem}
In the work \cite{TW18}, Tsatsoulis and Weber consider general odd power $n$ except for the support theorem which they only show in the case $n=3$. This is because they only prove the support theorem for the enhanced noise in the cubic case while they expect the result to be true for general power, see their remark $6.2$. Our approach covers their case as we only need that the support of the enhanced noise contains zero, which is straighforward to prove in their case.
\end{rem}

\medskip

\begin{prop}\label{Prop:lowerboundProba}
There exists $\lambda\in(0,1)$ and $t_\lambda>0$ such that
\begin{equation}
\|P_t^*\delta_{u_0}-P_t^*\delta_{\tilde u_0}\|_{\text{TV}}\le 1-\lambda
\end{equation}
for $t\ge t_\lambda$ and $u_0,\tilde u_0\in\CC^{-\kappa}$.
\end{prop}

\medskip

\begin{proof}
Let $\kappa_0>\kappa>0$ and consider
\begin{equation}
\mathcal{A}_R:=\{u_0\in\CC^{-\kappa_0}\ ;\ \|u_0\|_{\CC^{-\kappa}}\le R\}
\end{equation}
for any $R>0$ wich is a compact of $\CC^{-\kappa_0}$ since the embedding $\CC^{-\kappa}\hookrightarrow\CC^{-\kappa_0}$ is compact. By Corollary \ref{COR:StrongFeller}, for every $a>0$ there exists $r=r(a)>0$ such that
\begin{equation}
\|P_t^*\delta_{u_0}-P_t^*\delta_{\tilde u_0}\|_{\text{TV}}\le 1-a
\end{equation}
for any $u_0,\tilde u_0\in B_{\CC^{-\kappa}}(0,r)$ the ball of radius $r$ and center $0$ in $\CC^{-\kappa}$. Proposition \ref{Prop:ApproximateControlability} gives the existence of $T=T(a)>0$ such that
\begin{equation}
(P_T\IDC_{B_{\CC^{-\kappa}}(0,r)})(u_0)>0
\end{equation}
for any $u_0\in\mathcal{A}_R$. Since $\IDC_{B_{\CC^{-\kappa}}(0,r)}$ is a measurable function on $\CC^{-\kappa_0}$, the strong Feller property implies that $u_0\mapsto(P_T\IDC_{B_{\CC^{-\kappa}}(0,r)})(u_0)$ is a continuous function hence
\begin{equation}
\inf_{u_0\in\mathcal{A}_R}(P_T\IDC_{B_{\CC^{-\kappa}}(0,r)})(u_0)\ge b>0
\end{equation}
for a positive constant $b=b(a)>0$ since $\mathcal{A}_R$ is compact. Markov property gives
\begin{align}
(P_t\Phi)(u_0)-(P_t\Phi)(\tilde u_0)&=\E\big[\big(P_{t-T}\phi)(u(T\pv u_0)\big)\big]-\E\big[\big(P_{t-T}\phi)(u(T\pv\tilde u_0)\big)\big]\\ 
&=\int_{\CC^{-\kappa_0}}\big((P_t\Phi)(u_1)-(P_t\Phi)(\tilde u_1)\big)\P_T^{u_0}(\de u_1)\otimes\P_T^{\tilde u_0}(\de\tilde u_1)
\end{align}
for $t\ge T$ and $\Phi$ any borelian function on $\CC^{-\kappa_0}$ where $\P_T^{u_0}(\de u_1)$ denotes the law of a couples of independent random variables distributed as $u(T\pv u_0)$. For $u_0,\tilde u_0\in\mathcal{A}_R$ and $A=B_{\CC^{-\kappa}}(0,r)\times B_{\CC^{-\kappa}}(0,r)$, we get
\begin{align}
\|P_t^*\delta_{u_0}-P_t^*\delta_{\tilde u_0}\|_{\text{TV}}&\le(\P_T^{u_0}\otimes\P_T^{\tilde u_0})(A^C)+(1-a)(\P_T^{u_0}\otimes\P_T^{\tilde u_0})(A)\\
&\le 1-a(\P_T^{u_0}\otimes\P_T^{\tilde u_0})(A)\\ 
&\le 1-ab^2
\end{align}
for $t\ge T$ which is the claimed result. Taking $R>0$ large enough, we have
\begin{equation}
\inf_{u_0\in\CC^{-\kappa_0}}\inf_{t\ge 1}\mathbb{P}(\|u(t\pv u_0)\|_{\CC^{-\kappa}}\le R)>\frac{1}{2}
\end{equation}
using Markov ineguality with the uniform bound from Corollary \ref{BoundMoments} uniform in $t\ge1$ and $u_0\in\CC^{-\kappa_0}$. As above with $B=\mathcal{A}_R\times\mathcal{A}_R$, we get
\begin{align}
\|P_t^*\delta_{u_0}-P_t^*\delta_{\tilde u_0}\|_{\text{TV}}&\le(\P_1^{u_0}\otimes\P_1^{\tilde u_0})(B^C)+(1-ab^2)(\P_1^{u_0}\otimes\P_1^{\tilde u_0})(B)\\
&\le 1-ab^2(\P_1^{u_0}\otimes\P_1^{\tilde u_0})(B)\\ 
&\le 1-\frac{ab^2}{4}
\end{align}
for $t\ge1+T$ using that the two independent solutions belongs to $\A_R$ with probability greater than $\frac{1}{4}$ which completes the proof.
\end{proof}

We finally prove the main result of exponential convergence to a unique invariant measure. In particular, the measure is the Gibbs measure constructed directly in our previous work \cite{EMR24}.

\medskip

\begin{cor}\label{Cor:ExpMixing}
There exists a unique invariant measure $\mu$ on $\CC^{-\kappa}$ and we have
\begin{equation}
\|P_t^*\delta_{u_0}-\mu\|_{\text{TV}}\le(1-\lambda)^{\lfloor\frac{t}{t_\lambda}\rfloor}\|\delta_{u_0}-\mu\|_{\text{TV}}
\end{equation}
for any $u_0\in\CC^{-\kappa}$ and $t\ge t_\lambda$.
\end{cor}

\medskip

\begin{proof}
For any measure $\mu,\tilde\mu$ on $\CC^{-\kappa}$, we have
\begin{equation}
\|P_t^*\mu-P_t^*\tilde\mu\|_{\text{TV}}=\frac{1}{2}\sup_{\|\varphi\|_{L^\infty}\le 1}\int_{\CC^{-\kappa}\times\CC^{-\kappa}}|P_t\varphi(u_0)-P_t\varphi(\tilde u_0)|\P(\de u_0,\de\tilde u_0)
\end{equation}
for any coupling $\P$ of $\mu$ and $\tilde\mu$ and $t\ge t_\lambda$. Proposition \ref{Prop:lowerboundProba} gives
\begin{equation}
\|P_t^*\mu-P_t^*\tilde\mu\|_{\text{TV}}\le(1-\lambda)(1-\P(D))
\end{equation}
with $D=\{(u_0,u_0);u_0\in\CC^{-\kappa}\}$ the diagonal of $\CC^{-\kappa}\times\CC^{-\kappa}$. Since the total variation distance can be characterized as
\begin{equation}
\|\mu-\tilde\mu\|_{\text{TV}}=\inf\big\{1-\P(D)\ ;\ \P\text{ any coupling of }\mu\text{ and }\tilde\mu\big\}
\end{equation}
thus
\begin{equation}
\|P_t^*\mu-P_t^*\tilde\mu\|_{\text{TV}}\le(1-\lambda)\|\mu-\tilde\mu\|_{\text{TV}}
\end{equation}
which implies that the invariant measure is unique as two distinct invariant measure are singular, see \cite{DZ96}. Denoting as $\mu$ the unique invariant measure, we get
\begin{equation}
\|P_t^*\delta_{u_0}-\mu\|_{\text{TV}}\le(1-\lambda)\|P_{t-t_\lambda}^*\delta_{u_0}-\mu\|_{\text{TV}}
\end{equation}
for any $u_0\in\CC^{-\kappa}$ which completes the proof.
\end{proof}

\appendix

\section{The Anderson Hamiltonian}\label{App:the_anderson_hamiltonian}

We gather here results on $\H$ needed in this work, we refer to \cite{BDM25,M21} and references therein for details, see also Section $4$ in our previous work \cite{EMR24}. The condition for exponents to belongs to $(-1,1)$ is due to the roughness of the noise $\xi$. First, we need the continuity of the Anderson heat semigroup in Hölder spaces.

\medskip

\begin{prop}\label{Prop:ContinuityH}
For any $\alpha\in(-1,1)$, the semigroup $e^{-t\H}$ is continuous from $\CC^\alpha$ to itself.
\end{prop}

\medskip

This is Proposition $4.1$ in \cite{EMR24}, see also Proposition $4.6$ in \cite{BDM25}.

\medskip

\begin{prop}\label{Prop:SchauderH}
For any $\alpha,\beta\in(-1,1)$ such that $\alpha<\beta$, we have
\begin{equation}
\|e^{-t\H}u_0\|_{\CC^\beta}\lesssim t^{-\frac{\alpha-\beta}{2}}\|u_0\|_{\CC^\alpha}
\end{equation}
for any $t\in(0,1)$.
\end{prop}

\medskip

This follows from Proposition $4.12$ in \cite{BDM25}.

\medskip

\begin{prop}\label{Prop:GaussianBound}
There exists a constant $c>0$ such that
\begin{equation}
K_t(x,y)\lesssim\frac{1}{t}e^{-c\frac{|x-y|^2}{t}}
\end{equation}
for any $t\in(0,1)$ and $x,y\in\T^2$.
\end{prop}

\medskip

In the construction of the Anderson Hamiltonian, a crucial tool is the $\Gamma$ map defined as the inverse of the map
\begin{equation}
\Phi(v)=v-v\pl X_{>n}
\end{equation}
where $X=\Delta^{-1}\xi\in\CC^{1-\kappa}$ for any $\kappa>0$. This map is not invertible in general, one can introduce a truncation $X_{>n}$ as done in \cite{GUZ20} such that the map
\begin{equation}
\Phi_n(v)=v-v\pl X_{>n}
\end{equation}
is a small perturbation of the identity thus invertible for $n$ large enough depending on $\|X\|_{\CC^{1-\kappa}}$, see Section $2.1.1$ in \cite{GUZ20} or Section $2.2$ in \cite{M21}. In particular, $\Phi_n$ and $\Gamma_n$ converges to the identity as $n$ goes to infinity hence continuity estimates are uniform with respect to $n$ large enough.

\medskip

\begin{lem}\label{Lem:truncation}
There exists $n_0=n_0(\xi)>0$ such that the maps $\Phi_n$ and $\Gamma_n$ are continuous on $L^p$ for $p\ge1$ uniformly over $n\ge n_0$.
\end{lem}

\medskip

Recall that
\begin{equation}
\Gamma_n^{-1}\H\Gamma_nw= -\Delta w + \xi\ple w + R_n(w)
\end{equation}
with $R_n$ depending on the enhanced noise. The following result follows from the construction in \cite{M21}, see also Lemma $4.5$ in \cite{BDM25}. Note that the roughest term $\xi\ple w$ is not included in $R_n$ here thus the exponent $\kappa$ in the result.

\medskip

\begin{lem}\label{Lem:RemainderH}
The remainder term is a bounded operator
$$
R_n:B_{q,\infty}^{\kappa}\to B_{q,\infty}^{-\kappa'}
$$
for arbitrary $\kappa,\kappa'>0$ and $q\in[1,\infty]$.
\end{lem}

\section{Estimates in Besov spaces}{\label{APP:Besov}}

A natural and convenient setting is given by the Besov space $B_{p,q}^\alpha$ which can be defined using the Littlewood-Paley  decomposition, see for example Bahouri, Chemin and Danchin's book\cite{BCD11}. This decomposition can be stated  as follows 
\begin{equation*}
u=\sum_{n\ge0}\Delta_nu
\end{equation*}
with $\Delta_nu=\big(\mathcal{F}^{-1}\IDC_{|\cdot|\simeq 2^n}\mathcal{F}\big)u$, that is the projection of $u$ in frequencies on an annulus of size $2^n$. It is defined by 
\begin{equation*}
\big(\Delta_nu\big)(x)=\int_{\R^d}\chi_n(x-y)u(y)\de y:=2^{d(n-1)}\int_{\R^d}\chi\big(2^{n-1}(x-y)\big)u(y)\de y
\end{equation*} 
with $\chi\in\mathcal{S}(\R^d)$ and $\text{supp}\ \widehat\chi\subset\{\frac{1}{2}\le|z|\le 2\}$ for $n\ge1$ and
\begin{equation*}
\big(\Delta_0u\big)(x):=\int_{\R^d}\chi_0(x-y)u(y)\de y
\end{equation*}
with $\chi_0\in\mathcal{S}(\R^d)$ and $\text{supp}\ \widehat\chi_0\subset\{|z|\le 1\}$. We also denote $K=\widehat{\chi}$ such that $\Delta_nu=\big(\mathcal{F}^{-1}K(2^n\cdot)\mathcal{F}\big)u$. Then the Besov space $B_{p,q}^\alpha$ are distributions such that
\begin{equation*}
\|u\|_{B_{p,q}^\alpha}:=\Big(\ \sum_{n\ge0}2^{\alpha qn}\|\Delta_nu\|_{L^p(\T^d)}^q\ \Big)^{\frac{1}{q}}<\infty
\end{equation*}
where the sum is understood to be a supremum if $q=\infty$. The particular case $p=q=2$ corresponds to the Sobolev space $B_{2,2}^\alpha=H^\alpha$ and for $p=q=\infty$ with $\alpha\in\R_+\backslash\N$, one gets the usual Hölder spaces $B_{\infty,\infty}^\alpha=\CC^\alpha$. The following lemma gathers important properties such as embedding or product rule.

\medskip

\begin{lem}\label{LEM:Besov}
The Besov spaces $B^s_{p,q}$ defined above satisfy the following properties.\\
\textup{(i)} For any $s\in\R$ we have $B^s_{2,2}=H^s$, and more generally for any $2\le p <\infty$ and $\eps>0$ we have
\begin{align*}
\|f\|_{B^s_{p,\infty}}\les \|f\|_{W^{s,p}}\les \|f\|_{B^s_{p,2}} \les \|f\|_{B^{s+\eps}_{p,\infty}}.
\end{align*}
\textup{(ii)} Let $s\in\R$ and $1\leq p_1\leq p_2\leq \infty$ and $q_1,q_2\in [1,\infty]$. Set $r=s-2\big(\frac{1}{p_1}-\frac{1}{p_2}\big)$, then for any $f\in B^s_{p_1,q}$, we have the embedding
\begin{align*}
\|f\|_{B^{r}_{p_2,q}}\les \|f\|_{B^s_{p_1,q}}.
\end{align*}
\textup{(iii)} Let $\alpha,\beta\in \R$ such that $\alpha+\beta>0$ and $p_1,p_2,q_1,q_2\in [1,\infty]$ with
\begin{align*}
\frac1p = \frac{1}{p_1}+\frac{1}{p_2}\qquad\text{ and }\qquad\frac1q=\frac{1}{q_1}+\frac{1}{q_2}.
\end{align*}
Then for any $f\in B^{\alpha}_{p_1,q_1}$ and $g\in B^{\beta}_{p_2,q_2}$, we have $fg\in B^{\alpha\wedge \beta}_{p,q}$ and the following holds.
\begin{itemize}
    \item If $\alpha\wedge \beta<0$, then
    \begin{align*}
    \|fg\|_{B^{\alpha\wedge\beta}_{p,q}}\les \|f\|_{B^{\alpha}_{p_1,q_1}}\|g\|_{B^{\beta}_{p_2,q_2}}.
    \end{align*}
    \item If $\alpha\wedge \beta >0$, then
    \begin{align*}
    \|fg\|_{B^{\alpha\wedge\beta}_{p,q}}\les \|f\|_{B^{\alpha\wedge\beta}_{p_1,q_1}}\|g\|_{B^{\alpha\wedge\beta}_{p_2,q_2}}.
    \end{align*}
\end{itemize}
\textup{(iv)} Let $\alpha\in\R$ and $p,q\in[1,\infty]$. Then we have
\begin{equation}
\Big|\int_{\T^d}f(x)g(x)\de x\Big|\lesssim\|f\|_{B_{p,q}^\alpha}\|g\|_{B_{p',q'}^{-\alpha}}
\end{equation}
where $p',q'$ are the respective conjugate exponent of $p,q$.
\end{lem}

\medskip

From \cite{MW17a}, we gather the two following most useful estimates. The following lemma is a useful interpolation estimate.

\medskip

\begin{lem}{\label{LEM:InterpolationBesov}}
Let $\alpha_0,\alpha_1\in\R$, $p_0,p_1,q_0,q_1\in[1,+\infty]$ and set
$$
\alpha = \theta\alpha_1 + (1-\theta)\alpha_2\qquad \frac{1}{p}=\frac{\theta}{p_1}+\frac{1-\theta}{p_2}\qquad \frac{1}{q}=\frac{\theta}{q_1}+\frac{1-\theta}{q_2}
$$
for $\theta\in[0,1]$. Then
$$
\|f\|_{B^\alpha_{p,q}}\lesssim \|f\|_{B^{\alpha_1}_{p_1,q_1}}^{\theta}\|f\|_{B^{\alpha_2}_{p_2,q_2}}^{1-\theta}.
$$
\end{lem}

\medskip

This product rule in mixed Besov and Lebesgue spaces is useful for energy estimate.

\medskip

\begin{lem}{\label{LEM:PowerProductBesov}}
Let $\alpha>0$, $r\in\N$ and $p_1,p_2,q\in[1,+\infty]$. Set
$$
\frac{1}{p}=\frac{1}{p_1}+\frac{1}{p_2}
$$
then
$$
\|f^{r+1}\|_{B^\alpha_{p,q}}\lesssim \|f^r\|_{L^{p_1}}\|f\|_{B^\alpha_{p_2,q}}.
$$
\end{lem}

\medskip

We will also make use of the following fractional Leibniz rule in $L^2$-based Sobolev spaces.

\medskip

\begin{lem}\label{LEM:prodHugo}
Let $r,s<1$ with $r+s>0$, and $t=r+s-1$. Then there exists $C>0$ such that for any $u\in H^s$ and $v\in H^r$, we have the fractional Leibniz rule
$$
\|uv\|_{H^t}\le C\|u\|_{H^s}\|v\|_{H^r}.
$$
\end{lem}

\medskip

\begin{proof}
The corresponding estimate in $\R^d$ is proved in \cite[Corollary 2.1]{T01}. Then the one on the torus follows from the one in $\R^d$ through the use of a finite partition of unity.
\end{proof}


\vspace{2cm}

\noindent \textcolor{gray}{$\bullet$} H. Eulry -- École Normale Supérieure de Lyon, UMPA, UMR CNRS-ENSL 5669, 46,
allée d’Italie, 69364-Lyon Cedex 07, France\\
{\it E-mail}: hugo.eulry@ens-lyon.fr

\smallskip

\noindent \textcolor{gray}{$\bullet$} A. Mouzard -- Modal’X - UMR CNRS 9023, Université Paris Nanterre, 92000 Nanterre, France.\\
{\it E-mail}: antoine.mouzard@math.cnrs.fr


\begin{thebibliography}{10}

\bibitem{BCD11}
{\sc H.~Bahouri, J.-Y. Chemin, and R.~Danchin}, {\em Fourier {{Analysis}} and
  {{Nonlinear Partial Differential Equations}}}, vol.~343 of Grundlehren Der
  Mathematischen {{Wissenschaften}}, Springer, Berlin, Heidelberg, 2011.

\bibitem{BDFT23}
{\sc I.~Bailleul, N.~V. Dang, L.~Ferdinand, and T.~D. T{\^o}}, {\em $\phi^4_3$
  measures on compact {{Riemannian}} $3$-manifolds}, Oct. 2023.

\bibitem{BDM25}
{\sc I.~Bailleul, N.~V. Dang, and A.~Mouzard}, {\em Analysis of the anderson
  operator}, 2025.

\bibitem{CC18}
{\sc R.~Catellier and K.~Chouk}, {\em {Paracontrolled distributions and the
  3-dimensional stochastic quantization equation}}, The Annals of Probability,
  46 (2018), pp.~2621 -- 2679.

\bibitem{CF18}
{\sc K.~Chouk and P.~K. Friz}, {\em Support theorem for a singular {{SPDE}}:
  {{The}} case of {{gPAM}}}, Annales de l'Institut Henri Poincar{\'e},
  Probabilit{\'e}s et Statistiques, 54 (2018), pp.~202--219.

\bibitem{DD03}
{\sc G.~Da~Prato and A.~Debussche}, {\em Strong solutions to the stochastic
  quantization equations}, The Annals of Probability, 31 (2003),
  pp.~1900--1916.

\bibitem{DD20}
\leavevmode\vrule height 2pt depth -1.6pt width 23pt, {\em Gradient estimates
  and maximal dissipativity for the {{Kolmogorov}} operator in $\phi
  ^{4}_{2}$}, Electronic Communications in Probability, 25 (2020), pp.~1--16.

\bibitem{DZ96}
{\sc G.~Da~Prato and J.~Zabczyk}, {\em Ergodicity for {{Infinite Dimensional
  Systems}}}, London {{Mathematical Society Lecture Note Series}}, Cambridge
  University Press, Cambridge, 1996.

\bibitem{DZ97}
{\sc G.~Da~Prato and J.~Zabczyk}, {\em Differentiability of the {{Feynman-Kac}}
  semigroup and a control application},  (1997).

\bibitem{DZ14}
\leavevmode\vrule height 2pt depth -1.6pt width 23pt, {\em Stochastic
  {{Equations}} in {{Infinite Dimensions}}}, Encyclopedia of {{Mathematics}}
  and Its {{Applications}}, Cambridge University Press, Cambridge, 2~ed., 2014.

\bibitem{EMR24}
{\sc H.~Eulry, A.~Mouzard, and T.~Robert}, {\em Anderson stochastic
  quantization equation}, Jan. 2024.

\bibitem{GJ81}
{\sc J.~Glimm and A.~Jaffe}, {\em Quantum physics. {A} functional integral
  point of view}.
\newblock New {York} - {Heidelberg} - {Berlin}: {Springer}-{Verlag}. {XX}, 417
  p., 43 ill. \$ 26.40 (1981)., 1981.

\bibitem{GH21}
{\sc M.~Gubinelli and M.~Hofmanov{\'a}}, {\em A {PDE} construction of the
  {Euclidean} {{\(\Phi^4_3\)}} quantum field theory}, Commun. Math. Phys., 384
  (2021), pp.~1--75.

\bibitem{GIP15}
{\sc M.~Gubinelli, P.~Imkeller, and N.~Perkowski}, {\em {{PARACONTROLLED
  DISTRIBUTIONS AND SINGULAR PDES}}}, Forum of Mathematics, Pi, 3 (2015),
  p.~e6.

\bibitem{GUZ20}
{\sc M.~Gubinelli, B.~Ugurcan, and I.~Zachhuber}, {\em Semilinear evolution
  equations for the {{Anderson Hamiltonian}} in two and three dimensions},
  Stochastics and Partial Differential Equations: Analysis and Computations, 8
  (2020), pp.~82--149.

\bibitem{JP21}
{\sc A.~Jagannath and N.~Perkowski}, {\em A simple construction of the
  dynamical $\phi^4_3$ model}, Aug. 2021.

\bibitem{MW17a}
{\sc J.-C. Mourrat and H.~Weber}, {\em The {{Dynamic}} ${\Phi^4_3}$ {{Model
  Comes Down}} from {{Infinity}}}, Communications in Mathematical Physics, 356
  (2017), pp.~673--753.

\bibitem{MW17}
\leavevmode\vrule height 2pt depth -1.6pt width 23pt, {\em Global
  well-posedness of the dynamic $\phi^{4}$ model in the plane}, The Annals of
  Probability, 45 (2017), pp.~2398--2476.

\bibitem{M21}
{\sc A.~Mouzard}, {\em Weyl law for the {{Anderson Hamiltonian}} on a
  two-dimensional manifold}, July 2021.

\bibitem{PW81}
{\sc G.~Parisi and Y.-s. Wu}, {\em Perturbation {{Theory Without Gauge
  Fixing}}}, Sci. Sin., 24 (1981), p.~483.

\bibitem{S24}
{\sc D.~W. Stroock}, {\em Probability theory. {An} analytic view}, Cambridge:
  Cambridge University Press, 3rd edition~ed., 2024.

\bibitem{SV72}
{\sc D.~W. Stroock and S.~R.~S. Varadhan}, {\em On the support of diffusion
  processes with applications to the strong maximum principle}.
\newblock Proc. 6th {Berkeley} {Sympos}. math. {Statist}. {Probab}., {Univ}.
  {Calif}. 1970, 3, 333-359 (1972)., 1972.

\bibitem{T01}
{\sc J.~Tamba{\v c}a}, {\em Estimates of the {{Sobolev Norm}} of a {{Product}}
  of {{Two Functions}}}, Journal of Mathematical Analysis and Applications, 255
  (2001), pp.~137--146.

\bibitem{TW18}
{\sc P.~Tsatsoulis and H.~Weber}, {\em Spectral gap for the stochastic
  quantization equation on the 2-dimensional torus}, Annales de l'Institut
  Henri Poincar{\'e}, Probabilit{\'e}s et Statistiques, 54 (2018),
  pp.~1204--1249.

\end{thebibliography}
\end{document}